\documentclass[final,a4paper]{article}

\usepackage{ucs}
\usepackage[utf8x]{inputenc}
\usepackage[english]{babel}
\usepackage[spaces,hyphens]{url}

\usepackage{booktabs}
\usepackage{color}
\usepackage{float}%
  \floatstyle{ruled}%
  \newfloat{Algorithm}{tp}{lsa}%
\usepackage{graphicx}
\usepackage[procnames]{listings}
  \lstMakeShortInline{\"}%
  \newcommand{\fn}[1]{{\slshape\bfseries #1}}%
  \newcommand{\n}[1]{\begingroup\def\x{\csname label\endcsname{l:#1}}\x\endgroup}%
  \newcommand{\nr}[1]{\begingroup\def\x{\csname ref\endcsname{l:#1}}\x\endgroup}%
  \newcommand\q[1]{{\slshape #1}}%
  \newcommand{\var}[1]{$\langle$\textit{#1}$\rangle$}%
  \DeclareMathAlphabet{\mathitbf}{OML}{cmm}{b}{it}
  
\usepackage{mdwlist}
\usepackage[newitem,newenum]{paralist}
  \setdefaultenum{\slshape (1)}{\slshape (a)}{\slshape (i)}{A.}
  \setdefaultitem{\guillemotright}{\guilsinglright}{$-$}{*}
\usepackage{my/miscmath}
  \usepackage{amsthm}
  
  \newtheorem{theorem}{Theorem}[section]
  \newtheorem{lemma}{Lemma}[section]
  \theoremstyle{definition}
  \newtheorem{definition}{Definition}[section]
\usepackage{relsize}
\usepackage{rotating}
\usepackage{my/xy}

\usepackage{hyperref}

\usepackage{authblk}
\PrerenderUnicode{üä}
\title{Fatgraph Algorithms and the Homology of the Kontsevich Complex}
\author{Riccardo Murri}
\affil{%
    Grid Computing Competence Center,\\
    Organisch-Chemisches Institut,\\
    Universität Zürich,\\
    Winterthurerstrasse 190, CH-8006 Zürich,\\
    Switzerland.\\
    e-mail: \url{riccardo.murri@gmail.com}
}%
\date{Feb.~14, 2012}

\usepackage{relsize}
\newcommand\ifmonospace{%
  \ifdim\fontdimen3\font=0pt%
}
\newcommand\Cpp{%
  \ifmonospace%
  C++%
  \else%
  C\kern-.0667em\raise.30ex\hbox{\relsize{-2}{++}}%
  \fi%
  \spacefactor1000%
}

\begin{document}
\maketitle

\begin{abstract}
  Fatgraphs are multigraphs enriched with a cyclic order of the edges
  incident to a vertex.  This paper presents algorithms to:
  \begin{inparaenum}
  \item generate the set $\R_{g,n}$ of fatgraphs, given the genus $g$
    and the number of boundary cycles $n$;
  \item compute automorphisms of any given fatgraph;
  \item compute the homology of the fatgraph complex $\R_{g,n}$.
  \end{inparaenum}
  The algorithms are suitable for effective computer implementation.

  In particular, this allows us to compute the rational homology of
  the moduli space of Riemann surfaces with marked points. We thus
  compute the Betti numbers of $M_{g,n}$ with $(2g + n) \leq 6$,
  corroborating known results.
\end{abstract}

\section{Introduction}
\label{sec:intro}

This paper deals with algorithms for the enumeration of
fatgraphs and their automorphisms, and the computation of the homology
of the complex formed by fatgraphs of a given genus $g$ and number of
boundary components $n$.  

\begin{figure}
  \centering
  \includegraphics[width=\linewidth]{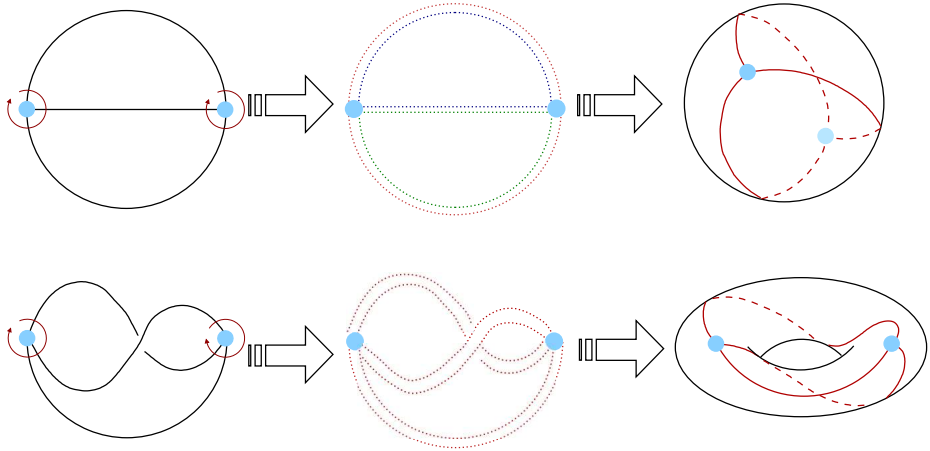}
  \caption{%
    Thickening of a fatgraph into a Riemann surface. \emph{Left
      column:} Starting fatgraph: the cyclic order at the vertices is
    given by the orientation of the ambient euclidean
    plane. \emph{Middle column:} Thickening of the fatgraph by gluing
    topological disks along the boundary components. The border of a
    cells is drawn as a dotted line; each topological disk has been
    given a different color. \emph{Right column:} The resulting
    Riemann surface with the embedded graph.  Note that the two
    starting graphs would be isomorphic when considered as ordinary
    multigraphs; they are distinguished by the additional cyclic
    structure at the vertices.}
  \label{fig:fattening}
\end{figure}

A fatgraph\footnote{%
  Fatgraphs have appeared independently in many different areas of
  mathematics: several equivalent definitions are known, with names
  such as ``ribbon graphs'', ``cyclic graphs'', ``maps'', ``dessins
  d'enfants'', ``rotation systems''.  See \cite{lando-zvonkin} for a
  comprehensive survey.  }\space%
is a multigraph enriched with the assignment, at each vertex $v$, of a
cyclic order of the edges incident to $v$.  Such graphs can be
``fattened'' into a smooth punctured oriented surface, by gluing
polygons along the edges in such a way that two adjacent edges on the
polygon boundary are consecutive in the cyclic order at the common
endpoint (see Figure~\ref{fig:fattening}); an additional assignment of
a length for each edge allows to define a conformal structure on the
surface.  The resulting Riemann surface is naturally marked, by
choosing the marking points to be the centers of the polygons. There
is thus a functorial correspondence between fatgraphs and marked
Riemann surfaces; a fatgraph $G$ is said to have genus $g$ and $n$
boundary components if it corresponds to a punctured Riemann surface
$S \in \M_{g,n}$.

In the papers \cite{kontsevich;1993} and \cite{kontsevich;feynman},
M.~Kontsevich introduced ``Graph Homology'' complexes that relate the
stable homology groups of certain infinite-dimensional Lie algebras to
various other topological objects.  In particular, the ``associative
operad'' variant of this construction results in a chain complex whose
homology is isomorphic to the (co)homology of the moduli space of
smooth Riemann surfaces $\M_{g,n}$: the graded module underlying the
complex is freely generated by the set $\R_{g,n}$ of fatgraphs of
genus $g$ and number of boundary components $n$, endowed with the
differential defined by edge contraction.  

The needed definitions and theorems about fatgraphs and their homology
complex are briefly recalled in Section~\ref{sec:mgn}; the interested
reader is referred to \cite{mondello:arXiv:0705.1792v1} and
\cite{lando-zvonkin} for proofs and context.

The bulk of this paper is concerned with finding an effectively
computable representation of fatgraphs (see
Section~\ref{sec:Fatgraph}), and presenting algorithms to:
\begin{enumerate}
\item compute automorphisms of any given fatgraph
  (Section~\ref{sec:isomorphism});
\item generate the set $\R_{g,n}$ of fatgraphs, given the genus $g$
  and number of boundary components $n$
  (Section~\ref{sec:generation});
\item compute the homology of the fatgraph complex $\R_{g,n}$
  (Section~\ref{sec:homology}).
\end{enumerate}
Note that, in contrast with other computational approaches to
fatgraphs (e.g., \cite{penner+knudsen+wiuf+andersen:2009}) which draw
on the combinatorial definition of a fatgraph, our computer model of
fatgraphs is directly inspired by the topological definition, and the
algorithm for enumerating elements of $\R_{g,n}$ is likewise backed by
a topological procedure.

Theorem~\ref{thm:fatgraph-homology} provides an effective way to
compute the (co)homology of $\M_{g,n}$.  The Betti numbers of
$\M_{g,n}$ can be computed from the knowledge of the dimension of
chain spaces $W_p$ of the fatgraph complex and the ranks of boundary
operators $D_p$; this computation can be accomplished in the following
stages:
\begin{enumerate}[I.]
\item Generate the basis set of $W_*$; by definition, the basis set is
  the set $\R_{g,n}$ of \emph{oriented} fatgraphs that correspond to
  surfaces in $\M_{g,n}$.
\item Work out the differential $D\colon W_* \to W_*$ as matrices
  $\+D^{(p)}$ mapping coordinates in the fatgraph basis of $W_p$ into
  coordinates relative to the fatgraph basis of $W_{p-1}$.
\item Compute the ranks of the matrices $\+D^{(p)}$.
\end{enumerate}

Stage~I needs just the pair $g,n$ as input; its output is the set of
orientable marked fatgraphs belonging in $\R_{g,n}$.  By definition,
marked fatgraphs are decorated abstract fatgraphs, and the decoration
is a simple combinatorial datum (namely, a bijection of the set of
boundary cycles with the set $\{ 1, \ldots, n \}$): therefore, the
problem can be reduced to enumerating abstract fatgraphs.  With a
recursive algorithm, one can construct \emph{trivalent}
$\M_{g,n}$-fatgraphs from trivalent graphs in $\M_{g-1,n}$ and
$\M_{g-1,n+1}$.  All other graphs in $\M_{g,n}$ are obtained by
contraction of non-loop edges.

The differential $D$ has a simple geometrical definition: $D(G)$ is a
sum of graphs $G'$, each gotten by contracting a non-loop edge of $G$. A
simple implementation of Stage~II would just compare each contraction
of a graph with $p$ edges with any graph with $p-1$ edges, and score a
$\pm 1$ (depending on the orientation) in the corresponding entry of
the matrix $\+D^{(p)}$.  However, this algorithm has quadratic complexity,
and the large number of graphs involved makes it very inefficient
already for $\M_{0,5}$.  The simple observation that contraction of
edges is defined on the topological fatgraph underlying a marked
fatgraph allows us to apply the naive algorithm to topological
fatgraphs only, which cuts complexity down by a factor~$O((n!)^2)$.  The
resulting matrix is then extended to marked fatgraphs by the action
of graph automorphism groups on the markings of boundary cycles.
This is the variant detailed in Section~\ref{sec:homology}.

Stage~III is conceptually the simplest: by elementary linear algebra,
the Betti numbers can be computed from the rank of matrices $\+D^{(p)}$ and
the dimension of their domain space.  The computational problem of
determining the rank of a matrix has been extensively studied; it
should be noted, however, that this step can actually be the most
computationally burdening.

It is worth mentioning that V.~Godin \cite{godin:homology} introduced
a slightly different fatgraph complex, which computes
the \emph{integral} (co)homology of $\M_{g,n}$; possible adaptation of
the algorithms to this complex and an outlook on the expected problems
is given in Section~\ref{sec:conclusions}.

An effective implementation (using the Python programming language
\cite{python:website}) of the algorithms presented here is available
at \url{http://code.google.com/p/fatghol}.  It has so far been used to
compute the Betti numbers of~$\M_{g,n}$ for~$(2g + n) \leq 6$.

Results are summarized in Table~\ref{tab:betti}: the values coincide
with results already published in the literature.  References are
given in the closing Section~\ref{sec:conclusions}, together with a
discussion on the implementation performance and possible future
directions for improving and extending the algorithms.
\begin{table}[hbt]
  \centering
  \begin{tabular}{cccccccccccccc}
    \toprule
           ~ & $b_0$ & $b_1$ & $b_2$ & $b_3$ & $b_4$ & $b_5$ & $b_6$ & $b_7$ & $b_8$ & $b_9$ & $b_{10}$ & $b_{11}$ & $b_{12}$ \\
    \midrule
    $\M_{0,3}$ & 1    &      &       &       &       &       &&&&&&& \\
    $\M_{0,4}$ & 1    & 2    &       &       &       &       &&&&&&& \\
    $\M_{0,5}$ & 1    & 5    & 6     &       &       &       &&&&&&& \\
    $\M_{0,6}$ & 1    & 9    & 26    & 24    &       &       &&&&&&& \\
    $\M_{1,1}$ & 1    &      &       &       &       &       &&&&&&& \\
    $\M_{1,2}$ & 1    &      &       &       &       &       &&&&&&& \\
    $\M_{1,3}$ & 1    &      &       & 1     &       &       &&&&&&& \\
    $\M_{1,4}$ & 1    &      &       & 4     & 3     &       &&&&&&& \\
    $\M_{2,1}$ & 1    &      & 1     &       &       &       &&&&&&& \\
    $\M_{2,2}$ & 1    &      & 2     &       &       & 1     &&&&&&& \\
    \bottomrule
  \end{tabular}
  \caption{\label{tab:betti}%
    Betti numbers of $\M_{g,n}$ for $2g+n\leq6$.  For readability, null
    values have been omitted and the corresponding entry left blank.
    See Section~\ref{sec:conclusions} for a discussion of these results.
  }
\end{table}

\subsection{Notation}
\label{sec:notation}

Algorithms are listed in pseudo-code reminiscent of the Python
language syntax (see \cite{python:reference27}); comments in the code
listings are printed in \emph{italics} font.  The word ``object'' is
used to denote an heterogeneous composite type in commentaries to the
code listings: for our purposes, an \emph{object} is just a tuple
`"($a_1$, $a_2$, $\ldots$, $a_N$)"', where each of the slots "$a_i$"
can be independently assigned a value;\footnote{%
  This is the definition of what is usually called a ``record'' in
  Computer Science literature, and lacks important features of what is
  generally meant by ``objects'' in a programming context.  However,
  the Python programming language only provides objects (i.e., records
  are implemented as objects with no methods), and our algorithm
  implementation relies on object-oriented programming features. We
  have thus decided to keep our choice of words closer to the actual
  code.%
} we write $X.a_i$ to denote the slot $a_i$ of object $X$.  Object
slots are mutable, i.e., they can be assigned different values over
the course of time.  Appendix~\ref{sec:pseudo-code} gives a complete
recap of the notation used and the properties assumed of syntax, data
structures, and operators.

A great deal of this paper is concerned with finding
computationally-effective representations of topological objects; in
general, we use boldface letters to denote the computer analog of a
mathematical object.  For instance, the letter $G$ always denotes a
fatgraph, and $\+G$ its corresponding computer representation as a
"Fatgraph" object.

Finally, if $A$ is a category of which $X$, $Y$ are objects, we use
Eilenberg's notation $A(X,Y)$ for the $\Hom$-set, instead of the more
verbose $\Hom_A(X,Y)$.

\section{Fatgraphs and marked Riemann surfaces}
\label{sec:mgn}

This section recaps the main definitions and properties of fatgraphs
and the relation of the fatgraph complex to the cohomology of
$\M_{g,n}$.  These results are well-known: a clear and comprehensive
account is given by G.~Mondello in \cite{mondello:arXiv:0705.1792v1};
the book by Lando and Zvonkin \cite{lando-zvonkin} provides a broad
survey of the applications of fatgraphs and an introduction accessible
to readers without a background in Algebraic Geometry.

``Fatgraphs'' take their name from being usually depicted as graphs
with thin bands as edges, instead of 1-dimensional lines; they have
also been called ``ribbon graphs'' in algebraic geometry literature.
Here, the two names will be used interchangeably.
\begin{definition}[Geometric definition of fatgraphs]
  \label{dfn:fatgraphs-geom}
  A fatgraph is a finite CW-complex of pure dimension 1, together
  with an assignment, for each vertex $v$, of a cyclic ordering of the
  edges incident at $v$.

  A morphism of fatgraphs is a cellular map $f\colon G\to G'$ such that, for
  each vertex $v$ of $G'$, the preimage $f^{-1}(V)$ of a small
  neighborhood %
  $V$ of $v$ is a small neighborhood of a tree in $G$ (i.e.,
  $f^{-1}(V)$ is a contractible connected graph).
\end{definition}
Unless otherwise specified, we assume that all vertices of a
fatgraph have valence at least 3.

If $G$ is a fatgraph, denote $\Vertices{G}$, $\Edges{G}$ and
$\Legs{G}$ the sets of vertices, unoriented edges and oriented edges
(equivalently called ``legs'' or ``half-edges'').

Let $G$ be a fatgraph, and $G'$ be the CW-complex obtained by
contracting an edge $\alpha \in \Edges{G}$ to a point.  If $\alpha$
connects two \emph{distinct} vertices (i.e.,  $\alpha$ is not a
loop) then $G'$ inherits a fatgraph structure from $G$: if $(\alpha <
\alpha_1 < \ldots < \alpha_k < \alpha)$ and $(\alpha < \alpha'_1 <
... < \alpha'_h < \alpha)$ are the cyclic orders at endpoints of
$\alpha$, then the vertex formed by collapsing $\alpha$ is endowed
with the cyclic order $(\alpha_1 < \ldots < \alpha_k < \alpha'_1 <
\ldots < \alpha'_h)$.  The graph~$G'$ is said to be obtained from~$G$
by contraction of~$\alpha$.

Contraction morphisms play a major role in manipulation of ribbon
graphs.
\begin{lemma}\label{lemma:contraction1}
  Any morphism of fatgraphs is a composition of isomorphisms and
  contractions of non-loop edges.
\end{lemma}
We can thus define functors $\Vertices{-}$, $\Edges{-}$ and $\Legs{-}$
that send morphisms of graphs to maps of their set of vertices,
(unoriented) edges, and oriented edges.

The following combinatorial description of a fatgraph will also be
needed:
\begin{definition}[Combinatorial definition of fatgraph]
  \label{dfn:fatgraphs-combi}
  A fatgraph is a $4$-tuple $(L, \sigma_0, \sigma_1, \sigma_2)$
  comprised of a finite set $L$, together with bijective
  maps $\sigma_0, \sigma_1, \sigma_2\colon L \to L$ such that:
  \begin{itemize}
  \item $\sigma_1$ is a fixed-point free involution: $\sigma_1^2 = \idmap$, and
  \item $\sigma_0 \circ \sigma_2 = \sigma_1$.
  \end{itemize}
\end{definition}
\begin{lemma}
  \label{lemma:comb-fatgraph}
  Definitions~\ref{dfn:fatgraphs-geom} and~\ref{dfn:fatgraphs-combi}
  are equivalent.
\end{lemma}

Any two of the maps $\sigma_0, \sigma_1, \sigma_2$ determine the
third, by means of the defining relation $\sigma_0 \circ \sigma_2 =
\sigma_1$; therefore, to give a ribbon graph it is sufficient to
specify only two out of three maps.

In the combinatorial description, $\Vertices{G}$ is the set $L_0$ of
orbits of $\sigma_0$, $\Edges{G}$ is the set $L_1$ of orbits of
$\sigma_1$, and $\Legs{G}$ is plainly the set $L$.

There is a functorial construction to build a topological surface
$S(G)$ from a fatgraph $G$; this is usually referred to as
``thickening'' or ``fattening'' in the literature.
\begin{lemma}
  \label{lemma:S-functor}
  There exists a functor $S$ that associates to every fatgraph $G$ a
  punctured Riemann surface $S(G)$, and to every morphism $f\colon G \to
  G'$ a continuous map $S(f)\colon S(G) \to S(G')$.
\end{lemma}

Denote by $\Holes{G}$ the set $L/\sigma_2$ of orbits of $\sigma_2$: in the
topological description, its elements are the support of 1-cycles in
$H^1(G)$ that correspond under a retraction to small loops around
the punctures in $S(G)$; they are called ``boundary cycles'' of
$G$.

The assignment $G \mapsto \Holes{G}$ extends to a functor $\Holes{-}$;
by Lemma~\ref{lemma:contraction1}, for any $f\colon G_1 \to G_2$ the
map $\Holes{f} \colon \Holes{G_1} \to \Holes{G_2}$ is a bijection.

The correspondence between fatgraphs and Riemann surfaces allows
us to give the following.
\begin{definition}\label{dfn:rg-genus}
  The number of boundary cycles of a graph $G$ is given by $n =
  \card{\Holes{G}}$, and is equal to the puncture number of the
  Riemann surface $S(G)$.

  If $S(G)$ has genus $g$ and $n$ boundary cycles, then:
  \begin{equation}
    \label{eq:euler-characteristics-of-fatgraphs}
    \chi(G) := \chi(S(G)) = 2 - 2g - n = 2 - 2g - \card{\Holes{G}},
  \end{equation}
  so we can define, for any fatgraph $G$, the \emph{genus} $g$, as
  given by the relation above.
\end{definition}
\begin{lemma}\label{lemma:contraction2}
  If $G'$ is obtained from $G$ by contraction of a non-loop edge, then
  $G$ and $G'$ share the same genus and number of boundary cycles.
\end{lemma}

\begin{definition}
  A \emph{marked fatgraph} is a fatgraph $G$ endowed with a
  bijection $\nu: \Holes{G} \to \{1,\ldots, n\}$.  The map $\nu$ is
  said to be the ``marking'' on $G$.

  A morphism $f\colon G_1 \to G_2$ of marked fatgraphs must preserve the
  marking of boundary cycles:
  \begin{equation*}
    \xymatrix{%
      \Holes{G_1} &  & \Holes{G_2}
      \\
      & \{ 1, \ldots, n \} &
      \ar^{f} "1,1";"1,3"
      \ar_{\nu_1} "1,1";"2,2"
      \ar^{\nu_2} "1,3";"2,2"
    }%
  \end{equation*}
\end{definition}
By a slight abuse of language, we shall usually omit mention of the
marking map~$\nu$ and just speak of ``the marked fatgraph $G$''.

\subsection{Moduli spaces of marked Riemann surfaces}
\label{sec:moduli-spaces}

Fix integers $g\geq 0$, $n>0$ such that $2 -2g - n < 0$. Let $S$ be a
smooth closed oriented surface of genus $g$ and $X = \{ x_1, \ldots,
x_n \}$ a set of points of $S$.
\begin{definition}\label{dfn:teichmuller}
  The Teichmüller space
  \begin{equation*}
    \T_{g,n} := \Conf(S) / \Diff^0(S, n)
  \end{equation*}
  is the quotient of the set of all conformal metrics on $S$ by the
  set of all diffeomorphisms homotopic to the identity and fixing the
  $n$ marked points.

  The mapping class group $\Gamma_{g,n}$ is the group of isotopy
  classes of self-diffeomorphisms that preserve orientation and fix
  marked points:
  \begin{equation*}
    \Gamma_{g,n} := \Diff^+ (S, n) / \Diff^0(S, n).
  \end{equation*}
  
  The topological space $\M_{g,n} := \T_{g,n} / \Gamma_{g,n}$ is the
  moduli space of (smooth) $n$-pointed algebraic curves of genus $g$.
  It parametrizes complex structures on $S$, up to
  diffeomorphisms that:
  \begin{inparaenum}
  \item are homotopic to the identity mapping on $S$,
  \item preserve the orientation of $S$, and
  \item fix the $n$ marked points.
  \end{inparaenum}
\end{definition}
The Teichmüller space $\T_{g,n}$ is an analytic space and is
homeomorphic to a convex domain in $\setC^{3g - 3 + n}$.  Since
$\T_{g,n}$ is an analytic variety and $\Gamma_{g,n}$ acts
discontinuously with finite stabilizers, $\M_{g,n}$ inherits a
structure of analytic orbifold of complex dimension $3g - 3 + n$.

Since $\T_{g,n}$ is contractible, its equivariant (co)homology with
rational coefficients is isomorphic to the rational (co)homology of
$\M_{g,n}$ (see \cite[VII.7.7]{brown}).

\subsection{The fatgraph cellularization of the moduli spaces of
  marked Riemann surfaces}
\label{sec:cellularization}

An embedding of a fatgraph $G$ is an injective continuous map $\iota\colon G\to
S$, that is, a homeomorphism of $G$ onto $\iota(G) \subseteq S$, such that the
orientation on $S$ induces the cyclic order at the vertices of $\iota(G)$.
\begin{definition}\label{dfn:embedded-rg}
  An embedded fatgraph is a fatgraph $G$ endowed with a
  homeomorphism ${\tilde \iota}$ between $S(G)$ and the ambient
  surface $S$, modulo the action of $\Diff^0(S)$.
\end{definition}
There is an obvious action of $\Gamma_{g,n}$ on the set $\ERG_{g,n}$
of fatgraphs embedded into $n$-marked Riemann surfaces of genus $g$.

If confusion is likely to arise, we shall speak of \emph{abstract}
fatgraphs, to mean the topological and combinatorial objects
defined in Definition~\ref{dfn:fatgraphs-geom}, as opposed to \emph{embedded}
fatgraphs as in Definition~\ref{dfn:embedded-rg} above.

\begin{definition}
  \label{dfn:metric-ribbon-graphs}
  A metric $\ell$ on a fatgraph $G$ is an assignment of a real
  positive number $\ell_\alpha$ for each edge $\alpha \in \Edges{G}$.
\end{definition}

Given a metric $\ell$ on a fatgraph $G$, the ``thickening''
construction for fatgraphs can be extended to endow the surface $S(G)$
with a conformal structure dependent on $\ell$.  Conversely, a theorem
due to Jenkins and Strebel guarantees that a metric can be defined on
each fatgraph embedded in a surface $S$, depending uniquely on the
conformal structure on $S$.

Let $G$ be a fatgraph (embedded or abstract) of genus $g$ with $n$
marked boundary components.  The set $\Delta(G) = \{ (G, \ell) \}$ of
metrics on $G$ has an obvious structure of topological cell; now glue
these cells by stipulating that $\Delta(G')$ is the face $\ell_\alpha
= 0$ of $\Delta(G)$ when $G'$ is obtained from $G$ by contraction of
the edge $\alpha$.  The topological spaces obtained by this gluing
instructions are denoted $\Tcomb_{g,n}$ (when using \emph{embedded}
fatgraphs), or $\Mcomb_{g,n}$ (when using \emph{abstract} fatgraphs).
The following theorem clarifies their relation to the Teichmüller and
the moduli space; details can be found, e.g., in
\cite[Section~4.1]{mondello:arXiv:0705.1792v1}.
\begin{theorem}\label{thm:HMT}
  The thickening construction induces orbifold isomorphisms:
  \begin{equation*}
    \T_{g,n} \times \setR^n \simeq \Tcomb_{g,n},
    \qquad
    \M_{g,n} \times \setR^n \simeq \Mcomb_{g,n},
  \end{equation*}
\end{theorem}

Call $M(G)$ the cell in $\Mcomb_{g,n}$ corresponding to an abstract
fatgraph $G$, and $T(\tilde G)$ the cell in $\Tcomb_{g,n}$
corresponding to an embedded fatgraph $\tilde G$.

The functorial action of $\Gamma_{g,n}$ on $\ERG_{g,n}$ induces an
action on $\Tcomb_{g,n}$, which permutes cells $T({\tilde G})$ by PL
isomorphisms.
\begin{lemma}
  \label{lemma:penner-kontsevich-bridge}
  $\Mcomb_{g,n}$ is the quotient space of $\Tcomb_{g,n}$ by the
  cellular action of the mapping class group $\Gamma_{g,n}$; the
  projection homomorphism commutes with the isomorphisms in
  Theorem~\ref{thm:HMT}.
\end{lemma}

\begin{lemma}
  \label{lemma:isotropy}
  The isotropy group $\Gamma_{{\tilde G}}$ of the cell $T({\tilde G}) \hookrightarrow
  \Tcomb_{g,n}$ is (isomorphic to) the automorphism group $\Aut G$ of
  the abstract fatgraph $G$ underlying ${\tilde G}$.
\end{lemma}
The action of $\Gamma_{g,n}$ commutes with the face operators, so
$M(G)$ is a face of $M(G')$ iff $G'$ is obtained from $G$ by
contraction of a non-loop edge.

\subsection{Equivariant homology of $\T_{g,n}$ and the complex of fatgraphs}
\label{sec:rg-complex}

\begin{definition}\label{dfn:orientation}
  An orientation of a fatgraph $G$ is an orientation of the vector
  space $\setQ\Edges{G}$, that is, the choice of an order of the edges of
  $G$, up to even permutations.
\end{definition}
Giving an orientation on $G$ (resp.~${\tilde G}$) is the same as
orienting the simplex $\Delta(G)$ (resp.~$T({\tilde G})$).

If $G$ is a fatgraph with $p$ edges, let $W_G := \bigwedge^p
\setQ\Edges{G}$ be the 1-dimensional vector space generated by the wedge
products $\alpha_1 \land \ldots \land \alpha_p$ of edges of~$G$.
Every $f \in \Aut G$ induces a map $f: \Edges{G} \to \Edges{G}$ on the
edges and thus a map $f_*: \alpha_1 \land \ldots \land \alpha_p
\mapsto f(\alpha_1) \land \ldots \land f(\alpha_p)$.  Trivially,
$f_*(\alpha_1 \land \ldots \land \alpha_p) = \pm \alpha_1 \land \ldots
\land \alpha_p$, depending on whether $f$ preserves or reverses the
orientation of $G$.
\begin{definition}\label{dfn:orientable}
  A fatgraph $G$ is \emph{orientable} iff it has no
  orientation-reversing automorphisms.
\end{definition}

Form a differential complex of orientable fatgraphs as follows.
\begin{definition}
  The complex $(W_*, D)$ of orientable fatgraphs is defined by:
  \begin{itemize}
  \item $W_p := \bigoplus_G W_G$, where $G$ runs over orientable fatgraphs
    with $(2g + n - 1 + p)$ edges;
  \item $D := \sum_1^p (-1)^i d_i$, where $d_i\colon W_p \to W_{p-1}$ is given by:
    \begin{equation*}
      d_i(\alpha_1 \land \ldots \land \alpha_p) :=
      \begin{cases}
        \alpha_1 \land \ldots \land \widehat{\,\alpha_i\,} \land
        \ldots \land \alpha_p
        &\parbox[t]{10em}{%
          if $\alpha_i$ is not a loop and $G/\alpha_i$ is orientable,}
        \\
        0 &\text{otherwise.}
      \end{cases}
    \end{equation*}
  \end{itemize}
\end{definition}


Every oriented fatgraph $(G, \omega)$ defines an element $\omega_G \in
W_G$ by taking the wedge product of edges of $G$ in the order given by
$\omega$; conversely, any $\alpha_1 \land \ldots \land \alpha_p \in
W_G$ defines an orientation on $G$ by setting $\omega := \alpha_1 <
\ldots < \alpha_p$.

\begin{theorem}\label{thm:fatgraph-homology}
  The $\Gamma_{g,n}$-equivariant homology of $\T_{g,n}$ with rational
  coefficients is computed by the complex of oriented fatgraphs
  $(W_*, D)$, i.e., there exists an isomorphism:
  \begin{equation*}
    H^{\Gamma_{g,n}}_*(\T_{g,n}, \setQ) \cong H_*(W_*, D).
  \end{equation*}
\end{theorem}
\begin{proof}
  The genus and number of boundary cycles will be fixed throughout, so for
  brevity, set $\Gamma := \Gamma_{g,n}$, $\T := \T_{g,n}$ and $\Tcomb :=
  \Tcomb_{g,n}$.  

  By Theorem~\ref{thm:HMT}, we have:
  \begin{equation*}
    H^\Gamma_*(\T, \setQ) = H^\Gamma_*(\Tcomb, \setQ). 
  \end{equation*}
  Recall that $H^\Gamma_*(\Tcomb, \setQ)$ can be defined as the homology of the
  double complex $P_* \otimes C_*(\Tcomb, \setQ)$, where $P_*$ is any projective
  resolution of $\setQ$ over $\setQ[\Gamma]$.  The spectral sequence $E^1_{pq} :=
  H_q(P_* \otimes_\Gamma C_p) = H_q(\Gamma, C_p)$ abuts to $H^\Gamma_{p+q}(\Tcomb)$ (see
  \cite[VII.5 and VII.7]{brown}).

  The space $\Tcomb_{g,n}$ has, by definition, an equivariant
  cellularization with cells indexed by embedded fatgraphs of
  genus $g$ with $n$ marked boundary components.  Let $R_p$ be a set
  of representatives for the orbits of $p$-cells under the action of
  $\Gamma$.  By Lemma~\ref{lemma:penner-kontsevich-bridge}, $R_p$ is in
  bijective correspondence with the set of abstract fatgraphs
  having $p$ edges, and the orientation of a cell translates directly
  to an orientation of the corresponding graph.  For each geometric
  simplex $T({\tilde G}) \subseteq \Tcomb$, let $\Gamma_{\tilde G}$ be its
  isotropy group, and let $\setQ{\tilde G}$ be the $\Gamma_{\tilde
      G}$-module consisting of the $\setQ$-vector space generated by an
  element $\Delta$ on which $\Gamma_{\tilde G}$ acts by the orientation
  character: $\tau \cdot \Delta = \pm \Delta$ depending on whether $\tau$ preserves or
  reverses the orientation of the cell $T({\tilde G})$.  By
  Lemma~\ref{lemma:isotropy}, there is an isomorphism between $\Gamma_{\tilde
      G}$ and $\Aut G$; if $\tau \in \Gamma_{\tilde G}$ reverses
  (resp.~preserves) orientation of $T({\tilde G})$, then the
  corresponding $f \in \Aut G$ reverses (resp.~preserves) orientation on
  $G$.  Therefore, $\setQ{\tilde G}$ and $W_G$ are isomorphic as $\Aut G =
  \Gamma_{T{\tilde G}}$ modules.

  Following \cite[p.~173]{brown}, let us decompose (as a $\Gamma$-module)
  \begin{equation*}
    C_*(\Tcomb, \setQ) = \bigoplus_{G \in R_p} W_G;
  \end{equation*}
  then, by Shapiro's lemma \cite[III.6.2]{brown}, we have:
  \begin{equation*}
    H_q(\Gamma, C_p) \cong \bigoplus_{G \in R_p} H_q(\Gamma_{\tilde G}, \setQ{\tilde G}) \cong
    \bigoplus_{G \in R_p} H_q(\Aut G, W_G). 
  \end{equation*}
  Since $\Aut G$ is finite and we take rational coefficients, then
  $H_q(\Aut G, W_G) = 0$ if $q > 0$ \cite[III.10.2]{brown}.  On the
  other hand, if $G$ is orientable then $\Aut G$ acts trivially on
  $W_G$, so:
  \begin{displaymath}
    H_0(\Aut G, W_G) =
    \begin{cases}
      0
      &\text{if $G$ has an orientation-reversing automorphism,}
      \\
      W_G
      &\text{if $G$ has \emph{no} orientation-reversing automorphisms.}
    \end{cases}
  \end{displaymath}
  Let $R'_p$ be the collection of all orientable fatgraphs with
  $p$ edges.  Substituting back into the spectral sequence, we see
  that only one column survives:
  \begin{align}
    \label{eq:11}
    E^1_{p,0} &= \bigoplus_{G \in R'_p} W_G = W_p,
    &
    \\
    \label{eq:12}
    E^1_{p,q} &= 0
    &\text{for all $q > 0$,}
  \end{align}
  In other words, $E^1_{pq}$ reduces to the complex $(E^1_{*,0}, d^1)$.

  Finally, we show that the differential $d^1\colon E^1_{p,0} \to E^1_{p-1,
    0}$ corresponds to the differential $D\colon W_p \to W_{p-1}$ under the
  isomorphism formula~\eqref{eq:11}; this will end the proof.  Indeed, we
  shall prove commutativity of the following diagram at the chain
  level:
  \begin{equation}
    \label{eq:13}
    \xymatrix{%
      P_* \otimes W_p
      &
      \bigoplus_{G \in R'_p} P_* \otimes W_G
      &
      \bigoplus_{G' \in R'_{p-1}} P_* \otimes W_{G'} 
      & 
      P_* \otimes W_{p-1} 
      \\
      & 
      P_* \otimes C_p(\Tcomb, \setQ)
      &
      P_* \otimes C_{p-1}(\Tcomb, \setQ) 
      &
      \ar^{\idmap_P \otimes D} "1,2";"1,3"
      \ar_{\idmap_P \otimes \partial} "2,2";"2,3"
      \ar_{\theta_p} "1,2";"2,2"
      \ar^{\theta_{p-1}} "1,3";"2,3"
      \ar@{=} "1,1";"1,2"
      \ar@{=} "1,3";"1,4"
    }
  \end{equation}
  which implies commutativity at the homology level:
  \begin{equation*}
    \xymatrix{%
      \bigoplus_{G \in R'_p} H_0(\Aut G, W_G)
      &
      \bigoplus_{G' \in R'_{p-1}} H_0(\Aut G', W_{G'})
      \\
      H_0(\Gamma, C_p(\Tcomb, \setQ))
      &
      H_0(\Gamma, C_{p-1}(\Tcomb, \setQ))
      \ar^{D} "1,1";"1,2"
      \ar_{d^1 = H_0(\Gamma, \partial)} "2,1";"2,2"
      \ar_{\cong} "1,1";"2,1"
      \ar^{\cong} "1,2";"2,2"
    }
  \end{equation*}
  whence the conclusion $E^1_{*,0} \cong (W_*, D)$.

  The vertical maps $\theta_p$, $\theta_{p-1}$ in~\eqref{eq:13} are
  the chain isomorphisms underlying the $\Gamma$-module decomposition
  $C_p(\Tcomb, \setQ) \cong \bigoplus_{G \in R_p} W_G$.  Taking the
  boundary of a cell $T({\tilde G}) \subseteq \Tcomb$ commutes with
  the $\Gamma$-action: $\partial T(\tau \cdot {\tilde G}) = \tau
  \cdot \partial T({\tilde G})$.  Furthermore, $T({\tilde G}')$ is a
  cell in $\partial T({\tilde G})$ iff ${\tilde G}'$ is obtained from
  ${\tilde G}$ by contraction of an edge; but ${\tilde G}'$ is a
  contraction of ${\tilde G}$ iff the underlying \emph{abstract}
  fatgraphs $G'$ and $G$ stand in the same relation.  Thus, the
  $\Gamma$-complexes $(C_*,
  \partial)$ and $(W_*, D)$ are isomorphic by $\theta_*$, so
  diagram~\eqref{eq:13} commutes, as was to be proved.
\end{proof}

\section{Computer representation of Fatgraphs}
\label{sec:Fatgraph}

Although the combinatorial definition of a fatgraph
(cf.~Lemma~\ref{lemma:comb-fatgraph}) lends itself to a computer
representation as a triple of permutations ---as used, e.g., in
\cite[Section 2.4]{penner+knudsen+wiuf+andersen:2009}---, the
functions that are needed by the generation algorithms (see
Section~\ref{sec:generation}) are rather topological in nature and
thus suggest an approach more directly related to the concrete
realization of a fatgraph.

\begin{definition}
  A "Fatgraph" object $\+G$ is comprised of the following data:
  \begin{itemize}
  \item A list "$\+G$.vertices" of "Vertex" objects.
  \item A list "$\+G$.edges" of "Edge" objects.
  \item A set "$\+G$.boundary_cycles" of "BoundaryCycle" objects.
  \item An orientation "$\+G$.orient".
  \end{itemize}
\end{definition}
The exact definition of the constituents of a "Fatgraph" object is the
subject of the following sections; informally, let us say that a
"Vertex" is a cyclic list of edges and that an "Edge" is a pair of
vertices and incidence positions.  A precise statement about the
correspondence of abstract fatgraphs and "Fatgraph" objects is made in
Section~\ref{sec:Fatgraph-category}.

There is some redundancy in the data comprising a "Fatgraph"
object: some of these data are inter-dependent and cannot be specified
arbitrarily.  Actually, all data comprising a "Fatgraph" object can be
computed from the vertex list alone, as the following sections show.



In what follows, the letters $l$, $m$ and $n$ shall denote the number
of vertices, edges and boundary cycles:
\begin{itemize}
\item "$l$ == $\card{V(G)}$ == size($\+G$.vertices)",
\item "$m$ == $\card{E(G)}$ == size($\+G$.edges)",
\item "$n$ == $\card{B(G)}$ == size($\+G$.boundary_cycles)".
\end{itemize}
For integers $\alpha$ and $k$, we use $(\alpha \% k)$ to denote the smallest
non-negative representative of $\alpha\mod k\;$.

\subsection{Vertices}
\label{sec:vertices}

\begin{figure}
  \centering
  \includegraphics{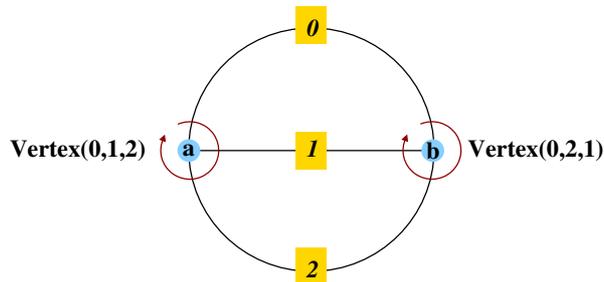}
  \caption{Representation of vertices as (cyclic) lists of edge
    labels; vertices are identified by lowercase Latin letters; edge
    labels are depicted as roman numerals on a yellow square
    background, sitting over the edge they label.  The representation
    of a vertex as a list is implicitly ciliated: here we use the
    convention that the edge closest to the tail of the arrow is the
    ciliated one.}
  \label{fig:vertex}
\end{figure}

We can represent a fatgraph vertex by assigning
labels\footnote{\label{fn:1}%
  Labels can be drawn from any finite set.  In actual computer
  implementations, two obvious choices are to use the set of machine
  integers, or the set of \q{Edge} objects themselves (i.e., label
  each fatgraph edge with the corresponding computer
  representation).}{\space}%
to all fatgraph vertices and mapping a vertex to the
cyclically-invariant list of labels of incident edges.
Figure~\ref{fig:vertex} gives an illustration.

\begin{definition}
  A vertex together with a choice of an attached edge is called a
  \emph{ciliated} vertex.  The chosen edge is called the
  \emph{cilium}.
\end{definition}
\begin{definition}
  If $v$ is a ciliated vertex and $e$ is a half-edge attached to it,
  define the \emph{attachment index} of~$e$ at~$v$ as the index of
  edge~$e$ relative to the cilium at~$v$: if $\alpha$ is the
  attachment index of~$e$ at~$v$, then $\sigma_0^\alpha$ takes the
  cilium at~$v$ onto~$e$.
\end{definition}
The attachment index at a vertex is unambiguously defined for all
edges which are not loops; the two half-edges comprising a loop have
distinct attachment indices.  For brevity, in the following we
shall slightly abuse the definition and speak of \emph{the} attachment index
of an edge at a vertex.

\begin{definition}
  A "Vertex" object "$\+v$ == Vertex($e_1$, $\ldots$, $e_z$)" is a
  list of the labels $e_1$, \ldots, $e_z$ of attached edges.

  Two "Vertex" objects are considered equal if one is equal (as a sequence)
  to the other rotated by a certain amount.
\end{definition}
Note that the definition of "Vertex" objects as plain lists
corresponds to \emph{ciliated} vertices in a fatgraph.  In order to
implement the cyclic behavior of fatgraph vertices, the requirement on
equality must be imposed; equality of "Vertex" objects can be tested
by an algorithm of quadratic complexity in the vertex valence.



If $\+v$ is a vertex object, let us denote "num_loops($\+v$)" the
number of loops attached to $\+v$; it is a vertex invariant and will be
used in the computation of fatgraph isomorphisms.  Implementations
of "num_loops" need only count the number of repeated edge labels in
the list defining the "Vertex" object $\+v$.

\subsection{Edges}
\label{sec:edges}

\begin{definition}
  An "Edge" object $\+e$ is an unordered pair of endpoints, so
  defined: each endpoint corresponds to a 2-tuple $(\+v, a)$, where $\+v$
  is a vertex, and $a$ is the index at which edge $\+e$ appears
  within vertex $\+v$ (the attachment index).
\end{definition}
It is clear how an "Edge" object corresponds to a fatgraph edge: a fatgraph
edge is made of two half-edges, each of which is uniquely identified
by a pair formed by the end vertex $\+v$ and the attachment index $a$.
In the case of loops, the two ends will have the form $(\+v, a)$,
$(\+v, a')$ where $a$ and $a'$ are the two distinct attachment indices
at $\+v$.  

The "other_end($\+e$, $\+v$, $a$)" function takes as input an edge
object $\+e$, a vertex $\+v$, and an attachment index $a$ and returns
the endpoint of $\+e$ opposite to $(\+v, a)$.

The notation "Edge("\var{endpoints}")" will be used for an "Edge" object
comprising the specified endpoints.

\begin{figure}
  \centering
  \includegraphics{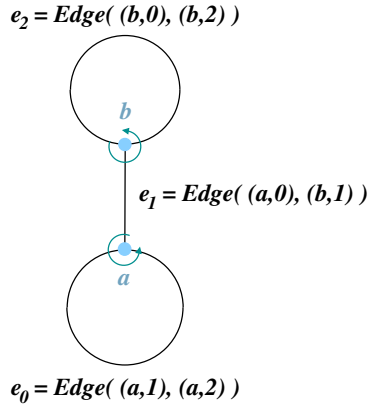}
  \caption{%
    Representation of fatgraph edges.  Each edge is identified with a
    pair of endpoints, where an endpoint is a vertex together with an
    attachment index.  In the figure, letters $a$ and $b$ denote the
    vertices; attachment indices are computed by assigning index $0$
    to the edge closest to the orientation arrow's tail.
  }%
  \label{fig:edges}
\end{figure}

Figure~\ref{fig:edges} provides a graphical illustration of the
representation of fatgraph edges as "Edge" objects.

\subsubsection{Computation of the edge list}
\label{sec:edge-list}

The edge list "$\+G$.edges" can be computed from the list of vertices as follows.

The total number $m$ of edges is computed from the sum of vertex
valences, and used to create a temporary array $P$ of $m$ lists (each
one initially empty). We then incrementally turn $P$ into a list of
edge endpoints (in the form $(\+v, a)$ where $\+v$ is a vertex and $a$
the attachment index) by just walking the list of vertices: $P[k]$ is
the list "[ $(\+v_k, 0)$, $\ldots$, $(\+v_k, z_k)$ ]" where $\+v_k$
(of valence $z_k$) is the $k$-th "Vertex" in "$\+G$.vertices".  The
list "$\+G$.edges" is just $P$ recast into "Edge" objects. In
pseudo-code:
\begin{lstlisting}
  $m$ = $(1/2) \cdot \sum_{\+v \in \+G.vertices} \text{valence}(\+v)$
  $P$ = array of $m$ empty lists
  for $\+v$ in $\+G$.vertices:
    for ($a$, $e$) in enumerate($\+v$):
      |append $(\+v, a)$ to $P[e]$|
  # wrap endpoints into "Edge" objects
  $\+G$.edges = [ Edge($p$) for $p$ in $P$ ]
\end{lstlisting}

\subsection{Boundary Cycles}
\label{sec:boundary-cycles}

\begin{figure}
  \label{fig:corners}
  \centering
  \includegraphics{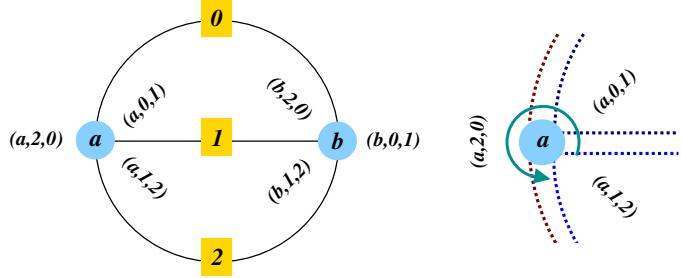}
  \caption{%
    Representation of fatgraph boundary cycles. 
    \emph{Left:} How the boundary cycles are represented with corners:
    each boundary component is identified with the set of triplets it
    encloses. Therefore the boundary cycles for the graph
    are represented by the sets $\{(a,0,1), (b,2,0)\}$, $\{(a,1,2),
    (b,1,2)\}$, and $\{(a,2,0), (b,0,1)\}$.
    \emph{Right:} Zoom around vertex $a$ in the left picture, to show
    the three corners identified with triples $(a, i, j)$.  The
    indices in the triple are attachment indices, i.e., displacement
    relative from the ciliated edge (the one closest to the arrow
    tail); they bear no relation to the labels on the edges (numbers
    on the light yellow background in the left picture).  
  }%
\end{figure}

\begin{definition}
  A "BoundaryCycle" object is a set of \emph{corners} (see
  Figure~\ref{fig:corners}).

  A corner object $\+C$ is a triple "(vertex, incoming, outgoing)",
  consisting of a vertex $\+v$ and two indices "$i$ == $\+C$.incoming",
  "$j$ == $\+C$.outgoing" of consecutive edges (in the cyclic order at
  $\+v$).  In order to have a unique representation of any
  corner, we impose the condition that either $j=i+1$, or $i$ and $j$
  are, respectively, the ending and starting indices of $\+v$
  (regarded as a list).
\end{definition}
It is easy to convince oneself that a "BoundaryCycle" object
corresponds to a boundary cycle as defined in
Section~\ref{sec:Fatgraph}.  Indeed, if $(L, \sigma_0, \sigma_1,
\sigma_2)$ is a fatgraph, then the boundary cycles are defined as the
orbits of $\sigma_2$ on the set $L$ of half-edges; a \emph{(endpoint
  vertex, attachment index)} pair uniquely identifies an half-edge and
can thus be substituted for it.  For computational efficiency reasons,
we add an additional successor index to form the corner triple
$(\+v,i,j)$ so that the action of $\sigma_2$ can be computed from corner
data alone, without any reference to the ambient fatgraph.\footnote{%
  This is important in order to share the same corner objects across
  multiple \q{BoundaryCycle} instances, which saves computer memory.
}%

\begin{figure}
  \label{fig:012012}
  \centering
  \includegraphics{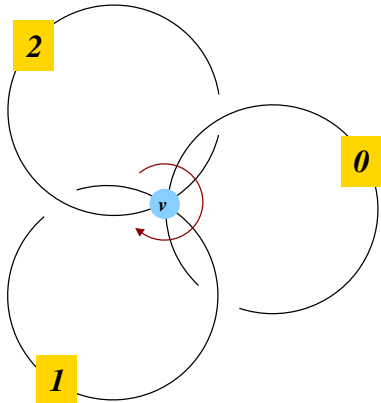}
  \caption{A fatgraph whose two boundary cycles are comprised of
    exactly the same edges; however, they give rise to disjoint sets
    of corners: $\{(\+v, 2, 3), (\+v, 4, 5), (\+v, 0, 1)\}$ versus
    $\{(\+v, 1, 2), (\+v, 3, 4), (\+v, 5, 0)\}$.}
\end{figure}

Since distinct orbits are disjoint, two "BoundaryCycle" objects are
either identical (they comprise the same corners) or have no intersection.
In particular, this representation based on corners distinguishes
boundary cycles made of the same edges: for instance, the boundary
cycles of the fatgraph depicted in Figure~\ref{fig:012012} are
represented by the disjoint set of corners $\{(\+v, 2, 3), (\+v, 4, 5),
(\+v, 0, 1)\}$ and $\{(\+v, 1, 2), (\+v, 3, 4), (\+v, 5, 0)\}$.


\subsubsection{Computation of boundary cycles}
\label{sec:compute-boundary-cycles}

The procedure for computing the set of boundary cycles of a given
"Fatgraph" object $\+G$ is listed in Algorithm~\ref{algo:compute-boundary-cycles}.
\begin{Algorithm}
  \caption{\label{algo:compute-boundary-cycles}%
    Output the set of boundary cycles of a \q{Fatgraph} object
    $\+G$. Input to the algorithm is a \q{Fatgraph} object $\+G$; the
    output is a list of \q{BoundaryCycle} objects.  The special
    constant \textsc{used} marks locations in the temporary array
    \q{corners} whose contents has already been assigned to a boundary
    cycle.}
\begin{lstlisting}
  def compute_boundary_cycles($\+G$):
    # build working array of corners
    corners = [ [ ($\+v$, $i$, ($(i+1)$ % $\card{\+v}$)) for $i$ in $0, \ldots, \card{\+v}-1$ ] |\n{bcy-def-corner-1}|
                for $\+v$ in $\+G$.vertices ]    |\n{bcy-def-corners-2}|
    result = |empty list|
    $l_0$ = 0
    $i_0$ = 0
    while True:                                      |\n{bcy-main-loop}|
      # locate the first unused corner
      for $l$ in $l_0$, $\ldots$, size(corners)-1:  |\n{bcy-ffwd-start}|
        $\+v$ = $\+G$.vertices[$l$]
        $i$ = first_index_not_used(corners[$\+v$], $i_0$)
        if $i$ is not None:                          |\n{bcy-quit-ffwd}|
          |exit ``for'' loop|                         |\n{bcy-ffwd-end}|
      if $l$ == size(corners)-1 and $i$ is None:         |\n{bcy-exit1}|
        # all corners used, mission accomplished
        return result                                    |\n{bcy-exit2}|
      else:
        $l_0$ = $l$
        $i_0$ = $i$
      # walk the boundary cycle and record corners
      start = $(\+v, i)$                             |\n{bcy-bcy-start}|
      triples = |empty list|
      while $(\+v,i)$ != start or size(triples) == 0:     |\n{bcy-loop}|
        triples.append(corners[$\+v$][$i$])
        $j$ = corners[$v$][$i$][2]                     |\n{bcy-follow1}|
        $\+e$ = $\+v$[$j$]                             |\n{bcy-follow2}|
        # mark location as ``used''
        corners[$\+v$][$i$] = |\textsc{used}|        |\n{bcy-mark-used}|
        $(\+v,i)$ = other_end($\+e$, $\+v$, $j$)       |\n{bcy-follow3}|
        $\+b$ = BoundaryCycle(triples)
      result.append($\+b$)                             |\n{bcy-bcy-end}|

  def first_index_not_used($L$):
      for index, item in enumerate($L$):
          if item is not |\textsc{used}|:
              return index
      return None
\end{lstlisting}
\end{Algorithm}
The algorithm closely follows a geometrical procedure: starting with
any corner, follow its ``outgoing'' edge to its other endpoint,
and repeat until we come back to the starting corner.  The list of
corners so gathered is a boundary cycle.  At each iteration, the used
corners are cleared out of the "corners" list by replacing them with
the special value \textsc{used}, so that they will not be picked up again in
subsequent iterations.

\begin{lemma}
  For any "Fatgraph" object $\+G$ representing a fatgraph $G$,
  the function "compute_boundary_cycles" in
  Algorithm~\ref{algo:compute-boundary-cycles} has the following
  properties:
  \begin{inparaenum}[\it 1)]
    \item terminates in finite time, and
    \item returns a list of "BoundaryCycle" objects that represent the
      boundary cycles of $G$.
  \end{inparaenum}
\end{lemma}
\begin{proof}
  The algorithm works on a temporary array "corners": as it walks
  along a boundary cycle (lines \nr{bcy-loop}--\nr{bcy-follow3}),
  corner triples are moved from the working array to the "triples"
  list and replaced with the constant \textsc{used}; when we're back
  to the starting corner, a "BoundaryCycle" object $\+b$ is
  constructed from the "triples" list and appended to the result.

  The "corners" variable is a list, the $n$-th item of which is
  (again) a list holding the corners around the $n$-th vertex (i.e.,
  "$\+G$.vertices[$n$]"), in the order they are encountered when
  winding around the vertex.  By construction, "corners[$\+v$][$i$]"
  has the the form $(\+v,i,j)$ where $j$ is the index following $i$ in
  the cyclic order, i.e., $(\+v,i,j)$ represents the corner formed by
  the ``incoming'' $i$-th edge and the ``outgoing'' $j$-th edge.

  The starting corner for each walk along a boundary cycle is
  determined by scanning the "corners" list (lines
  \nr{bcy-ffwd-start}--\nr{bcy-ffwd-end}): loop over all indexes
  $\+v$, $i$ in the "corners" list, and quit looping as soon as
  "corners[$\+v$][$i$]" is not \textsc{used}
  (line~\nr{bcy-quit-ffwd}).  If all locations in the corners list are
  \textsc{used}, then the all corners have been assigned to a boundary
  cycle and we can return the result list to the caller.
\end{proof}

\subsection{Orientation}
\label{sec:orientation}

According to Definition~\ref{dfn:orientation}, orientation is given by
a total order of the edges (which directly translates into an
orientation of the associated orbifold cell).

\begin{definition}\label{sec:orientation-1}
  The orientation "$\+G$.orient" is a list that
  associates each edge with its position according to the order given
  by the orientation.  
  Two such lists are equivalent if they differ by an even permutation.
\end{definition}
If $\+e_1$ and $\+e_2$ are edges in a "Fatgraph" object $\+G$, then
$\+e_1$ precedes $\+e_2$ iff 
"$\+G$.orient[$\+e_1$] < $\+G$.orient[$\+e_2$]"; 
this links the fatgraph orientation from
Definition~\ref{dfn:orientation} with the one above.

If a "Fatgraph" object is derived from another "Fatgraph" instance
(e.g., when an edge is contracted), the resulting graph must derive
its orientation from the ``parent'' graph, if we want the edge
contraction to correspond to taking cell boundary in the
orbicomplex~$\Mcomb_{g,n}$.

When no orientation is given, the trivial one is (arbitrarily) chosen:
edges are ordered in the way they are listed in the "$\+G$.edges"
list, i.e., "$\+G$.orient[$\+e$]" is the position at which $\+e$ appears in
"$\+G$.edges".

According to Definition~\ref{dfn:orientable}, a fatgraph is orientable
iff it has no orientation-reversing automorphism.  The author knows of
no practical way to ascertain if a fatgraph is orientable other than
enumerating all automorphisms and checking if any one of them reverses
orientation:
\begin{lstlisting}
  def is_oriented($\+G$):
    for $\+a$ in automorphisms($\+G$):
      if is_orientation_reversing($\+a$):
        return False
    # no orientation-reversing automorphism found, $\+G$ is orientable
    return True

\end{lstlisting}

\subsection{A category of \q{Fatgraph} objects}
\label{sec:Fatgraph-category}

\subsubsection{Isomorphisms of \q{Fatgraph} objects}
\label{sec:iso1}

In this section, we shall only give the definition of "Fatgraph"
isomorphisms and prove the basic properties; the algorithmic
generation and treatment of "Fatgraph" isomorphisms is postponed to
Section~\ref{sec:isomorphism}.

\begin{definition}
  \label{dfn:iso}
  An isomorphism of "Fatgraph" objects $\+G_1$ and $\+G_2$ is a triple 
  $\+f =$"(pv, rot, pe)" where:
  \begin{itemize}
  \item "pv" is a permutation of the vertices: vertex $\+v_1$ of
    $\+G_1$ is sent to vertex "pv[$\+v$]" of $\+G_2$, and rotated by
    "rot[$\+v$]" places leftwards;
  \item "pe" is a permutation of the edge labels: edge $\+e$ in
    $\+G_1$ is mapped to edge "pe[$\+e$]" in $\+G_2$.
  \end{itemize}
  The adjacency relation must be preserved by isomorphism triples: if
  $\+v_1$ and $\+v_2$ are endpoint vertices of the edge $\+e$, then
  "pv[$\+v_1$]" and "pv[$\+v_2$]" must be the endpoint
  vertices of edge "pe[$\+e$]" in $\+G_2$.
\end{definition}
Since a vertex in a "Fatgraph" instance is essentially the list of labels
of edges attached to that vertex, we can dually state the
compatibility condition above as requiring that, for any vertex $\+v$
in "$\+G_1$.vertices" and any valid index $j$ of an edge of $\+v$, we
have:
\begin{equation}
  \label{eq:fg-iso}
  \text{\lstinline:$\+G_2$.vertices[pv[$\+v$]][$j$+rot[$\+v$]] == pe[$\+G_1$.vertices[$\+v$][$j$]]:}
\end{equation}
The above formula~\eqref{eq:fg-iso} makes the parallel between "Fatgraph" object
isomorphisms and fatgraph maps (in the sense of
Definition~\ref{dfn:fatgraphs-geom}) explicit.

\begin{lemma}
  \label{lemma:iso}
  Let $G_1$, $G_2$ be fatgraphs, represented respectively by $\+G_1$
  and $\+G_2$. Every isomorphism of fatgraphs $f\colon G_1 \to G_2$ lifts
  to a corresponding isomorphism "$\+f$ == (pv, rot, pe)" on the computer
  representations. Conversely, every triple "(pv, rot, pe)"
  representing an isomorphism between the "Fatgraph" instances induces
  a (possibly trivial) fatgraph isomorphism between $G_1$ and $G_2$.
\end{lemma}
\begin{proof}
  Every isomorphism $f\colon G_1 \to G_2$ naturally induces bijective maps
  $f_V\colon V(G_1) \to V(G_2)$ and $f_E\colon E(G_1) \to E(G_2)$ on vertices
  and edges.  Given a cilium on every vertex, $f$ additionally
  determines, for each vertex $v \in V(G)$, the displacement
  $f_{\textsl{rot}}(v)$ of the image of the cilium of $v$ relative to
  the cilium of $f_V(v)$.
  Similarly, $f_E$ determines a bijective mapping of edge labels, and
  is completely determined by it.  This is exactly the data collected
  in the triple "(pv, rots, pe)", and the compatibility
  condition~\eqref{eq:fg-iso} holds by construction.

  Conversely, assume we are given a triple "(pv, rots, pe)",
  representing an isomorphism of "Fatgraph" instances.
  We can construct maps $f_V$, $f_E$ as follows: $f_V$ sends a vertex
  $\+v \in \+G_1$ to the vertex corresponding to "pv[$\+v$]"; $f_E$ maps
  the cilium of $\+v$ to the edge attached to "pv[$\+v$]" at "rot[$\+v$]"
  positions away from the cilium; the compatibility condition
  \eqref{eq:fg-iso} guarantees that $f_E$ is globally well-defined.
\end{proof}

\begin{lemma}
  \label{lemma:1}
  Let $\+G_1$, $\+G_2$ be "Fatgraph" objects, and $\eta$ a bijective
  map between "$\+G_1$.edges" and "$\+G_2$.edges" that preserves the
  incidence relation.  Then there is a unique "Fatgraph" isomorphism
  $\+f$ that extends $\eta$ (in the sense that "$\+f$.pe == $\eta$").
\end{lemma}
\begin{proof}
  Start constructing the "Fatgraph" morphism $\+f$ by setting
  "$\+f$.pe == $\eta$".  If $\+e_1$,~\ldots,~$\+e_{z_k}$ are the edges
  incident to "$\+v_k \in \+G_1$.vertices", then there is generally one
  and only one endpoint $\+v_k'$ common to edges $\eta(\+e_k)$; define
  "$\+f$.pv[$\+v_k$] == $\+v_k'$".  There is only one case in which
  this is not true, namely, if \emph{all} edges share the same two
  endpoints:\footnote{%
    So there are only two vertices in total, and the
    corresponding fatgraph belongs in $\R_{0,m}$.
  }\space%
  in this case, however, there is still only one choice of
  "$\+f$.pv[$\+v_k$]" such that the cyclic order of edges at the
  source vertex matches the cyclic order of edges at the target
  vertex. Finally, choose "$\+f$.rot[$\+v_k$]" as the displacement
  between the cilium at $\+v_k'$ and the image of the cilium of
  $\+v_k$.

  It is easy to check that eq.~\eqref{eq:fg-iso} holds, so $\+f$ is a
  well-defined isomorphism.
\end{proof}

\subsubsection{Contraction morphisms}
\label{sec:contraction}

Recall from the definition in Section~\ref{sec:mgn} that
contraction produces a ``child'' fatgraph from a ``parent'' fatgraph
and a chosen regular (i.e., non-looping) edge. 

The "Fatgraph.contract" method (see Algorithm~\ref{algo:contraction})
thus needs only take as input the ``parent'' graph $\+G$ and the edge
$\+e$ to contract, and produces as output the ``child'' fatgraph
$\+G'$.  The contraction algorithm proceeds in the following way:
\begin{itemize}
\item The two end vertices of the edge $\+e$ are fused into one: the
  list "$\+G'$.vertices" is built by copying the list
  "$\+G$.vertices", removing the two endpoints of $\+e$, and adding
  the new vertex (resulting from the collapse of $\+e$) at the end.
\item Deletion of an edge also affects the orientation: the
  orientation "$\+G'$.orient" on the ``child'' fatgraph keeps the
  edges in the same order as they are in the parent fatgraph.
  However, since "$\+G'$.orient" must be a permutation of the edge
  indices, we need to renumber the edges and shift the higher-numbered
  edges down one place.
\item The ``child'' graph $\+G'$ is constructed from the list
  "$\+G'$.vertices" and the derived orientation "$\+G'$.orient"; 
  the list of ``new'' edges is constructed according to the procedure
  given in Section~\ref{sec:edge-list}.
\end{itemize}
Listing~\ref{algo:contraction} summarizes the algorithm applied.
\begin{Algorithm}
  \caption{Construct a new \q{Fatgraph} object $\+G'$ obtained by
    contracting the edge $\+e$ in $\+G$. The renumbering function $s$
    is the identity on numbers in the range $0$, \ldots, $e-1$, and
    shifts numbers in range $e+1$, \ldots, $m$ down by 1. 
    Function \q{rotated}$(L,p)$ returns a copy of list $L$ shifted
    leftwards by $p$ places.}
  \label{algo:contraction}
\begin{lstlisting}
  def contract($\+G$, $\+e$):
    |let ($\+v_1$, $a_1$), ($\+v_2$, $a_2$) be the endpoints of $\+e$|
    $V'$ = [ Vertex($x$ for $x$ in $\+v$ if $x$ != $\+e$) 
            for $\+v$ in $\+G$.vertices if $\+v$ != $\+v_1$ and $\+v$ != $\+v_2$ ]
    # append the fused vertex at end of list $V$
    $\+v'$ = Vertex(rotated($\+v_1$, $a_1$) + rotated($\+v_2$, $a_2$))
    $V'$.append($\+v'$)
    $\omega'$ = [ $s$($\+G'$.orient[$x$]) for $x$ in $\+G'$.edges if $x$ != $e$ ] |\n{contract-omega}|
    return Fatgraph(vertices = $V'$; orient = $\omega'$)

\end{lstlisting}
\end{Algorithm}

The vertex resulting from the contraction of $\+e$ is formed as follows.
Assume $\+v_1$ and $\+v_2$ are the endpoint vertices of the
contracted edge.  Now fuse endpoints of the contracted edge:
\begin{enumerate}
\item Rotate the lists $\+v_1$, $\+v_2$ so that the given edge $\+e$ appears
  \emph{last} in $\+v_1$ and \emph{first} in $\+v_2$. 
\item Form the new vertex $\+v$ by concatenating the two rotated lists
  (after expunging vertices $\+v_1$ and $\+v_2$).
\end{enumerate}
Note that this changes the attachment indices of all edges incident to
$\+v_1$ and $\+v_2$, therefore the edge list of $\+G'$ needs to be
recomputed from the vertex list.

The ``child'' fatgraph $G'$ inherits an orientation from the
``parent'' fatgraph, which might differ from its default orientation.
Let $\alpha_1, \ldots, \alpha_h, \ldots, \alpha_m$ be the edges of the
parent fatgraph $G$, with $e = \alpha_h$ being contracted to create
the ``child'' graph $G'$.  If $\alpha_{k(1)} < \alpha_{k(2)} < \ldots
< \alpha_{k(m)}$ is the ordering on $E(G)$ that induces the
orientation on $G$ and $h = k(j)$, then $\alpha_{k(1)} < \ldots <
\alpha_{k(j-1)} < \alpha_{k(j+1)} < \ldots < \alpha_{k(m)}$ descends
to a total order on the edges of $G'$ and induces the correct
orientation.\footnote{%
  That is to say, the orientation that corresponds to the
  orientation induced on the cell $\Delta(G')$ as a face of
  $\Delta(G)$.
}%

Orientation is represented in a "Fatgraph" object as a list, mapping
edge labels to a position in the total order; using the notation
above, the orientation of $G$ is given by $\omega := k^{-1}$.
The orientation on $G'$ is then given by $\omega'$ defined as follows:
\begin{equation*}
  \omega'(i) :=
  \begin{cases}
    \omega(i)    &\text{if $\omega(i)<h$,} 
    \\
    \omega(i)-1  &\text{if $\omega(i)>h$.}
  \end{cases}
\end{equation*}
Alternatively we can write:
\begin{equation*}
  \omega' = s \circ \omega, 
  \qquad
  s(x) :=
  \begin{cases}
    x   &\text{if $x < h$,}
    \\
    x-1 &\text{if $x > h$.}
  \end{cases}
\end{equation*}
This corresponds exactly to the assignment in
Algorithm~\ref{algo:contraction}.

The above discussion can be summarized in the following.
\begin{lemma}
  \label{lemma:K-contraction}
  If $\+G$ and $\+G'$ represent fatgraphs $G$ and $G'$, and 
  "$\+G$ == contract($\+G'$, $\+e$)", then $G$ is obtained from $G'$
  by contraction of the edge $e$ represented by $\+e$.
\end{lemma}

\paragraph{The \q{contract\_boundary\_cycle} function.}
The boundary cycles of the ``child'' "Fatgraph" object $\+G'$ can also
be computed from those of $\+G$.  The implementation (see
Listing~\ref{algo:compute-boundary-cycles}) is quite straightforward:
we copy the given list of corners and alter those who refer to the two
vertices that have been merged in the process of contracting the
specified edge.
\begin{Algorithm}
  \caption{Return a new \q{BoundaryCycle} instance, image of $\+b$ under
    the topological map that contracts the edge with index $\+e$.}
  \label{algo:contract-bcy}
\begin{lstlisting}
  def contract_boundary_cycle($\+G$, $\+b$, $\+e$):
    |let ($\+v_1$, $a_1$), ($\+v_2$, $a_2$) be the endpoints of $\+e$|
    $z_1$ = valence($\+v_1$)
    $z_2$ = valence($\+v_2$)
    # ``child'' boundary cycle $\+b'$ starts off as an empty list
    $\+b'$ = [ ]                                                   |\n{cbcy:1}|
    for corner in $\+b$:                                           |\n{cbcy:2}|
      if corner[0] == $\+v_1$:
        if $a_1$ == corner.incoming:                                   |\n{cbcy:4}|
          continue |with next \q{corner}|                          |\n{cbcy:5}|
        else: 
          $i_1$ = (corner.incoming - $a_1$ - 1) % $z_1$                |\n{cbcy:8}|
          $i_2$ = (corner.outgoing - $a_1$ - 1) % $z_1$                |\n{cbcy:9}|
          |append corner ($\+v_1$, $i_1$, $i_2$) to $\+b'$|
      elif corner[0] == $\+v_2$:
        if $a_2$ == corner.incoming:                                   |\n{cbcy:6}|
          continue |with next \q{corner}|                          |\n{cbcy:7}|
        if $a_2$ == corner.outgoing:
          |append $(\+v_1, z_1+z_2-3, 0)$ to $\+b'$|                 |\n{cbcy:12}|
        else:
          $i_1$ = $z_1 - 1$ + ((corner.incoming - $a_2$ - 1) % $z_2$)  |\n{cbcy:10}|
          $i_2$ = $z_1 - 1$ + ((corner.outgoing - $a_2$ - 1) % $z_2$)  |\n{cbcy:11}|
          |append $(v_1, i_1, i_2)$ to $b'$|
      else:
        # keep corner unchanged
        |append \q{corner} to $\+b'$|                                |\n{cbcy:3}|
    return BoundaryCycle($\+b'$)
\end{lstlisting}
\end{Algorithm}

Let $\+v_1$ and $\+v_2$ be the end vertices of the edge to be contracted,
and $a_1$, $a_2$ be the corresponding attachment indices.  Let $z_1$ and
$z_2$ be the valences of vertices $\+v_1$, $\+v_2$.  We build the list of
corners of the boundary cycle in the ``child'' graph incrementally:
the $\+b'$ lists starts empty (line~\nr{cbcy:1}), and is then added
corners as we run over them in the loop between lines~\nr{cbcy:2}
and~\nr{cbcy:3}.

\begin{figure}
  \centering
  \includegraphics{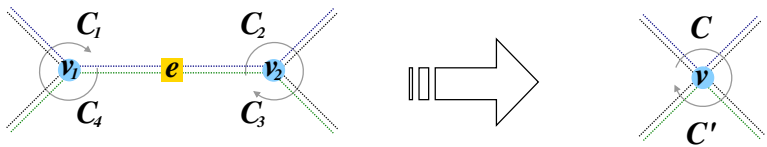}
  \caption{%
    How corners are modified by edge contraction.
    \emph{Left:} Four distinct corners are formed at the endpoints
    $\+v_1$, $\+v_2$ of edge $\+e$, which is to be contracted: $\+C_1
    = (\+v_1, 0, 2)$, $\+C_2 = (\+v_2, 0, 1)$, $\+C_3 = (\+v_2, 1,
    2)$, and $\+C_4 = (\+v_1, 0, 1)$. Edges are shown thickened, and
    (potentially) distinct boundary cycles are drawn in different colors.
    \emph{Right:} After contraction of $\+e$, corners $\+C_1$ and
    $\+C_2$ are fused into $C = (\+v, 0, 1)$, and $\+C_3$, $\+C_4$ are
    fused into $\+C' = (\+v, 2, 3)$.
  }%
  \label{fig:contract-bcy}
\end{figure}

There are four distinct corners that are bounded by the edge $\+e$ to
be contracted; denote them by $\+C_1$, $\+C_2$, $\+C_3$, $\+C_4$.  These map
onto two distinct corners $\+C$, $\+C'$ after contraction.
Assume that $\+C_1$ and $\+C_2$ map to $\+C$: then $\+C_1$ and $\+C_2$ lie
``on the same side'' of the contracted edge, i.e., any boundary cycle
that includes $\+C_1$ will include also $\+C_2$ and viceversa. (See
Figure~\ref{fig:contract-bcy} for an illustration.)  Since they both
map to the same corner $\+C$ in the ``child'' graph, we only need
to keep one: we choose to keep (and transform) the corner that has the
contracted edge at the \emph{second} index
(lines~\nr{cbcy:4}--\nr{cbcy:5}); similarly for $\+C_3$ and $\+C_4$ in
mapping to $\+C'$ (lines~\nr{cbcy:6}--\nr{cbcy:7}).

Recall that, when contracting an edge with endpoints $\+v_1$ and $\+v_2$,
the new vertex is formed by concatenating two
series of edges:
\begin{inparaenum}
\item edges attached to the former $\+v_1$, starting with the successor
  (in the cyclic order) of the contracted edge;
\item edges attached to the former $\+v_2$, starting with the successor
  of the contracted edge.
\end{inparaenum}
Therefore:
\begin{enumerate}
\item The image of a corner rooted in vertex $\+v_1$ will have its
  attachment indices rotated leftwards by $a_1 + 1$ positions: the
  successor of the contracted edge has now attachment index 0
  (lines~\nr{cbcy:8}--\nr{cbcy:9}).  Note that the highest attachment
  index belonging into this group is $z_1 - 2$: position $z_1 - 1$
  would correspond to the contracted edge.
\item The image of a corner rooted in vertex $\+v_2$ has its attachment
  indices rotated leftwards by $a_2 + 1$ positions, and shifted up by
  $z_1 - 1$ (lines~\nr{cbcy:10}--\nr{cbcy:11}).  As a special case, when
  the contracted edge is in second position we need to map the corner
  to the corner having attachment index 0 in second position
  (line~\nr{cbcy:12}).
\end{enumerate}
Any other corner is copied with no alterations (line~\nr{cbcy:3}).

\subsubsection{The category of \q{Fatgraph} objects}
\label{sec:fatgraph-cat}

We can now formally define a category of "Fatgraph" objects and
their morphisms. 
\begin{definition}
  $\R^\#$ is the category whose objects are "Fatgraph" objects, and
  whose morphisms are compositions of "Fatgraph" isomorphisms (as
  defined in Section~\ref{sec:isomorphism}) and edge contraction maps.
\end{definition}
More precisely, if $\+G$ and $\+G'$ are isomorphic "Fatgraph" objects,
then the morphism set $\R^\#(\+G, \+G')$ is defined as the set of
"Fatgraph" isomorphisms in the sense of Section~\ref{sec:isomorphism};
otherwise, let $m$ and $m'$ be the number of edges of $\+G$, $\+G'$,
and set $k := m - m'$: each element in $\R^\#(\+G, \+G')$ has the form
$\+a' \circ (\pi_1 \circ \cdots \circ \pi_k) \circ \+a$ where $\+a$, $\+a'$
are automorphisms of $\+G$, $\+G'$ and $\pi_1$, \ldots, $\pi_k$ are
non-loop edge contractions.

\begin{figure}
  \label{fig:ctor}
  \centering
  \includegraphics{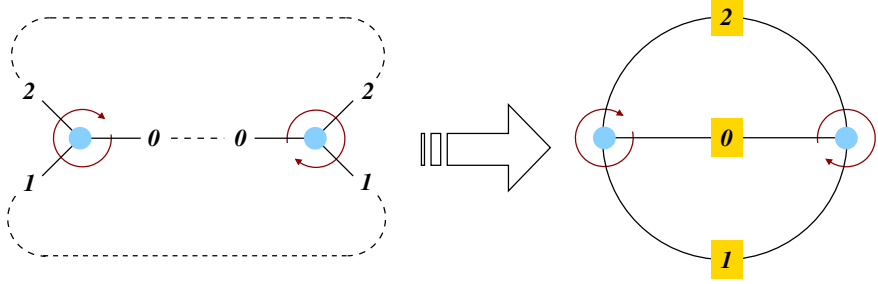}
  \caption{Construction of a fatgraph out of a set of
    \q{Vertex} instances: half-edges tagged with the same
    (numeric) label are joined together to form an edge.}
\end{figure}

\begin{theorem}
  \label{thm:K}
  There exists a functor $K$ from the category $\R^\#$ of "Fatgraph"
  objects to the category $\R$ of abstract fatgraphs, which is
  surjective and full.
\end{theorem}
\begin{proof}
  Given a "Fatgraph" $\+G$, its constituent "Vertex" objects determine
  cyclic sequences $v_0 = (e_0^0, e^0_1, \ldots, e^0_{z_0})$, \ldots,
  $v_l = (e^l_0, \ldots, e^l_{z_l})$, such that $$\{e^0_0, \ldots,
  e^0_{z_0}, e^1_0, \ldots, e^{l-1}_{z_{l-1}}, e^l_0, \ldots,
  e^l_{z_l}\} = \{0, \ldots, m-1\}.$$  Fix a starting element for each
  of the cyclic sequences $v_0$, \ldots, $v_l$.  Then set:
  \begin{equation*}
    L := \{ (e, i, v) : v = v_j \in \{v_0, \ldots, v_l\},\; 
    e = e^j_i \in v \},
  \end{equation*}
  and define maps $\sigma_0, \sigma_1, \sigma_2\colon L \to L$ as follows:
  \begin{itemize}
  \item $\sigma_0$ sends $(e, i, v_j)$ to $(e', i', v_j)$ where $i' =
    (i+1) \% z_j$ and $e' = e^j_{i'}$ is the successor of $e$ in the
    cyclic order at $v_j$;
  \item $\sigma_1$ maps $(e, i, v)$ to the unique other triplet $(e',
    i', v') \in L$ such that $e = e'$;
  \item finally, $\sigma_2$ is determined by the constraint $\sigma_0
    \circ \sigma_2 = \sigma_1$.
  \end{itemize}
  Then $K(\+G) = (L, \sigma_0, \sigma_1, \sigma_2)$ is a fatgraph.
  Figure~\ref{fig:ctor} provides a graphical illustration of the way
  a "Fatgraph" object is constructed out of such combinatorial data.

  Now let $G$ be an abstract fatgraph; assuming $G$ has $m$ edges, assign
  to each edge a ``label'', i.e., pick a bijective map $e: \Edges{G}
  \to E$, where $E$ is an arbitrary finite set.  Each vertex $v \in
  V(G)$ is thus decorated with a cyclic sequence of edge labels; the
  set of which determines a "Fatgraph" object $\+G$; it is clear that
  $G = K(\+G)$.

  This proves that $K$ is surjective; since every fatgraph morphism
  can be written as a composition of isomorphisms and edge
  contractions (Lemma~\ref{lemma:contraction1}), it is also full.  
  It is clear that every edge contraction is the image of an edge
  contraction in the corresponding "Fatgraph" objects, and the
  assertion for isomorphisms follows as a corollary of
  Lemma~\ref{lemma:iso}.
\end{proof}

\begin{definition}
  If $G = K(\+G)$ then we say that the "Fatgraph" object $\+G$
  \emph{represents} the abstract fatgraph $G$.
\end{definition}
It is clear from the construction above that there is a considerable
amount of arbitrary choices to be made in constructing a
representative "Fatgraph"; there are thus many representatives for the
same fatgraph, and different choices lead to equivalent "Fatgraph"
objects.

\begin{lemma}
  Two distinct "Fatgraph" objects representing the same abstract
  fatgraph are isomorphic.
\end{lemma}
\begin{proof}
  Assume $\+G_1$ and $\+G_2$ both represent the same abstract fatgraph
  $G = K(\+G_1) = K(\+G_2)$.  Let $\eta_1$, $\eta_2$ be the maps that
  send "Edge" objects in $\+G_1$, $\+G_2$ to the corresponding edges
  in $G$; then $\eta = \eta_1\inv \circ \eta_2$ maps edges of $\+G_1$
  into edges of $\+G_2$ and respects the incidence relation, therefore
  it is the edge part of a "Fatgraph" isomorphism by
  Lemma~\ref{lemma:1}.
\end{proof}

\begin{theorem}
  The categories $\R^\#$ and $\R$ are equivalent.
\end{theorem}
\begin{proof}
  The functor $K$ is surjective and full by Theorem~\ref{thm:K}; that
  it is also faithful follows from the following argument.  Any
  fatgraph morphism is a composition of edge contractions and
  isomorphisms.  Any isomorphism determines, in particular, a map on the
  set of edges, and there is one and only one "Fatgraph" isomorphism
  induced by this map (Lemma~\ref{lemma:1}).  Any edge contraction is
  uniquely determined by the contracted edge: if $f\colon G_1 \to G_2$ is
  the morphism contracting edge $e$ and $G_i = K(\+G_i)$, then $\+f$,
  contraction of the "Edge" object $\+e$ representing $e$, is the sole
  morphism of $\+G_1$ into $\+G_2$ that maps onto~$f$.
\end{proof}

\section{Fatgraphs isomorphism and equality testing}
\label{sec:isomorphism}

The isomorphism problem on computer representations of fatgraphs
consists in finding out when two distinct "Fatgraph" instances
represent isomorphic fatgraphs (in the sense of
Definition~\ref{dfn:fatgraphs-geom}) or possibly the same fatgraph.  Indeed,
the procedure for associating a "Fatgraph" instance to an abstract
fatgraph (see Theorem~\ref{thm:K}) involves labeling all
edges, choosing a starting edge (cilium) on each vertex and
enumerating all vertices in a certain order; for each choice, we get a
different "Fatgraph" instance representing the same (abstract)
fatgraph.

The general isomorphism problem for (ordinary) graphs is a well-known
difficult problem. However, the situation is much simpler for
fatgraphs, because of the following property.
\begin{lemma}[Rigidity Property]\label{lemma:rigidity}
  Let $G_1$, $G_2$ be connected fatgraphs, and $f\colon G_1 \to G_2$ an
  isomorphism.  For any vertex $v \in V(G_1)$, and any edge $x$
  incident to $v$, $f$ is uniquely determined (up to homotopies fixing
  the vertices of $G_i$) by its restriction to $v$ and $x$.
\end{lemma}
In particular, an isomorphism of graphs with ciliated vertices is
completely determined once the image $w=f(v)$ of a vertex $v$ is
known, together with the displacement (relative to the cyclic order at
$w$) of the image of the cilium of $v$ relative to the cilium of
the image vertex $w$.
\begin{proof}
  Consider $f$ as a CW-complex morphism: $f = (f_0, f_1)$ where
  $f_i$ is a continuous map on the set of $i$-dimensional cells.

  Let $U$ be a small open neighborhood of $v \in V(G_1)$.  Given
  $f|_U$, incrementally construct a CW-morphism $f'\colon G_1 \to G_2$ as
  follows.  Each edge $x'$ incident to $v$ can be expressed as $x' =
  \sigma_0^\alpha x$ for some "$0 \leq \alpha < $valence($v$)".  Let
  $w = f(v)$ and $y=f(x)$, and define:
  \begin{align*}
    f_1'(x) &:= y,
    \\
    f_0'(v) &:= w,
    \\
    f_1'(x') &:= \sigma_0^\alpha y  = \sigma_0^\alpha f_1'(x)
    \quad
    \text{if $x' = \sigma_0^\alpha x$},
    \\
    f_0'(v'_\alpha) &:= w'_\alpha,
  \end{align*}
  where:
  \begin{itemize}
  \item $v'_\alpha$ is the endpoint of $x' = \sigma_0^\alpha x$ ``opposite''
    to $v$,
  \item $w'_\alpha$ is the endpoint of $y' = \sigma_0^\alpha y$ ``opposite''
    to $w$.
  \end{itemize}
  Then $f'$ extends $f$ on an open set $U' \supsetneq \overline{U}$,
  which contains the subgraph formed by all edges attached to $v$ and $v'$.
  In addition:
  \begin{itemize}
  \item $f'_1(x') = f(x')$ up to a homotopy fixing the endpoints since
    $f$ commutes with $\sigma_0$,
  \item $f'_0(v_\alpha) = f(v_\alpha)$ since $f$ preserves adjacency.
  \end{itemize}

  By repeating the same construction about the vertices $v'_\alpha$
  and $w'_\alpha$, one can extend $f'$ to a CW-morphism that agrees
  with $f$ on an open set $U'' \supsetneq \overline{U'}$.  

  Recursively, by connectedness, we can thus extend $f'$
  to agree with $f$ (up to homotopy) over all of $G_1$.
\end{proof}

\subsection{Enumeration of \q{Fatgraph} isomorphisms}

The stage is now set for presenting the algorithm to enumerate the
isomorphisms between two given "Fatgraph" objects.  Pseudo-code is
listed in Algorithm~\ref{algo:iso}; as this procedure is quite
complex, a number of auxiliary functions have been used, whose purpose
is explained in Section~\ref{sec:iso-aux}.  Function \fn{isomorphisms}, given two
"Fatgraph" objects $\+G_1$ and $\+G_2$, returns a list of triples
"(pv, rot, pe)", each of which determines an isomorphism.  If there is
no isomorphism connecting the two graphs, then the empty list "[ ]" is
returned. 
\begin{Algorithm}
  \caption{Enumerate isomorphisms between two \q{Fatgraph} objects
    $\+G_1$ and~$\+G_2$: output of the algorithm is a list of triples
    \q{(pv, rot, pe)}.  If there is no isomorphism connecting the two
    input fatgraphs, the empty list is returned.}
  \label{algo:iso}
\begin{lstlisting}
def isomorphisms($\+G_1$, $\+G_2$):
  # immediately rule out easy cases of no isomorphisms
  if |graphs invariants differ|:                   |\n{iso-easy-1}|
    return [ ]                                     |\n{iso-easy-2}|
  result = [ ]
  vs1 = valence_spectrum($\+G_1$)
  vs2 = valence_spectrum($\+G_2$)
  (valence, vertices) = starting_vertices($\+G_2$)       |\n{iso:1}|
  $\+v_1$ = vs1[valence][0]                              |\n{iso:2}|
  for $\+v_2$ in compatible_vertices($\+v_1$, vertices):  |\n{iso-outer}|
    for rot in $0, \ldots, $valence:                      |\n{iso-inner}|
      |Initialize |pv, rots, pe| as empty maps|
      pv[$\+v_1$] = $\+v_2$
      rots[$\+v_1$] = rot
      extend_map(pe, $\+v_1$, rotated($\+v_2$, rot))
      if |extension failed|:
        |continue with next| rot
      # breadth-first search to extend the mapping over corresponding vertices
      nexts = neighbors(pv, pe, $\+G_1$, $\+v_1$, $\+G_2$, $\+v_2$) |\n{iso-nb1}|
      while size(pv) < $\+G_1$.num_vertices:
        neighborhood = [ ]
        for ($\+v_1'$, $\+v_2'$, $r$) in nexts:
          (pv, rots, pe) = extend_iso(pv, rots, pe, $\+G_1$, $\+v_1'$, $r$, $\+G_2$, $\+v_2'$)
          if |cannot extend|:
            |exit ``while'' loop and continue with next| rot
          |append| neighbors(pv, pe, $\+G_1$, $\+v_1'$, $\+G_2$, $\+v_2'$) |to| neighborhood |\n{iso-nb2}|
        nexts = neighborhood
      # isomorphism found, record it
      result.append((pv, rots, pe))                |\n{iso:end}|
  return result

\end{lstlisting}
\end{Algorithm}

By the rigidity lemma~\ref{lemma:rigidity}, any fatgraph isomorphism
is uniquely determined by the mapping of a small neighborhood of any
vertex.  The overall strategy of the algorithm is thus to pick a pair
of ``compatible'' vertices and try to extend the map as in the proof
of of lemma~\ref{lemma:rigidity}.  

We wish to stress the difference with
isomorphism of ordinary graphs: since an isomorphism $f$ is uniquely
determined by any pair of corresponding vertices, the initial choice
of candidates $v$, $f(v)$ either yields an isomorphism or it does not:
there is no backtracking involved.  

Since the isomorphism computation is implemented as an exhaustive
search, it is worth doing a few simple checks to rule out cases of
non-isomorphic graphs (lines \nr{iso-easy-1}--\nr{iso-easy-2}).  One
has to weigh the time taken to compute a graph invariant versus the
potential speedup obtained by not running the full scan of the search
space; experiments run using the Python code show that the
following simple invariants already provide some good speedup:
\begin{itemize}
\item the number of vertices, edges, boundary cycles;
\item the total number of loops;
\item the set of valences;
\item the number of vertices of every given valence.
\end{itemize}

Since an isomorphism is uniquely determined by its restriction to
\emph{any} vertex, one can restrict to considering just pairs of the
form $(\+v_1, \+v_2)$ where $\+v_1$ is a chosen vertex in $\+G_1$.
Then the algorithm tries all possible ways (rotations) of mapping
$\+v_1$ into a compatible vertex $\+v_2$ in $\+G_2$.  The body of the
inner loop (line~\nr{iso-inner} onwards) mimics the construction in
the proof of Lemma~\ref{lemma:rigidity}.

The starting vertex $\+v_1$ should be selected so to
minimize the number of mapping attempts performed; this is currently
done by minimizing the product of valence and number of vertices of
that valence on $\+G_2$ (line~\nr{iso:1}), and then picking a
vertex of the chosen valence in $\+G_1$ as $\+v_1$
(line~\nr{iso:2}).\footnote{The checks already performed ensure that
  \lstinline`$\+G_1$` and \lstinline`$\+G_2$` have the same ``valence
  spectrum'', so \lstinline`$\+G_1$` has at least one vertex of the
  chosen valence.}

First, given the target vertex $\+v_2$ and a rotation "rot", a new
triple "(pv,rots,pe)" is created; "pv" is set to represent the initial
mapping of $\+v_1$ onto $\+v_2$, rotated leftwards by "rot" positions,
and "pe" maps edges of $\+v_1$ into corresponding edges of the rotated
$\+v_2$.  If this mapping is not possible (e.g., $\+v_1$ has a loop
and $\+v_2$ does not, or not in a corresponding position), then the
attempt is aborted and execution continues from line~\nr{iso-inner}
with the next candidate "rot".

The mapping defined by "(pv,rots,pe)" is then extended to neighbors of
the vertices already inserted. This entails a breadth-first
search\footnote{The variables \q{nexts} and \q{neighborhood} play the
  role of the FIFO list in the usual formulation of breadth-first search:
  vertices are added to \q{neighborhood} during a loop, and the
  resulting list is then orderly browsed (as \q{nexts}) in the next
  iteration.} over pairs of corresponding vertices, starting from
$\+v_1$ and $\+v_2$.  Note that, in this extension step, not only 
the source and target vertices, but also the rotation to be
applied is uniquely determined: chosen a vertex $\+v_1'$ connected to
$\+v_1$ by an edge $\+e$, there is a unique rotation $r$ on $\+v_2'$
such that "pv[$\+e$]" has the same attachment index to $\+v_2'$ that
$\+e$ has to $\+v_1$.  If, at any stage, the extension of the current
triple "(pv, rots, pe)" fails, it is discarded and execution continues
from line~\nr{iso-inner} with the next value of "rot".

When the loop started at line~\nr{iso-outer} is over, execution
reaches the end of the "isomorphisms" function, and returns the
(possibly empty) list of isomorphisms to the caller.

\begin{theorem}
  Given "Fatgraph" objects $\+G_1$, $\+G_2$, function "isomorphisms"
  returns all "Fatgraph" isomorphisms from $\+G_1$ to $\+G_2$.
\end{theorem}
\begin{proof}
  Given an isomorphism $f : \+G_1 \to \+G_2$, restrict $f$ to the
  starting vertex~$\+v_1$: then $f$ will be output when
  Algorithm~\ref{algo:iso} examines the pair $\+v_1$, $f(\+v_1)$;
  since Algorithm~\ref{algo:iso} performs an exhaustive search, $f$
  will not be missed.

  Conversely, since equation~\eqref{eq:fg-iso} holds by construction
  for all the mappings returned by "isomorphisms", then each returned
  triple $\+f=$"(pv, rots, pe)" is an isomorphism.
\end{proof}

\subsubsection{Auxiliary functions}
\label{sec:iso-aux}

Here is a brief description of the auxiliary functions used in the
listing of Algorithm~\ref{algo:iso} and~\ref{algo:neighbors}.  Apart
from the \fn{neighbors} function, they are all straightforward to
implement, so only a short specification of the behavior is given, with no
accompanying pseudo-code.

\paragraph{The \q{neighbors} function.}
\begin{definition}
  Define a candidate extension as a triplet
  "($\+v_1'$, $\+v_2'$, $r$)", where:
  \begin{itemize}
  \item $\+v_1'$ is a vertex in $\+G_1$, connected to
    $\+v_1$ by an edge $\+e$;
  \item $\+v_2'$ is a vertex in $\+G_2$, connected to $\+v_2$ by
    edge "$\+e'$ == pe[$\+e$]";
  \item $r$ is the rotation to be applied to $\+v_2'$ so that edge
    $\+e$ and $\+e'$ have the same attachment index, i.e., they are
    incident at corresponding positions in $\+v_1'$ and $\+v_2'$.
  \end{itemize}
\end{definition}

Function "neighbors" lists candidate extensions that extend map
"pv" in the neighborhood of given input vertices $\+v_1$ (in the
domain fatgraph $\+G_1$) and $\+v_2$ (in the image fatgraph
$\+G_2$).  It outputs a list of triplets $(\+v_1', \+v_2', r)$, each
representing a candidate extension.
\begin{Algorithm}
  \caption{Enumerate the candidate extensions of the given
    \q{pv} and \q{pe} in the neighborhood of input vertices $\+v_1$
    and $\+v_2$.}
  \label{algo:neighbors}
\begin{lstlisting}
    def neighbors(pv, pe, $\+G_1$, $\+v_1$, $\+G_2$, $\+v_2$):
      result = [ ]
      for |each non-loop edge $\+e$ attached to $\+v_1$|:
        |let $(\+v_1', a_1)$ be the endpoint of $\+e$ distinct from $\+v_1$|
        if $\+v_1'$ |already in| pv |domain|:
          |continue with next $\+e$|
        |let $(\+v_2', a_2)$ be the endpoint of $\+e'$=pe[$\+e$] distinct from $\+v_2$|
        if $\+v_2'$ |already in| pv |image|:
          |continue with next $\+e$|
        result.append($(\+v_1', \+v_2', a_1-a_2)$)
      return result
      
\end{lstlisting}
\end{Algorithm}

A sketch of this routine is given in Algorithm~\ref{algo:neighbors}.  
Two points are worth of notice:
\begin{enumerate}
\item By the time \fn{neighbors} is called (at lines~\nr{iso-nb1}
  and~\nr{iso-nb2} in Algorithm~\ref{algo:iso}), the map "pe" has already
  been extended over all edges incident to $\+v_1$, so we can
  safely set $\+e'=$"pe[$\+e$]" in \fn{neighbors}.
\item Algorithm~\ref{algo:iso} only uses \fn{neighbors} with the
  purpose of extending "pv" and "pe", so \fn{neighbors} ignores
  vertices that are already in the domain or image of "pv".
\end{enumerate}

\paragraph{The \q{valence\_spectrum} function.}
The auxiliary function \fn{valence\_spectrum}, given a "Fatgraph" instance
$\+G$, returns a mapping that associates to each valence $z$ the list
$V_z$ of vertices of $\+G$ with valence $z$.

\paragraph{The \q{starting\_vertices} function.}
For each pair $(z, V_z)$ in the valence spectrum, define its \emph{intensity}
as the product $z \cdot \card{V_z}$ (valence times the number of vertices
with that valence).  The function \fn{starting\_vertices} takes as input a
"Fatgraph" object $\+G$ and returns the pair $(z, V_z)$ from the valence
spectrum that minimizes intensity.  In case of ties, the pair with the
largest $z$ is chosen.

\paragraph{The \q{compatible} and \q{compatible\_vertices} functions.}
Function \fn{compatible} takes a pair of vertices $\+v_1$ and $\+v_2$ as
input, and returns boolean "True" iff $\+v_1$ and $\+v_2$ have the
same invariants.  (This is used as a short-cut test
to abandon a candidate mapping before trying a full adjacency list
extension, which is computationally more expensive.)  The sample code
uses valence and number of loops as invariants.

The function \fn{compatible\_vertices} takes a vertex $\+v$ and a list of
vertices $L$, and returns the list of vertices in $L$ that are
compatible with $\+v$ (i.e., those which $\+v$ could be mapped to).

\paragraph{The \q{extend\_map} and \q{extend\_iso} functions.}
The \fn{extend\_map} function takes as input a mapping "pe" and a pair of
\emph{ciliated} vertices $\+v_1$ and $\+v_2$, and alters "pe" to map
edges of $\+v_1$ to corresponding edges of $\+v_2$: the cilium to the
cilium, and so on: "pe[$\sigma_0^\alpha(\+e)$] == $\sigma_0^\alpha$(pe[$\+e$])".
If this extension is not possible, an error is signaled to the caller.

The \fn{extend\_iso} function is passed a "(pv, rots, pe)" triplet, a
vertex $\+v_1'$ of $\+G_1$, a vertex $\+v_2'$ of $\+G_2$ and a
rotation $r$; it alters the given "(pv,rots,pe)" triple by adding a
mapping of the vertex $\+v_1'$ into vertex $\+v_2'$ (and rotating the
target vertex by "r" places rightwards).  If the extension is
successful, it returns the extended map "(pv, rot, pe)"; otherwise,
signals an error.

\subsection{Operations with \q{Fatgraph} Isomorphisms}
\label{sec:iso}

\paragraph{Compare pull-back orientation.}
The \fn{compare\_orientations} function takes an isomorphism triple 
"(pv, rots, pe)" and a pair of "Fatgraph" objects $\+G_1$ , $\+G_2$, 
and returns $+1$ or $-1$ depending on whether the orientations of the
target "Fatgraph" pulls back to the orientation of the source
"Fatgraph" via the given isomorphism.

Recall that for a "Fatgraph" object $\+G$, the orientation is
represented by a mapping "$\+G$.orient" that associates an edge $\+e$ with
its position in the wedge product that represents the
orientation; therefore, the pull-back orientation according to an
isomorphism "(pv, rots, pe)" from $\+G$ to $\+G'$ is simply given by the map 
"$\+e$ $\mapsto$ $\+G'$.orient[pe[$\+e$]]".
Thus, the comparison is done by constructing the permutation that maps
"$\+G$.orient[$\+e$]" to "$\+G'$.orient[pe[$\+e$]]" and taking its
sign (which has linear complexity with respect to the number of edges).

\paragraph{The \q{is\_orientation\_reversing} function.}
Determining whether an automorphism reverses orientation is crucial
for knowing which fatgraphs are orientable.
Function \fn{is\_orientation\_reversing} takes a "Fatgraph" object and
an isomorphism triple "(pv, rots, pe)" as input, and returns boolean
"True" iff the isomorphism reverses orientation.  This amounts to
checking whether the given orientation and that of the pull-back one agree,
which can be done with the comparison method discussed above.

\paragraph{Transforming boundary cycles under an isomorphism.}
The function \fn{transform\_boundary\_cycle} is used when comparing
\emph{marked} fatgraphs: as the marking is a function on the
boundary cycles, we need to know exactly which boundary cycle of the
target graph corresponds to a given boundary cycle in the source
graph.  

Recall that "BoundaryCycle" instances are defined as list of
\emph{corners}; function \fn{transform\_boundary\_cycle} takes a
"BoundaryCycle" $\+b$ and returns a new "BoundaryCycle" object $\+b'$,
obtained by transforming each corner according to a graph isomorphism.
Indeed, \fn{transform\_boundary\_cycle} is straightforward loop over the
corners making up $\+b$: For each corner
"($\+v$,i,j)", a new one is constructed by transforming the vertex
according to map "pv", and displacing indices "i" and "j" by the
rotation amount indicated by "rot[\+v]" (modulo the number of edges
attached to $\+v$).

\section{Generation of fatgraphs}
\label{sec:generation}

Let \fn{MgnGraphs} be the function which, given two integers $g$,
$n$ as input, returns the collection of $\R_{g,n}$ graphs.  Let us
further stipulate that the output result will be represented as a list
$R$: the $0$-th item in this list is the list of graphs with the
maximal number $m$ of edges; the $k$-th item $R[k]$ is the list of
graphs having $m - k$ edges.  There are algorithmic advantages in this
subdivision, which are explained below.

Graphs with the maximal number of edges are trivalent graphs; they are
computed by a separate function \fn{MgnTrivalentGraphs}, described in
Section~\ref{sec:trivalent}.

We can then proceed to generate all graphs in $\R_{g,n}$ by
contraction of regular edges: through contracting one edge in
trivalent graphs we get the list $R[1]$ of all graphs with
$m-1$ edges; contracting one edge of $G \in R[1]$, we
get $G' \in R[2]$ with $m-2$ edges, and so on.
\begin{Algorithm}
  \caption{Function \fn{MgnGraphs} returns all connected fatgraphs
    having prescribed genus $g$ and number of boundary cycles $n$.
    Actual output of the function is a list $R$, whose $k$-th element
    $R[k]$ is itself a list of graphs in $\R_{g,n}$ with $m-k$ edges.}
  \label{algo:MgnGraphs}
  \begin{lstlisting}
def MgnGraphs($g$,$n$):
  $m$ = $4g + 2n - 5$                        # maximum number of edges
  $R$ = |array of |$m$| empty lists|
  $R[0]$ = MgnTrivalentGraphs($g$,$n$)     # first item contains all 3-valent graphs
  for $k$ in $1$, \ldots, $m-1$:                     |\n{mgngraphs:1}|
    |Initialize $R[k]$ as an empty list|
    for $\+G$ in $R[k-1]$:                            |\n{mgngraphs:3}|
      for $\+e$ in edge_orbits($\+G$):                |\n{mgngraphs:4}|
        if |$\+e$ is a loop|:                         |\n{mgngraphs:5}|
          |continue with next $\+e$|
        $\+G'$ = contract($\+G$, $\+e$)               |\n{mgngraphs:6}|
        if |$\+G'$ not already in $R[k]$|:            |\n{mgngraphs:7}|
          |append $\+G'$ to $R[k]$|                   |\n{mgngraphs:8}|
  return $R$
\end{lstlisting}
\end{Algorithm}
Pseudo-code for \fn{MgnGraphs} is shown in
Algorithm~\ref{algo:MgnGraphs}.  The loop at
lines~\nr{mgngraphs:4}--\nr{mgngraphs:8} is the core of the function:
contract edges of the fatgraph $G$ (with $m-k+1$ edges) to generate
new fatgraphs with $m-k$ edges.  However, we need not contract every
edge of a fatgraph: if $a \in \Aut G$ is an automorphism and $e \in
E(G)$ is an edge, then the contracted graphs $G' = G/e$ and $G'' =
G/a(e)$ are isomorphic.  Hence, we can restrict the computation to
only one representative edge per orbit of the action induced by $\Aut
G$ on the set $E(G)$; the "edge_orbits" function referenced at
line~\nr{mgngraphs:4} should return a list of representative edges,
one per each orbit of $\Aut(G)$ on $E(G)$.

Lines~\nr{mgngraphs:7}--\nr{mgngraphs:8} add $\+G'$ to $R[k]$ \emph{only
  if it is not already there}.  This is the most computationally
expensive part of the \fn{MgnGraphs} function: we need to perform a
comparison between $\+G'$ and each element in $R[k]$; testing equality
of two fatgraphs requires computing if there are isomorphisms
between the two, which can only be done by attempting enumeration of
such isomorphisms.  Fatgraph isomorphism is discussed in detail in
Section~\ref{sec:isomorphism}.

If $N_k$ is the number of elements in $R[k]$ and $T_\text{iso}$ is the
average time needed to determine if two graphs are isomorphic, then
evaluating whether $\+G'$ is already contained in $R[k]$ takes $O(N_k
\cdot T_\text{iso})$ time: thus, the subdivision of the output $R$
into lists, each one holding graphs with a specific number of edges,
reduces the number of fatgraph comparisons done in the innermost loop
of \fn{MgnGraphs}, resulting in a substantial shortening of the total
running time.

Note that the top-level function \fn{MgnGraphs} is quite independent of
the actual implementation of the "Fatgraph" type of objects: all is
needed here, is that we have methods for enumerating edges of a
"Fatgraph" object, contracting an edge, and testing two graphs for
isomorphism.

\begin{lemma}
  If \fn{MgnTrivalentGraphs}$(g,n)$ returns the complete list of
  \emph{trivalent} fatgraphs in $\R_{g,n}$, then the function
  \fn{MgnGraphs} defined above returns the complete set of fatgraphs
  $\R_{g,n}$.
\end{lemma}
\begin{proof}
  By the above dissection of the algorithm, all we need to prove is
  that any fatgraph in $\R_{g,n}$ can be obtained by a chain of edge
  contractions from a trivalent fatgraph.  This follows immediately
  from the fact that any fatgraph vertex $v$ of valence $z \geq 3$ can
  be expanded (in several ways) into vertices $v_1$, $v_2$ of valences
  $z_1$, $z_2$ such that $z = (z_1 -1) + (z_2 -1)$, plus a connecting
  edge.
\end{proof}

\subsection{Generation of Trivalent Fatgraphs}
\label{sec:trivalent}

Generation of trivalent graphs can be tackled by an inductive
procedure: given a trivalent graph, a new edge is added, which joins
the midpoints of two existing edges.  
In order to determine which graphs should be input to this ``edge
addition'' procedure, one can follow the reverse route, and
ascertain how a trivalent graph is transformed by \emph{deletion} of
an edge.

Throughout this section, $l$ and $m$ stand for the number of vertices
and edges of a graph; it will be clear from the context, which
exact graph they are invariants of.

\subsubsection{Removal of edges}
\label{sec:removal}

Let $G \in \R_{g,n}$ be a \emph{connected} trivalent graph. Each edge
$x \in E(G)$ falls into one of the following categories:
\begin{enumerate}[\slshape A)]
\item $x$ is a loop: both endpoints of $x$ are
  attached to a single vertex $v$, another edge $x'$ joins $v$ with a
  distinct vertex $v'$;
\item $x$ joins two distinct vertices $v, v'
  \in V(G)$ and separates two distinct boundary cycles $\beta, \beta'
  \in B(G)$;
\item $x$ joins two distinct vertices $v, v'
  \in V(G)$ but belongs to only one boundary cycle $\beta \in B(G)$,
  within which it occurs twice (once for each orientation).
\end{enumerate}
Deletion of edge $x$ requires different adjustments in order to get a
trivalent graph again in each of the three cases above; it also yields
a different result in each case.
\begin{figure}
  \centering
  \includegraphics{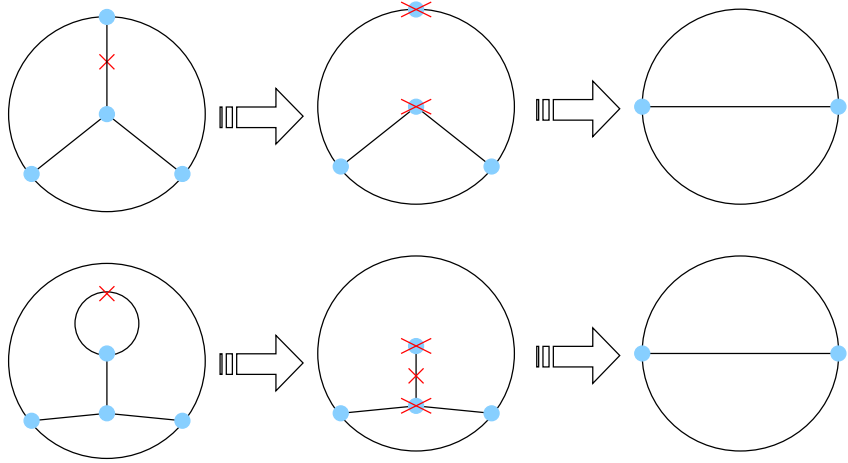}
  \caption{Graphical illustration of fatgraph edge removal.  \emph{Top
      row:} a regular edge (crossed) is removed from an $\R_{0,4}$
    graph; its endpoints are further removed; the remaining edges are
    joined and the resulting graph is a trivalent fatgraph in
    $\R_{0,3}$.  \emph{Bottom row:} a loop is removed from a trivalent
    $\R_{0,4}$ graph; the stem together with its endpoints has to be
    removed as well; the remaining edges are joined, and we end up
    with a trivalent fatgraph in $\R_{0,3}$.}
  \label{fig:removal}
\end{figure}

Case {\slshape A)}: If $x$ is a loop attached to $v$, then, after deletion
of $x$, one needs to also delete the loose edge $x'$ and the
vertex $v'$ (that is, join the two other edges attached to $v'$; see
Figure~\ref{fig:removal}, bottom row).  The resulting fatgraph $G'$ has:
\begin{itemize}
\item two vertices less than $G$: $v$ and $v'$ have been deleted;
\item three edges less: $x$, $x'$ have been deleted and two other
  edges merged into one;
\item one boundary cycle less: the boundary cycle totally bounded by
  $x$ has been removed.
\end{itemize}
Therefore:
\begin{align*}
  2 - 2g' &= \chi(G') = l' - m' + n' 
  \\
  &= (l-2) -(m-3) + (n-1)
  \\
  &= l - m + n = \chi(G) = 2 - 2g,
\end{align*}
hence $g=g'$, and 
\begin{equation}
G' \in \R_{g,n-1}.
\label{eq:A}\tag{A}
\end{equation}

In case {\slshape B)}, $x$ joins distinct vertices $v$,
$v'$ and separates distinct boundary cycles (see Figure~\ref{fig:removal},
top row).  Delete $x$ and merge the two edges attached to each of
the two vertices $v$ and $v'$; in the process, the two boundary cycles
$\beta, \beta'$ also merge into one. The resulting fatgraph $G'$ is
connected. Indeed, given any two vertices $u, u' \in V(G')$, there is
a path $(x_1, \ldots, x_k)$ connecting $u$ with $u'$ in $G$.  If this
path passes through $x$, one can replace the occurrence of $x$ with
the perimeter ---excluding $x$--- of one of the two boundary cycles
$\beta, \beta'$ to get a path joining $v$ and $v'$ which avoids $x$,
and thus projects to a path in $G'$.  Again we see that $G'$ has:
\begin{itemize}
\item two vertices less than $G$: $v$ and $v'$ have been deleted;
\item three edges less: $x$ has been deleted and four other
  edges merged into two, pair by pair;
\item one boundary cycle less: the boundary cycles $\beta$, $\beta'$
  have been merged into one.
\end{itemize}
Therefore $g=g'$, and 
\begin{equation}
G' \in \R_{g,n-1}.
\label{eq:B}\tag{B}
\end{equation}

In case {\slshape C)}, $x$ joins distinct vertices $v$,
$v'$ but belongs into \emph{one} boundary cycle $\beta \in B(G)$ only.
Delete edge $x$ and the two vertices $v$, $v'$, joining the attached
edges two by two as in case {\slshape B)}.  We distinguish two
cases, depending on whether the resulting fatgraph is connected.
\begin{enumerate}
\item[\slshape C')] If the resulting fatgraph $G'$ is connected, then
  $\beta \in B(G)$ has been split into two distinct boundary cycles
  $\beta', \beta'' \in B(G')$.  Indeed, write the boundary cycle
  $\beta$ as an ordered sequence of \emph{oriented} edges: $y_0 \to y_1 \to
  \ldots \to y_k \to y_0$.  Assume the $y_*$ appear in this sequence in the
  exact order they are encountered when walking along $\beta$ in the
  sense given by the fatgraph orientation. The oriented edges $y_i$
  are pairwise distinct: if $y_i$ and $y_j$ share the same supporting
  edge, then $y_i$ and $y_j$ have opposite orientations. By the
  initial assumption of case~{\slshape C)}, edge $x$
  must appear \emph{twice} in the list: if $\bar x$ and $\underline x$
  denote the two orientations of $x$, then $y_i = \bar x$ and $y_j =
  \underline x$.  Deleting $x$ from $\beta$ is (from a homotopy point
  of view) the same as replacing $y_i = \bar x$ with $\bar x \to
  \underline x$, and $y_j = \underline x$ with $\underline x \to \bar
  x$ when walking a boundary cycle. Then we see that $\beta$ splits
  into two disjoint cycles:
\begin{align*}
  \beta' &= y_0 \to y_1 \to \cdots \to y_{i-1} \to \bar{x} \to
  \underline{x} \to y_{j+1} \to \cdots \to y_k \to y_0,
  \\
  \beta'' &= y_{i+1} \to \cdots \to y_{j-1} \to \underline{x} \to
  \bar{x} \to y_{i+1}.
\end{align*}
In this case, $G'$ has:
\begin{itemize}
\item two vertices less than $G$: $v$ and $v'$ have been deleted;
\item three edges less: $x$ has been deleted and four other
  edges merged into two, pair by pair;
\item one boundary cycle \emph{more}: the boundary cycle $\beta$ has
  been split in the pair $\beta'$, $\beta''$.
\end{itemize}
Therefore $g'=g-1$ and $n'=n+1$, so:
\begin{equation}
G' \in \R_{g-1,n+1}.
\label{eq:C'}\tag{C'}
\end{equation}

\item[\slshape C'')] $G'$ is a disconnected union of fatgraphs $G'_1$ and
  $G'_2$; for this statement to hold unconditionally, we temporarily allow a
  single circle into the set of connected fatgraphs (consider it a fatgraph
  with one closed edge and no vertices) as the one and only element of
  $\R_{0,2}$.  As will be shown in Lemma~\ref{lemma:no-c2}, this is
  irrelevant for the \fn{MgnTrivalentGraphs} algorithm.  Now:
  \begin{equation*} 
    l'_1 + l'_2 = l - 2, 
    \qquad m'_1 + m'_2 = m - 3,
    \qquad n'_1 + n'_2 = n + 1,
  \end{equation*} hence:
\begin{align*} 
  (2- 2g'_1) + (2-2g'_2) &= (l-2) - (m-3) + (n+1) 
  \\ 
  &= (l-m+n) + 2 = 4 - 2g
\end{align*} 
So that $g'_1 + g'_2 = g + 2$, $n'_1 + n'_2 = n+1$, and:\footnote{%
  Here we use $\otimes$ to indicate juxtaposition of graphs: $G_1
  \otimes G_2$ is the (non-connected) fatgraph having two connected
  components $G_1$ and $G_2$.
}%
\begin{equation} 
  G'= G'_1 \otimes G'_2 \in \R_{g'_1, n'_1} \otimes \R_{g'_2, n'_2}.
  \label{eq:C''}\tag{C''}
\end{equation}
\end{enumerate}

\subsubsection{Inverse construction}
\label{sec:addition}

If $x \in E(G)$ is an edge of a fatgraph $G$, denote $\bar{x}$ and
$\underline{x}$ the two opposite orientations of $x$.

In the following, let $\R'_{g,n}$ be the set of fatgraphs with a
selected oriented edge:
\begin{equation*}
  \R'_{g,n} := \{ (G,\bar x) : G \in \R_{g,n}, \bar{x} \in L(G)\}.
\end{equation*}
Similarly, let $\R''_{g,n}$ be the set of fatgraphs with two
chosen oriented edges:
\begin{equation*}
  \R''_{g,n} := \{ (G, \bar x, \bar y) : G \in \R_{g,n}, 
                   \bar{x}, \bar{y} \in L(G) \}.
\end{equation*}
The following abbreviations are convenient:
\begin{equation*}
\R = \cup \R_{g,n}, 
\qquad
\R' = \cup \R'_{g,n}, 
\qquad
\R'' = \cup \R''_{g,n}.
\end{equation*}

\begin{figure}
  \centering
  \includegraphics{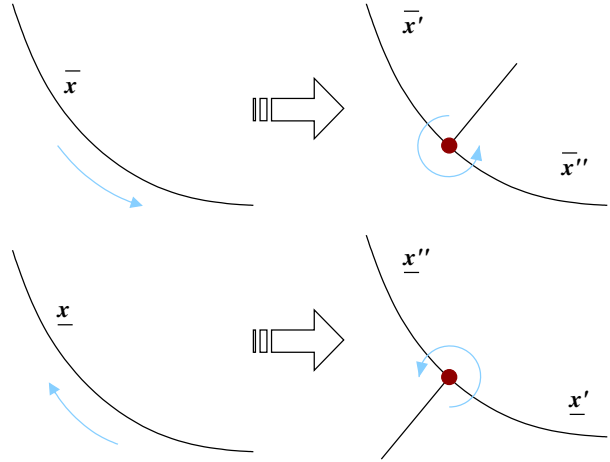}
  \caption{When adding a new vertex in the middle of an edge $x$, the cyclic order depends on the oriented edge: the two orientations $\bar{x}$ and $\underline{x}$ get two inequivalent cyclic orders.}
  \label{fig:adding}
\end{figure}

Define the attachment of a new edge to a fatgraph in the following
way.  Given a fatgraph $G$ and an \emph{oriented} edge $\bar{x}$, we
can create a new trivalent vertex $v$ in the midpoint of~$x$, and
attach a new edge to it, in such a way that the two halves of~$x$
appear, in the cyclic order at $v$, in the same order induced by the
orientation of $\bar{x}$.  Figure~\ref{fig:adding} depicts the process.

We can now define maps that invert the constructions {\slshape A)}, {\slshape B)},
{\slshape C')} and {\slshape C'')} defined in the previous section.
\begin{sidewaysfigure}
  \centering
  \includegraphics[width=\textwidth]{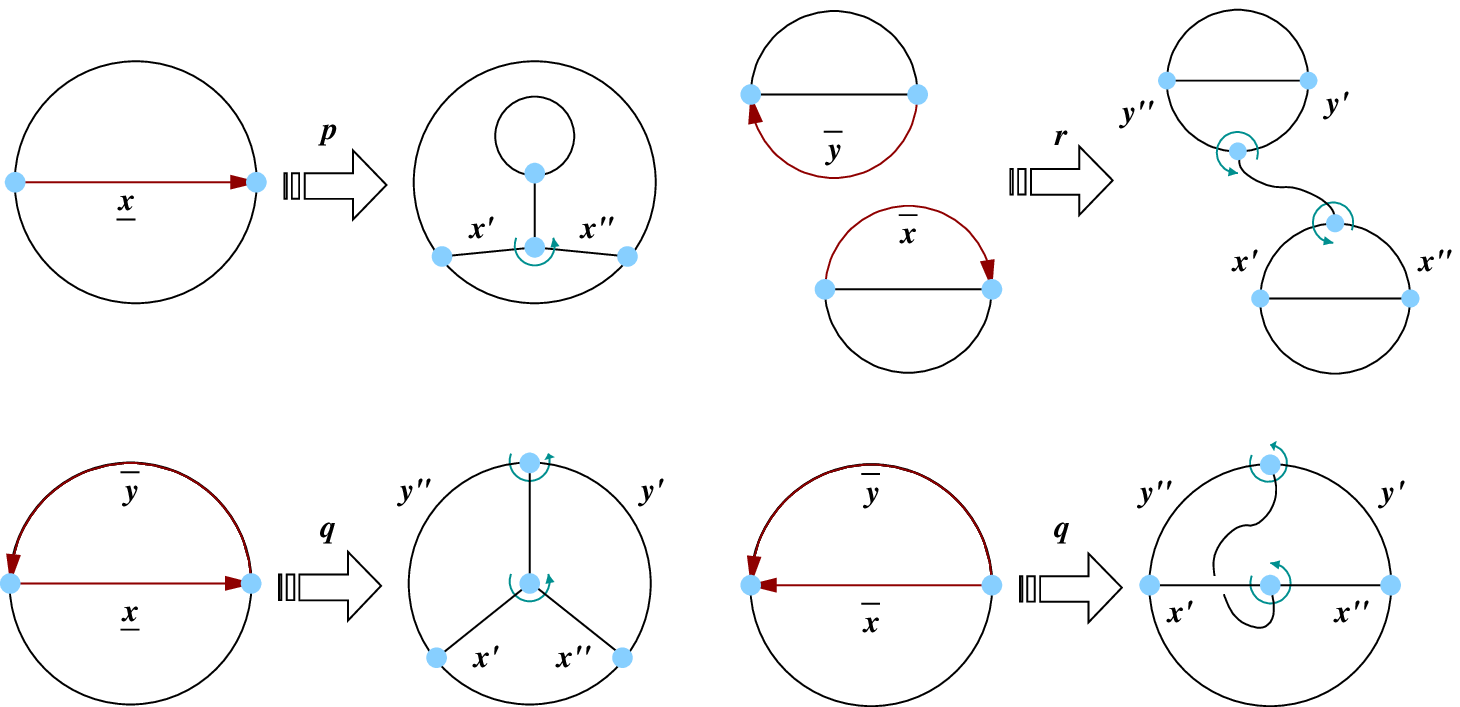}
  \caption{Graphical illustration of maps $p$, $q$, $r_{g,n}$.  Top
    left: $p(G,\bar{x})$ attaches a ``slipknot'' to edge $\bar{x}$.
    Top right: $r_{2,5}(G_1, \bar{x}, G_2, \bar{y})$ joins fatgraphs
    $G_1$ and $G_2$ with a new edge. Bottom: $q(G,\bar{x}, \bar{y})$
    (left) and $p(G, \underline{x}, \bar{y})$ (right); it is shown how
    changing the orientation of an edge can lead to different
    results.}
  \label{fig:pqr}
\end{sidewaysfigure}

Let $p_{g,n} : \R'_{g,n-1} \to \R_{g,n}$ be the map that creates a
fatgraph $p(G,\bar{x})$ from a pair $(G, \bar{x})$ by attaching the
loose end of a ``slip knot''\footnote{A single 3-valent vertex with
  one loop attached and a regular edge with one loose end.} to the
midpoint of~$x$. The map $p: \R' \to \R$ defined by $p|_{\R'_{g,n}} := p_{g,n}$
is ostensibly inverse to {\slshape A)}.

To invert {\slshape B)} and {\slshape C')}, define a map $q : \R'' \to \R$ that
operates as follows:
\begin{itemize}
\item Given $(G, \bar{x}, \bar{y})$ with $\bar{x} \not= \bar{y}$, the
  map $q$ attaches a new edge to the midpoints of $x$ and $y$; again
  the cyclic order on the new midpoint vertices is chosen such that
  the two halves of $x$ and $y$ appear in the order induced by the
  orientations $\bar x$, $\bar y$.
\item When $\bar{x} = \bar{y}$, let us further stipulate that the
  construction of $q(G, \bar{x}, \bar{x})$ happens in two steps:
  \begin{enumerate}
  \item a new trivalent vertex is created in the midpoint of $x \in
    E(G)$ and a new edge $\xi$ is attached to it,
  \item create a new trivalent vertex in the middle of the half-edge
    which comes first in the ordering induced by the
    orientation $\bar{x}$; attach the loose end of the new edge $\xi$
    to this new vertex.
  \end{enumerate}
  It is clear that the above steps give an unambiguous definition of
  $q$ in all cases where $\bar{x}$ and $\bar{y}$ are orientations of
  the same edge of $G$, that is, $(G, \bar{x}, \bar{x})$, $(G,
  \bar{x}, \underline{x})$, $(G, \underline{x}, \bar{x})$, and $(G,
  \underline{x}, \underline{x})$.
\end{itemize}
Ostensibly, $q$ inverts the edge removal in cases {\slshape B)} and
{\slshape C')}: the former applies when a graph $G \in \R_{g,n}$ is
sent to $q(G) \in \R_{g,n+1}$, the latter when $G \in \R_{g,n}$ is
sent to $q(G) \in \R_{g+1,n-1}$.

Finally, to invert {\slshape C'')}, let us define
\begin{equation*}
r_{g,n}: \bigoplus_{\substack{g'_1+g'_2=g+2 \\ n'_1+n'_2=n}} 
\R'_{g'_1, n'_1} \times \R'_{g'_2, n'_2} \to \R.
\end{equation*}
From $(G', \bar{x}', G'',
\bar{x}'')$, construct a new fatgraph by bridging $G'$ and $G''$ with
a new edge, whose endpoints are in the midpoints of $x'$ and $x''$;
again, stipulate that the cyclic order on the new vertices is chosen
such that the two halves of $x'$, $x''$ appear in the order induced by
the orientations $\bar{x}'$, $\bar{x}''$.

Summing up, any fatgraph $G \in \R_{g,n}$ belongs to the image of one
of the above maps $p$, $q$, and $r$. There is considerable
overlap among the different image sets: in fact, one can prove that
$r$ is superfluous.
\begin{lemma}\label{lemma:no-c2}
  Any fatgraph obtained by inverting construction {\slshape C'')} lies in
  the image of maps $p$ and $q$.
\end{lemma}
\begin{proof}
  Assume, on the contrarily, that $G$ lies in the image of $r$ only.
  Then, deletion of any edge $x$ from $G$ yields a disconnected graph
  $G' \otimes G''$.  Both subgraphs $G'$ and $G''$ enjoy the same
  property, namely, that deletion of any edge disconnects: otherwise,
  if the removal of $y \in E(G')$ does not disconnect $G'$, then
  neither does it disconnect $G = r_{g,n}(G', G'')$, contrary to the
  initial assumption. As long as $G'$ or $G''$ has more than 3 edges, we can
  delete another edge; by recursively repeating the process, we end up
  with a fatgraph $G^*$ with $l^* \leq 3$ edges, which is again
  disconnected by removal of any edge.  Since $G^*$ is trivalent, $3
  \cdot m^* = 2 \cdot l^*$, therefore $G^*$ must have exactly 3 edges
  and 2 vertices. But all such fatgraphs belong to $\R_{0,3}$ or
  $\R_{1,1}$, and it is readily checked that there is no way to add an
  edge such that the required property holds, that any deletion
  disconnects.
\end{proof}

\subsubsection{The \fn{MgnTrivalentGraphs} algorithm}
\label{sec:MgnTrivalentGraphs}

The stage is now set for implementing the recursive generation of
trivalent graphs. Pseudo-code is listed in Algorithm~\ref{algo:MgnTrivalentGraphs}.
\begin{Algorithm}
  \caption{Return a list of all connected trivalent fatgraphs with
    prescribed genus $g$ and number of boundary cycles $n$.  A
    fatgraph if ``admissible'' iff it has the prescribed genus $g$ and
    number of boundary cycles $n$.}
  \label{algo:MgnTrivalentGraphs}
\begin{lstlisting}
def MgnTrivalentGraphs($g$, $n$):
  # avoid infinite recursion in later statements
  if $n$==0 or $(g,n)<(0,3)$:
    return |empty list|

  # Induction base: $\M_{0,3}$ and $\M_{1,1}$
  if $(g,n)$ == $(0,3)$:
    return |list of fatgraphs in $\R_{0,3}$|
  elif $(g,n)$ == $(1,1)$:
    return |list of fatgraphs in $\R_{1,1}$|

  # general case
  else:                                          |\n{mgntriv:0}|
    R = |empty list|

    # case A): hang a circle to all edges of graphs in $M_{g,n-1}$
    for $G$ in MgnTrivalentGraphs($g$, $n-1$):
      for $x$ in edge_orbits($G$):
        |add $p(G, \bar{x})$ to |R| if admissible|
        |add $p(G, \underline{x})$ to |R| if admissible|

    # case B): bridge all edges of a single graph in $M_{g,n-1}$
    for $G$ in MgnTrivalentGraphs($g$, $n-1$):
      for $(x,y)$ in edge_pair_orbits($G$):
        |add $q(G, \underline{x}, \underline{y})$ to |R| if admissible|
        |add $q(G, \underline{x}, \bar{y})$) to |R| if admissible|
        |add $q(G, \bar{x}, \underline{y})$ to |R| if admissible|
        |add $q(G, \bar{x}, \bar{y})$ to |R| if admissible|

    # case C'): bridge all edges of a single graph in $M_{g-1,n+1}$
    for $G$ in MgnTrivalentGraphs($g-1$, $n+1$):
      for $(x,y)$ in edge_pair_orbits($G$):
        |add $q(G, \underline{x}, \underline{y})$ to |R| if admissible|
        |add $q(G, \underline{x}, \bar{y})$ to |R| if admissible|
        |add $q(G, \bar{x}, \underline{y})$ to |R| if admissible|
        |add $q(G, \bar{x}, \bar{y})$ to |R| if admissible| |\n{mgntriv:last}|

    |remove isomorphs from R|
    return R
\end{lstlisting}
\end{Algorithm}

\begin{lemma}
  \fn{MgnTrivalentGraphs}$(g,n)$ generates all trivalent fatgraphs for
  each given $g$, $n$.  Only one representative per isomorphism class
  is returned.
\end{lemma}
\begin{proof}
  The function call \fn{MgnTrivalentGraphs}$(g, n)$ recursively calls
  itself to enumerate trivalent graphs of $\R_{g,n-1}$ and
  $\R_{g-1,n+1}$.  In particular, \fn{MgnTrivalentGraphs} must:
  \begin{itemize}
  \item provide the full set of fatgraphs $\R_{0,3}$ and $\R_{1,1}$ as
    induction base.
  \item return the empty set when called with an invalid $(g,n)$ pair;
  \end{itemize}
  The general case is then quite straightforward:
  \begin{inparaenum}
  \item apply maps $p$, $q$ to every fatgraph in $\R_{g,n-1}$, and $q$
    to every fatgraph in $\R_{g-1,n+1}$;
  \item discard all graphs that do not belong to $\R_{g,n}$; and
  \item take only one graph per isomorphism class into the result set.
  \end{inparaenum}

  To invert construction {\slshape A)}, map $p$ is applied to all
  fatgraphs $G \in \R_{g,n-1}$; if $a \in \Aut G$, then $p(a(G), a(x))
  = p(G, x)$, therefore we can limit ourselves to one pair $(G,x)$ per
  orbit of the automorphism group, saving a few computational cycles.
  Similarly, since $q$ is a function of $(G, \bar{x}, \bar{y})$, which
  is by construction invariant under $\Aut G$, we can again restrict
  to considering only one $(G, \bar{x}, \bar{y})$ per $\Aut G$-orbit;
  this is computed by the "edge_pair_orbits($G$)" function.
\end{proof}
Note that there is no way to tell from $G$ if fatgraphs $p(G,x)$ and
$q(G,x,y)$ belong to $\R_{g,n}$: one needs to check $g$ and $n$ before
adding the resulting fatgraph to the result set $R$.

The selection of only one representative fatgraph per isomorphism
class can be done by removing duplicates from the collection of
generated graphs in the end, or by running the isomorphism test before
adding each graph to the working list $R$.  The computational
complexity is quadratic in the number of generated graphs in both
cases, but the latter option requires less memory.  In any case, this
isomorphism test is the most computationally intensive part of
\fn{MgnTrivalentGraphs}.

For an expanded discussion of the size of the result set $R$, and a
comparison with other generation algorithms, see Appendix~\ref{sec:N}.
It would be interesting to re-implement the trivalent generation
algorithm using the technique outlined in \cite{mckay:isomorph-free},
and compare it with the current (rather naive) algorithm.

\subsubsection{Implementing maps $p(G,x)$ and $q(G,x,y)$}
\label{sec:bridge}

Implementation of both functions is straightforward
and pseudo-code is therefore omitted;\footnote{The interested reader is
referred to  the publicly-available code at
\url{http://fatghol.googlecode.com} for details.} the only question is
how to represent the ``oriented edges'' that appear in the signature
of maps $p$ and $q$.

In both $p$ and $q$, the oriented edge $\bar{x}$ or $\underline{x}$ is
used to determine how to attach a new edge to the midpoint of the
target (unoriented) edge $x$.  We can thus represent an oriented edge
as a pair $(\+e, s)$ formed by a "Fatgraph" edge $\+e$ and a ``side''
$s$: valid values for $s$ are $+1$ and $-1$, interpreted as follows.
The parameter~$s$ controls which of the two inequivalent cyclic orders
the new trivalent vertex will be given.  Let $\+a$, $\+b$, $\+c$ be
the edges attached to the new vertex in the middle of $\+e$, where
$\+a$,$\+b$ are the two halves of $\+e$.  If $s$ is $+1$, then the new
trivalent vertex will have the cyclic order $\+a < \+b < \+c < \+a$;
if $s_1$ is $-1$, then the edges $\+a$ and $\+b$ are swapped and the
new trivalent vertex gets the cyclic order $\+b < \+a < \+c < \+b$
instead.

\section{The homology complex of marked fatgraphs}
\label{sec:homology}

Betti numbers of a complex $(W_*, D_*)$ can be reckoned (via a
little linear algebra) from the matrix form $\+D^{(k)}$ of the boundary
operators $D_k$.  Indeed, given that $b_k := \dim H_k(W, D)$ and
$H_k(W, D) := Z_k(W,D) / B_k(W,D) = \Ker D_k / D_{k-1}(W_{k-1})$, by the
rank-nullity theorem we have $\dim \Ker D_k = \dim W_k - \rank \+D^{(k)}$
hence $b_k = \dim \Ker D_k - \dim D_{k-1}(W_{k-1}) = \dim W_k - \rank \+D^{(k)} -
\rank \+D^{(k-1)}$.

In order to compute the matrix $\+D^{(k)}$, we need to
compute the coordinate vector of $D_kx^{(k)}_j$ for all vectors
$x^{(k)}_j$ in a basis of $W_k$.
If $(W_*, D_*)$ is the fatgraph complex, then the basis vectors
$x^{(k)}_j$ are marked fatgraphs with $k$ edges, and the differential
$D_k$ is defined as an alternating sum of edge contractions.
Therefore, in order to compute the coordinate vector of
$D_kx^{(k)}_j$, one has to find the unique fatgraph $x^{(k-1)}_h$
which is isomorphic to a given contraction of $x^{(k)}_j$ and score a
$\pm 1$ coefficient depending on whether orientations agree or not.

Although this approach works perfectly, it is practically inefficient.
Indeed, lookups into the basis set $\{ x^{(k-1)}_{h=1,\ldots,N} \}$ of
$W_{k-1}$ require on average $O(N^2)$ isomorphism checks.
Still, we can take a shortcut:
if two topological fatgraphs $G$ and $G'$ are not isomorphic, so are
any two marked fatgraphs $(G,\nu)$ and $(G',\nu')$.
Indeed, rearrange the rows and columns of the boundary operator matrix
$\+D^{(k)}$ so that marked fatgraphs $(G,\nu)$ over the same
topological fatgraph $G$ correspond to a block of consecutive
indices. Then there is a rectangular portion of $\+D^{(k)}$ that is
uniquely determined by a pair of topological fatgraphs $G$ and $G'$.
The main function for computing the boundary operator matrix can thus loop
over pairs of topological fatgraphs, and delegate computing the each
rectangular block to specialized code.  There are $n! / |\Aut G|$
marked fatgraphs per given topological fatgraph $G$, so this
approach can cut running time down by $O((n!)^2)$.

The generation of inequivalent marked fatgraphs (over the same topological
fatgraph $G$) can be reduced to the (computationally easier)
combinatorial problem of finding cosets of a subgroup of the symmetric
group $\Perm{n}$.  In addition, the list of isomorphisms between $G$
and $G'$ can be cached and re-used for comparing all pairs of marked
fatgraphs $(G,\nu)$, $(G',\nu')$.  This strategy is implemented by two
linked algorithms:
\begin{enumerate}
\item\label{nfp:1} \fn{MarkedFatgraphPool}: Generate all inequivalent
  markings of a given topological fatgraph $G$.
\item\label{nfp:2} \fn{compute\_block}: Given topological fatgraphs $G$
  and $G'$, compute the rectangular block of a boundary operator
  matrix whose entries correspond to coordinates of $D(G,\nu)$
  w.r.t. $(G',\nu')$.
\end{enumerate}

\subsection{Generation of inequivalent marked fatgraphs}
\label{sec:markings}

For any marked fatgraph $(G, \nu)$, denote $[G, \nu]$ its isomorphism
class; recall that $B(-)$ is the functor associating a fatgraph with
the set of its boundary cycles.  Let $N(G)$ be the sets of all
markings over $G$ and $\tilde N(G)$ the set of isomorphism classes
thereof:
\begin{align*}
  N(G) &:= \{\, (G,\nu) \;|\; \nu\colon B(G) \to \{1,\ldots,n\} \,\},
  \\
  \tilde N(G) &:= \{\, [G,\nu] \;|\; \nu\colon B(G) \to \{1,\ldots,n\} \,\}.
\end{align*}
Let $(G, \bar\nu)$ be a chosen marked fatgraph.  Define a group
homomorphism:
\begin{equation}\label{eq:Phi}
\Phi: \Aut(G) \ni a \longmapsto \bar\nu \circ \Holes{a} \circ
\bar\nu^{-1} \in \Perm{n}.
\end{equation}
The set $P = \Phi(\Aut G)$ is a subgroup of $\Perm{n}$.  
\begin{lemma}\label{lemma:sigma}
  The marked fatgraphs $(G, \bar\nu)$ and $(G, \sigma\bar\nu)$ are
  isomorphic if and only if $\sigma \in P$.
\end{lemma}
\begin{proof}
  Let $\sigma \in P$, then $\sigma\inv \in P$ and there
  exists $a \in \Aut G$ such that:
  \begin{equation*}
    \sigma^{-1} = \bar\nu \circ \Holes{a} \circ \bar\nu^{-1},
  \end{equation*}
  whence:
  \begin{equation*}
    (\sigma \circ \bar\nu) \circ \Holes{a} \circ \bar\nu^{-1} = \idmap,
  \end{equation*}
  therefore $a$ induces a marked fatgraph isomorphism
  between $(G, \bar\nu)$ and $(G, \sigma \circ \bar\nu)$.

  Conversely, let $\hat\nu = \sigma\bar\nu$ and assume $(G, \bar\nu)$ and
  $(G, \hat\nu)$ are isomorphic as marked fatgraphs: then there exists
  $a \in \Aut G$ such that $\hat\nu \circ \Holes{a} \circ \bar\nu^{-1}$ is
  the identity. Given any $\bar\nu \circ \Holes{a'} \circ \bar\nu^{-1} \in
  P$ we have:
  \begin{multline*}
    P \ni \bar\nu \circ \Holes{a'} \circ \bar\nu^{-1} 
    = \bar\nu \circ (\hat\nu^{-1} \circ \hat\nu) 
      \circ \Holes{a} \circ \Holes{a}^{-1}
      \circ \Holes{a'} \circ \bar\nu^{-1} 
    \\
    = (\bar\nu \circ \hat\nu^{-1}) \circ \hat\nu \circ \Holes{a}
      \circ (\bar\nu^{-1} \circ \bar\nu) \circ \Holes{a^{-1}} 
      \circ \Holes{a'} \circ \bar\nu^{-1} 
    \\
    = \sigma\inv \circ (\hat\nu \circ
      \Holes{a} \circ \bar\nu^{-1}) \circ \bar\nu \circ \Holes{a^{-1} \circ
      a'} \circ \bar\nu^{-1}
    \\
    = \sigma\inv \circ \left( \bar\nu \circ
      \Holes{a^{-1} \circ a'} \circ \bar\nu^{-1} \right)
    \in \sigma\inv P,
  \end{multline*}
  therefore $P = \sigma\inv P$, so $\sigma \in P$.
\end{proof}

Define a transitive action of $\Perm{n}$ over $N(G)$ by $\sigma \cdot
(G,\nu) := (G, \sigma\nu)$; this descends to a transitive action of
$\Perm{n}$ on $\tilde N(G)$.  By the previous Lemma, $P$ is the
stabilizer of $[G, \bar\nu]$ in $\tilde N(G)$.
\begin{lemma}\label{lemma:markings}
  The action of $\Perm{n}$ on $\tilde N(G)$ induces a bijective
  correspondence between isomorphism classes of marked fatgraphs and
  cosets of $P$ in $\Perm{n}$.
\end{lemma}
\begin{proof}
  Given isomorphic marked fatgraphs $(G,\nu)$ and $(G,\nu')$, let
  $\sigma, \sigma' \in \Perm{n}$ be such that $\nu = \sigma \circ
  \bar\nu$ and $\nu' = \sigma' \circ \bar\nu$.  By definition of marked
  fatgraph isomorphism, there is $a \in \Aut G$ such that the
  following diagram commutes:
  \begin{equation*}
    \xymatrix{
      \Holes{G} & & \Holes{G} \\
      & \{1, \ldots, n\} &
      \ar^{\Holes{a}} "1,1";"1,3"
      \ar_{\sigma\circ\bar\nu = \nu} "1,1";"2,2"
      \ar^{\nu' = \sigma' \circ \bar\nu} "1,3";"2,2"
      }
  \end{equation*}
  Hence commutativity of another diagram follows:
  \begin{equation*}
    \xymatrix{
      \Holes{G} &&  \Holes{G} \\
      \{1,\ldots,n\} && \{1,\ldots,n\}
      \ar^{\Holes{a}} "1,1";"1,3"
      \ar_{\bar\nu} "1,1";"2,1"
      \ar^{\bar\nu} "1,3";"2,3"
      \ar_{\sigma\inv \sigma'} "2,3";"2,1"
    }
  \end{equation*}
  Thus $(G, \bar\nu)$ is isomorphic to $(G, \sigma\inv \sigma' \circ
  \bar\nu)$; therefore $\sigma\inv \sigma' \in P$, hence,
  \mbox{$\sigma' \in \sigma P$}, i.e., $\sigma$ and $\sigma'$ belong
  into the same coset of $P$.
  
  Conversely, let $\tau, \tau' \in \sigma P$; explicitly:
  \begin{equation*}
    \tau = \sigma \circ \bar\nu \circ \Holes{a} \circ \bar\nu\inv,
    \qquad
    \tau' = \sigma \circ \bar\nu \circ \Holes{a'} \circ \bar\nu\inv.
  \end{equation*}
  Set $\nu = \tau \circ \bar\nu$, $\nu' = \tau' \circ \bar\nu$;
  substituting back the definition of $\tau$, we have:
  \begin{equation*}
    \nu =  \sigma \circ \bar\nu \circ \Holes{a} \circ \bar\nu\inv \circ \bar\nu 
        =  \sigma \circ \bar\nu \circ \Holes{a},
  \end{equation*}
  whence $\bar\nu = \sigma\inv \circ \nu \circ \Holes{a}\inv$, and:
  \begin{equation*}
    \nu' = \sigma \circ \bar\nu \circ \Holes{a'}
    = \sigma \circ \left( \sigma\inv \circ \nu \circ B(a)\inv \right) \circ \Holes{a'}
    = \nu \circ \Holes{a\inv \circ a'},
  \end{equation*}
  therefore $a\inv \circ a'$ is an isomorphism between the marked
  fatgraphs $(G,\nu)$ and $(G,\nu')$.
\end{proof}

The following is an easy corollary of the transitivity of the
action of $\Perm{n}$ on $\tilde N(G)$.
\begin{lemma}\label{lemma:sigma+aut}
  Given any marking $\nu$ on the fatgraph $G$, there exist $\sigma
  \in \Perm{n}$ and $a \in \Aut G$ such that: $\nu = \sigma \circ
  \bar\nu \circ a$.
\end{lemma}
\begin{proof}
  By Lemma~\ref{lemma:markings}, there exists $\sigma \in \Perm{n}$
  such that $[G, \nu] = [G, \sigma \circ \bar\nu]$, i.e., $(G, \nu)$ is
  isomorphic to $(G, \sigma \circ \bar\nu)$.  If $a \in \Aut G$ is this
  fatgraph isomorphism, then the following diagram commutes:
  \begin{equation*}
    \xymatrix{
      \Holes{G} & & \Holes{G} \\
      & \{1, \ldots, n\} &
      \ar^{\Holes{a}} "1,1";"1,3"
      \ar_{\nu} "1,1";"2,2"
      \ar^{\sigma \circ \bar\nu} "1,3";"2,2"
    }
  \end{equation*}
  Therefore $\nu = \sigma \circ \bar\nu \circ \Holes{a}$.
\end{proof}

\paragraph{The \q{MarkedFatgraphPool} algorithm.}
Given a fatgraph $G$ and a "Fatgraph" object "$\+G$" representing it,
let us stipulate that $\bar\nu$ be the marking on $G$ that enumerates
boundary cycles of $G$ in the order they are returned by the function
"compute_boundary_cycles($\+G$)".  By Lemma~\ref{lemma:sigma+aut},
every $(G, \nu^{(j)})$ can then be expressed (up to isomorphism) as
$(G, \sigma^{(j)} \circ \bar\nu)$ with $\sigma^{(j)} \in \Perm{n}$.  The
set $\{ \sigma^{(j)} \}$ enumerates all distinct isomorphism classes
of marked fatgraphs over $G$ iff $\{ \sigma^{(j)}P \}$ runs over all
distinct cosets of $P$ in $\Perm{n}$ (by Lemma~\ref{lemma:markings}).

The "MarkedFatgraphPool" function computes the set $\tilde N(G)$ of
isomorphism classes $[G,\nu]$.
\begin{theorem}
\label{thm:MarkedFatgraphPool}
Given a "Fatgraph" $\+G$ as input, the output of
"MarkedFatgraphPool($\+G$)", as computed by
Algorithm~\ref{algo:MarkedFatgraphPool}, is a tuple 
"(graph, $P$, $A$, markings, orientable)", 
whose components are defined as follows:
  \begin{itemize}
  \item The "graph" item is the underlying "Fatgraph" object~"$\+G$".
  \item The "$P$" slot holds a list of all elements in the group $P =
    \Phi(\Aut G)$.
  \item A corresponding set of pre-image representatives (each element
    is an automorphism of "$\+G$") is stored into "$A$": permutation
    "$P[i]$" is induced by automorphism "$A[i]$", i.e., if $\pi=P[i]$
    and $a=A[i]$ then $\pi = \Phi(a)$.
  \item The "markings" item holds the list $\{ \sigma^{(j)} \}$ of
    distinct cosets of $P$ (representing inequivalent markings).
  \item "orientable" is a boolean value indicating whether any
    $(G,\nu)$ in the pool is orientable.\footnote{It is an immediate
      corollary of Lemma~\ref{lemma:sigma+aut} that if \emph{one} marked
      fatgraph $(G, \nu^*)$ has an orientation-reversing
      automorphisms, then \emph{every} marked fatgraph $(G, \nu)$ over
      the same topological fatgraph $G$ has an orientation-reversing
      automorphism.}
  \end{itemize}
\end{theorem}
We need a separate boolean variable to record the orientability of the
family of marked fatgraphs $N(G) = \{ (G, \nu) \}$, because the
automorphism group of a marked fatgraph $\Aut (G, \nu)$ can be a
proper subgroup of $\Aut G$: hence, $(G,\nu)$ can be orientable even
if $G$ is not.
  \begin{Algorithm}
    \caption{Compute the distinct markings of a given fatgraph.  Input
      to the algorithm is a \q{Fatgraph} object $\+G$; final result is
      a tuple \q{($\+G$, $P$, $A$, markings, orientable)} which
      represents the set $\tilde N(G)$ of isomorphism classes of
      marked fatgraphs.}
    \label{algo:MarkedFatgraphPool}
\begin{lstlisting}
def phi($a$, $\+G$):
  $\pi$ = |array of $n$ elements|                                  
  for src_index, src_cycle in enumerate($\+G$.boundary_cycles):
    dst_cycle = $a$.transform_boundary_cycle(src_cycle)
    if dst_cycle not in $\+G$.boundary_cycles:
      |abort and signal error to caller|
    else:
      dst_index = |index of |dst_cycle| in |$\+G$.boundary_cycles
    $\pi$[src_index] = dst_index                                   
  return $\pi$

def MarkedFatgraphPool($\+G$):
  P = |empty list|
  A = |empty list|
  # assume $(G,\nu)$ is orientable until we have counter-evidence
  orientable = True
  # step (1): loop over $\Aut G$
  for $a$ in $\+G$.automorphisms():                  |\n{nb1:start}|
    try:
      $\pi$ = phi($a$, $\+G$)
    except phi |failed|:
      |continue with next |$a$                           |\n{nb1:4}|
    if |permutation |$\pi$| is identity|:                |\n{nb1:5}|
      # found a new automorphism:
      #   - does it reverse orientation?
      if $a$.is_orientation_reversing():
        orientable = False                             |\n{nb1:end}|
      #   - does it define a new marking?
      if $\pi$ not in P:
        |append |$\pi$| to $P$|
        |append |$a$| to $A$|
  # step (2): enumerate cosets of $P$
  markings = [ ]
  for $\sigma$ in $\Perm{n}$:  
    for $\pi$ in $P$:                              |\n{nb2:unseen1}|
      if $\pi \circ \sigma$ in markings:
        |continue with next $\sigma$|              |\n{nb2:unseen2}|
      |add $\sigma$ to |markings
\end{lstlisting}
  \end{Algorithm}

\begin{proof}
  Generation of all inequivalent markings over $G$ is a direct
  application of Lemma~\ref{lemma:markings}, performed in two steps:
  \begin{enumerate}
  \item In the first step, for each automorphism $a \in \Aut G$
    compute the permutation $\Phi(a)$ it induces on the set of
    boundary components, and form the subgroup $P$.  The subgroup $P$
    and the associated set of automorphisms $A \subseteq \Aut G$ are
    stored in variables "P" and "A".
  \item In the second step, compute cosets of $P$ by exhaustive
    enumeration.  They are recast into the list $\{ \sigma^{(j)} \}$,
    which is stored into the "markings" variable.
  \end{enumerate}
  As an important by-product of the computation, the automorphism
  group $\Aut (G, \bar\nu)$ is computed, and used to determine if the
  marked fatgraphs in the pool are orientable.

  The auxiliary function \fn{phi} computes the permutation $\Phi(a) =
  \bar\nu \circ \Holes{a} \circ \bar\nu^{-1}$.  A permutation "$\pi$"
  is created and returned; it is represented by an array with $n$
  slots, which is initially empty and is then stepwise constructed by
  iterating over boundary cycles.  Indeed, the boundary cycle
  "src_cycle" is transformed according to $\Holes{a}$ and its position
  in the list of boundary cycles of $\+G$ is then looked up.  Note
  that this lookup may fail: there are in fact cases, in which the
  \fn{Fatgraph.isomorphisms} algorithm finds a valid mapping, that
  however does not preserve the markings on boundary cycles; such
  failures need to be dealt with by rejecting $a$ as a "Fatgraph"
  automorphism.

  Step {\sl (1)} of the computation is performed in lines
  \nr{nb1:start}--\nr{nb1:end}:
  \begin{itemize}
  \item Computation of the permutation $\pi$ (induced by $a$ on the
    boundary cycles of $\+G$) may fail; if this happens, the algorithm ignores
    $a$ and proceeds with another automorphism.
  \item If $a$ preserves the boundary cycles pointwise, then it
    induces an automorphism of the marked graph and we need to test
    whether it preserves or reverses orientation.
  \item There are $\card{\Ker\Phi}$ distinct automorphisms inducing
    the same permutation on boundary cycles: if "$\pi$" is already in
    $P$, discard it and continue with the next $a$.
  \end{itemize}

  By Lemma~\ref{lemma:markings}, there are as many distinct markings
  as there are cosets of "P" in $\Perm{n}$.  Step~{\sl (2)} of the
  algorithm proceeds by simply enumerating all permutations in
  $\Perm{n}$, with "marking" initially set to the empty list; for each
  permutation $\sigma$ a test is made as to whether $\sigma P$
  intersects the list "markings" (lines
  \nr{nb2:unseen1}--\nr{nb2:unseen2}); if it does not, then the
  marking induced by $\sigma$ is added to the list.
\begin{lstlisting}
\end{lstlisting}
\end{proof}

A constructive version of Lemma~\ref{lemma:sigma+aut} can now be
implemented: the following function \fn{index\_and\_aut}, given a
"Fatgraph" object $\+G$ and a marking, returns the permutation (by
index number "$j$" in "$\+G$.markings") and fatgraph automorphism $a =
\+G.A[i]$ such that the topological fatgraph "$\+G$" decorated with
"marking" is isomorphic (through $a$) to the same graph decorated with
"$\+G$.markings[$j$]".
\begin{lstlisting}
def index_and_aut($\+G$, marking):
  for $(i, \pi)$ in enumerate($\+G.P$):
    $\tau$ = $\sigma \circ \pi$
    if $\tau$ in $\+G$.markings:
      $j$ = |index of $\tau$ in |$\+G$.markings
      return $(j, \+G.A[i])$
    else:
      |continue with next $\pi$|
\end{lstlisting}
The algorithm enumerates all permutations $\pi \in P$, and compares
$\sigma \circ \pi$ to every element of "$\+G$.markings": by
Lemma~\ref{lemma:markings}, we know that one must match.

\subsection{Computing boundary operator matrix blocks}
\label{sec:blocks}

The differential $D(G,\nu)$ is computed by summing contractions of
regular edges in $G$ (with alternating signs); likewise, the matrix
block corresponding to coordinates of the families of marked fatgraphs
$\{(G,\nu)\}$ and $\{(G',\nu')\}$ can be decomposed into a sum of
blocks, each block representing the coordinates of $\{(G/e,\nu) \}_{e
  \in \Edges{G}}$ projected on the linear span of $\{(G',\nu')\}$.

More precisely, given any two fatgraphs $G_1$ (with $m$ edges) and
$G_2$ (with $m-1$ edges), let $X_1, X_2 \subseteq \RG_{g,n}$ be the
linear span of $N(G_1)$ and $N(G_2)$ respectively, and denote by
$\mathrm{pr}_{X_2}$ the linear projection on subspace $X_2$.  Recall
that, for any fatgraph $G$, we have $D(G) = \sum \pm d^{(e)}(G)$, where the
sum is taken over all regular edges $e$ of $G$, and $d^{(e)}$ is the
contraction of edge $e$.

Let $G$ be the fatgraph obtained by contracting the
chosen edge $e$ in $G_1$.
If $G_2$ and $G$ are isomorphic, then the three graphs are related
by the following diagram of fatgraph morphisms, where $f_1$ is the
contraction map and $f_2$ is a fatgraph isomorphism:
\begin{equation}
  \label{eq:fgdiag}
  \xymatrix{
    G_1 & 
    \\
    G   &  G_2
    \ar_{f_1} "1,1";"2,1"
    \ar_{f_2}^{\sim} "2,1";"2,2"
  }
\end{equation}

The above diagram~\eqref{eq:fgdiag} functorially induces a diagram on
the set of boundary cycles:
\begin{equation}
  \label{eq:nbdiag}
  \xymatrix{
    \Holes{G_1} & \{1, \ldots, n\}
    \\
    \Holes{G}   &  \Holes{G_2}
    \ar_{\Holes{f_1}}^{=} "1,1";"2,1"
    \ar_{\Holes{f_2}}^{\sim} "2,1";"2,2"
    \ar^{\nu_1} "1,1";"1,2"
    \ar_{\nu} "2,1";"1,2"
    \ar_{\nu_2} "2,2";"1,2"
  }
\end{equation}
Diagram~\eqref{eq:nbdiag} commutes iff $f_1$, $f_2$ can be extended to
morphisms of marked fatgraphs $\hat f_1: (G_1, \nu_1) \to (G, \nu)$ and
$\hat f_2 : (G, \nu) \to (G_2, \nu_2)$.

Now choose "Fatgraph" objects $\+G_1$, $\+G$, $\+G_2$ representing
$G_1$, $G$, $G_2$.

Let $\bar\nu_1$, $\bar\nu$, $\bar\nu_2$ be the markings on $G_1$, $G$,
$G_2$ that enumerate boundary cycles in the order they are returned by
the function "compute_boundary_cycle" applied to "$\+G_1$", "$\+G$", "$\+G_2$"
respectively.  Define $\phi_1, \phi_2 \in \Perm{n}$ by:
\begin{equation}
  \label{eq:phi-maps}
  \phi_1 := \bar\nu \circ \Holes{f_1} \circ \bar\nu_1^{-1},
  \qquad
  \phi_2 := \bar\nu_2 \circ \Holes{f_2} \circ \bar\nu^{-1}.
\end{equation}

\begin{lemma}\label{lemma:unique-nu}
  Given any marking $\nu_1$ on $G_1$, choose $\sigma_1 \in \Perm{n}$ such that
  \mbox{$\nu_1 = \sigma_1 \circ \bar\nu_1$} and define:
  \begin{equation}\label{eq:nu2}
    \nu_2 := \sigma_1 \circ \phi_1\inv \circ \phi_2\inv \circ \bar\nu_2 ,
  \end{equation}
  Then $\nu_2$ is the unique marking on $G_2$ such that
  diagram~\eqref{eq:nbdiag} commutes.
\end{lemma}
\begin{proof}
  Let $\sigma_2 := \sigma_1 \circ \phi_1\inv \circ \phi_2\inv$.  We
  need to prove that the external square in diagram~\eqref{eq:nbdiag}
  is commutative; indeed, we have:
  \begin{multline*}
    \sigma_2 = \sigma_1 \circ (\bar\nu_1 \circ \Holes{f_1}\inv \circ
    \bar\nu\inv) \circ ( \bar\nu \circ \Holes{f_2}\inv \circ
    \bar\nu_2)
    \\
    = \sigma_1 \bar\nu_1 \circ \Holes{f_2 \circ f_1}\inv \circ
    \bar\nu_2\inv,
  \end{multline*}
  so that:
  \begin{multline*}
    \nu_2 \circ \Holes{f_2} \circ \Holes{f_1} = \sigma_2 \circ \nu_2
    \circ \Holes{f_2 \circ f_1}
    \\
    = \sigma_1 \bar\nu_1 \circ \Holes{f_2 \circ f_1}\inv \circ
    \bar\nu_2\inv \circ \nu_2 \circ \Holes{f_2 \circ f_1}
    \\
    = \sigma_1 \circ \bar\nu_1 = \nu_1.
  \end{multline*}
  The uniqueness assertion is of immediate proof, since maps
  $\Holes{f_1}$ and $\Holes{f_2}$ are invertible.
\end{proof}

Let $\+p_1$, $\+p_2$ be the "MarkedFatgraphPool" output corresponding
to $\+G_1$, $\+G_2$, and let $\{\nu_1^{(j)}\}_{j=1,\ldots,N_1}$,
$\{\nu_2^{(k)}\}_{k=1,\ldots,N_2}$ be the enumeration of fatgraph
markings corresponding to items in the lists "$\+p_1$.markings" and
"$\+p_2$.markings" respectively.
\begin{lemma}
  \label{lemma:unique-k-and-s}
  For any regular edge $e$ of $G_1$, and any choice of $j \in \{1,
  \ldots, N_1\}$, there exist unique $k \in \{1, \ldots, N_2\}$ and $s
  \in \{-1, 0, +1\}$ such that:
  \begin{equation}
    \label{eq:ComputeBlock}
    \mathrm{pr}_{X_2} \left( d^{(e)} [G_1, \nu_1^{(j)}] \right) 
    = s\cdot [G_2, \nu_2^{(k)}].
  \end{equation}
\end{lemma}
\begin{proof}
  If $G_2$ and $G = G_1/e$ are not isomorphic, then, for any marking
  $\nu_1$, $d(G_1, \nu_1)$ has no component in the subspace $X_2 = \{
  (G_2, \nu_2) \}$, so the assertion is true with $s = 0$.

  Otherwise, by Lemma~\ref{lemma:unique-nu}, given $\nu_1 =
  \nu_1^{(j)}$ there is a unique $\nu_2$ such that $s$ can be
  non-null; by Lemma~\ref{lemma:sigma+aut}, there exist $\nu_2^{(k)}
  := \sigma_2^{(k)} \circ \bar\nu_2$ and $a \in \Aut G$ such that:
  \begin{enumerate}
  \item the marked fatgraph $(G_2, \nu_2^{(k)})$ is a representative
    of the isomorphism class $[G_2,\nu_2]$;
  \item $a$ gives the isomorphism between marked
    fatgraphs $(G_2,\nu_2)$ and $(G_2, \nu_2^{(k)})$;
  \item $\nu_2^{(k)}$ is the marking on $G_2$ represented by $k$-th
    item in list "$\+p_2$.markings".
  \end{enumerate}
  The coefficient $s$ must then be $\pm 1$ since both $(G_2,
  \nu_2^{(k)})$ and $d^{(e)}(G_1, \nu_1^{(j)})$ are (isomorphic to)
  elements in the basis of $X_2$.
\end{proof}

\begin{theorem}
  \label{thm:triplets}
  Given "MarkedFatgraphPool" objects $\+p_1$, $\+p_2$, and a chosen
  edge $\+e$ of $\+G_1$, the function "compute_block" in
  Algorithm~\ref{algo:ComputeBlock} returns the set $S$ of all
  triplets $(j, k, s)$ with $s = \pm 1$ such that:
  \begin{equation}\label{eq:triplets}
    \mathrm{pr}_{X_2} \left(d^{(e)}[G_1, \nu_1^{(j)}] \right)
    = s \cdot [G_2, \nu_2^{(k)}].
  \end{equation}
\end{theorem}
\begin{Algorithm}
  \caption{%
    Return the set $S$ of triplets $(j,k,s)$ such that
    eq.~\eqref{eq:triplets} holds for $(G_1, \nu_1^{(j)})$ and $(G_2,
    \nu_2^{(k)})$ obtained by contracting $\+e$ in all marked graphs
    in $\+p_1$ and projecting onto graphs in the $\+p_2$ family.}
  \label{algo:ComputeBlock}
\begin{lstlisting}
def compute_block($\+p_1$, $\+e$, $\+p_2$):
  $\+G_1$ = $\+p_1$.graph
  $\+G_2$ = $\+p_2$.graph
  $\+G$ = contract($\+G_1$, $\+e$)
  if $\+G$ and $\+G_2$ |are not isomorphic|:  |\n{bl:noiso1}|
    return |empty list|                       |\n{bl:noiso2}|
  else:
    result = |empty list|
    $\+f_2$ = |first isomorphism computed by |Fatgraph.isomorphisms($\+G$, $\+G_2$) |\n{bl:2}|
    $\phi_1\inv$ = compute_phi1_inv($\+G$, $\+G_1$, $\+e$)
    $\phi_2\inv$ = compute_phi2_inv($\+G$, $\+G_2$, $\+f_2$)
    for $(j, \sigma)$ in enumerate($\+p_1$.markings):  |\n{bl:final1}|
      $(k, \+a)$ = index_and_aut($\+p_2$, $\sigma \circ \phi_1\inv \circ \phi_2\inv$)
      $p$ = $\+G_1$.orient[$\+e$]
      $s$ = $(-1)^p$ * compare_orientations($\+f_2$) * compare_orientations($\+a$) |\n{bl:s}|
      |append $(j, k, s)$ to |result                   |\n{bl:final2}|
    return result
  
def compute_phi1_inv($\+G$, $\+G_1$, $\+e$):
  $\tau$ = |empty array of $n$ elements|
  for $i$, $\+b$ in enumerate($\+G_1$.boundary_cycles):
    $\+b'$ = contract_boundary_cycle($\+G_1$, $\+b$, $\+e$)
    $i'$ = |index of $\+b'$ in |$\+G$.boundary_cycles
    $\tau[i']$ = $i$
  return $\tau$

def compute_phi2_inv($\+G$, $\+G_2$, $f_2$):
  $\tau'$ = |empty array of $n$ elements|
  for $i$, $\+b$ in enumerate($\+G_2$.boundary_cycles):
    $\+b'$ = transform_boundary_cycle($f_2$, $\+b$)
    $i'$ = |index of $\+b'$ in |$\+G$.boundary_cycles
    $\tau'[i']$ = $i$
  return $\tau'$
\end{lstlisting}
\end{Algorithm}
\begin{proof}
  The algorithm closely follows the computation done before
  Lemmas~\ref{lemma:unique-nu} and in the proof of
  Lemma~\ref{lemma:unique-k-and-s}.

  If $G_2$ and $G = G_1/e$ are not isomorphic, then $d^{(e)}[G_1,
  \nu_1]$ has no component in the subspace $X_2$ generated by $\{
  [G_2, \nu_2] \}$, whatever the marking $\nu_1$.  The assertion is
  thus satisfied by $S = \varnothing$, i.e., an empty list of triplets
  $(j, k, s)$ (lines~\nr{bl:noiso1}--\nr{bl:noiso2} in
  Algorithm~\ref{algo:ComputeBlock}).

  If $G_2$ is isomorphic to $G = G_1/e$ through $f_2$, then
  Lemma~\ref{lemma:unique-nu} provides the explicit formula
  $\nu_2^{(k)} = \sigma_1^{(j)} \circ \phi_1\inv \circ \phi_2\inv
  \circ \bar\nu_2$, where $\sigma_1^{(j)} = \nu_1^{(j)} \circ
  \bar\nu_1\inv$.

  By assumption, $\bar\nu_1$ numbers the boundary cycles on $G_1$ in
  the order they are returned by running function
  "compute_boundary_cycles" on $\+G_1$, so $\sigma_1^{(j)}$ is the
  permutation corresponding to the $j$-th element in "$\+p_1$.markings".

  The map $\phi_1$ is easy to compute: again, given that both
  $\bar\nu$ and $\bar\nu_1$ number the boundary cycles of $G$ and
  $G_1$ in the order they are returned by "compute_boundary_cycles",
  the auxiliary function "compute_phi1_inv" incrementally builds the
  result by looping over "$\+G_1$.boundary_cycles", contracting the
  target edge, and mapping the corresponding indices.
 
  Computation of the map $\phi_2$ depends on the isomorphism $f_2$;
  however, two different choices for $f_2$ will not change the outcome
  of the algorithm: in the final loop at lines~\nr{bl:final1}--\nr{bl:final2},
  only the sign of $f_2$ is used, and the sign is constant across all
  isomorphisms having the same source and target fatgraphs (iff they
  are both orientable).
  Computation of $\phi_2\inv$ (in the auxiliary function
  "compute_phi2_inv") is done in the same way as the computation of
  $\phi_1\inv$, except we transform $b$ to $b'$ by means of
  "transform_boundary_cycle($f_2$, $-$)", i.e., $\Holes{f_2}$.

  Finally, for every marking "$\sigma_1^{(j)}$" in "$\+p_1$.markings"
  (representing $\nu_1^{(j)}$), we know by
  Lemma~\ref{lemma:unique-k-and-s} that there is a unique index $k$
  and $a \in \Aut G_2$ such that: $\sigma_1^{(j)} \circ \phi_1\inv
  \circ \phi_2\inv = \sigma_2^{(k)}$ is the $k$-th item in
  "$\+p_2$.markings" (representing $\nu_2^{(k)}$), and such that the
  following chain:
  \begin{equation*}
    \xymatrix{
      G_1 & G & G_2 & G_2
      \ar_{f_1}"1,1";"1,2"
      \ar_{f_2}^{\simeq}"1,2";"1,3"
      \ar_{a}^{\simeq}"1,3";"1,4"
    }
  \end{equation*}
  extends to a marked fatgraph morphism:
  \begin{equation*}
    \xymatrix{
      (G_1, \nu_1^{(j)}) & (G, \nu) & (G_2, \nu_2)  & (G_2, \nu_2^{(k)}).
      \ar_{\hat f_1}"1,1";"1,2"
      \ar_{\hat f_2}^{\simeq}"1,2";"1,3"
      \ar_{\hat a}^{\simeq}"1,3";"1,4"
    }
  \end{equation*}

  The sign $s$ is then obtained by comparing the orientation
  $\omega_2$ of $(G_2, \nu_2^{(k)})$ with the push-forward orientation
  ${(a \circ f_2 \circ f_1)}_* \omega_1$, where $\omega_1$ is the
  orientation on $(G_1, \nu_1^{(j)})$, and multiplying by the
  alternating sign from the homology differential.  There are four
  components that make up $s$:
  \begin{itemize}
  \item the sign given by the contraction $f_1$: this is $+1$
    by definition since the ``child'' fatgraph $G$ inherits the
    orientation from the ``parent'' fatgraph $G_1$;
  \item the sign given by the isomorphism $f_2$: this is obtained by
    comparing $(f_2)_* \omega$ with $\omega_2$, which is implemented
    for a generic isomorphism by the function "compare_orientations";
  \item the sign of the automorphism $a$ of "$G_2$" which
    transforms the push-forward marking into the chosen representative
    in the same orbit: this again can be computed by comparing $(a)_*
    \omega_2$ with $\omega_2$ and only depends on the action of $a$ on
    edges of $G_2$;
  \item the alternating sign from the homology differential, which
    only depends on the position $p$ of edge $e$ within the order
    $\omega_1$.
  \end{itemize}
  The product of the three non-trivial components is returned as the
  sign $s$ (line~\nr{bl:s}).
\end{proof}

\subsection{Matrix form of the differential $D$}
\label{sec:MgnChainComplex}

The "compute_boundary_operators" function
(Algorithm~\ref{algo:total-matrix}) computes the matrix form $\+D^{(m)}$
of the differential $D$ restricted to the linear space generated by
fatgraphs with $m$ edges.  

Input to the function are the number $m$ and the list of graphs,
divided by number of fatgraph edges: "graphs[$m$]" is the list of
fatgraphs with $m$ edges.

The output matrix $\+D^{(m)}$ is constructed incrementally: it starts
with all entries set to $0$, and is then populated blockwise. Indeed,
for every pair of "MarkedFatgraphPool" objects $\+p_1$ (from a graph
with $m$ edges) and $\+p_2$ (with $m-1$ edges), and every non-loop
edge $\+e$, the rectangular matrix block whose upper-left corner is at
indices $j_0, k_0$ is summed the block resulting from
"compute_block($\+p1$, $\+e$, $\+p_2$)". 
\begin{Algorithm}
  \caption{Compute the boundary operator matrix, block by block.}
  \label{algo:total-matrix}
\begin{lstlisting}
def compute_boundary_operator($m$, graphs):
  $N_1$ = |number of graphs with $m$ edges|
  $N_2$ = |number of graphs with $m-1$ edges|
  $D^{(m)}$ = |$N_1 \times N_2$ matrix, initially null|
  $j_0$ = 0
  for $\+G_1$ in graphs[m]:                       
    $\+p_1$ = MarkedFatgraphPool($\+G_1$)
    $k_0$ = 0
    for $\+G_2$ in graphs[m-1]:
      $\+p_2$ = MarkedFatgraphPool($\+G_2$)
      for $\+e$ in $\+G_1$.edges:
        if $\+e$ is a loop:
          continue |with next $\+e$|
        for $(j, k, s)$ in compute_block($\+p_1$, $\+e$, $\+p_2$):
          |add $s$ to entry $D^{(m)}[k + k_0, j + j_0]$|
      |increment $k_0$ by the number of inequivalent markings in $\+p_2$|
    |increment $j_0$ by the number of inequivalent markings in $\+p_1$|
  return $D^{(m)}$
\end{lstlisting}
\end{Algorithm}

\section{Conclusions}
\label{sec:conclusions}

A Python implementation\footnote{%
  Code publicly available at \url{http://code.google.com/p/fatghol}.
} of the algorithms presented in this paper has been actually used to
compute the Betti numbers of all $\M_{g,n}$ with $2g+n \leq 6$.  The
results are summarized in Table~\ref{tab:betti}.  Corresponding
calculations based on theoretical results are scattered across a wide
array of publications.  For $g \geq 1$, the groups $H^1(\M_{g,n},
\setQ)$ are known from the works of Mumford \cite{mumford:1967} and
Harer \cite{harer;second-homology}; $H^2(\M_{g,n}, \setQ)$ has been
computed also by Harer in \cite{harer;second-homology}; a
comprehensive statement with a new proof is given by Arbarello and
Cornalba in \cite{arbarello-cornalba:2009} (where a minor mistake in
Harer's statement is corrected).  The complete homology of $\M_{1,2}$
and $\M_{2,1}$ has been published in Godin's paper
\cite{godin:homology}.  The homology of the $\M_{0,*}$ spaces is
computed in \cite[Corollary, 3.10]{getzler:1995}; see
\cite{algori:mathoverflow-q-38968} for an alternative approach using
results from \cite{totaro:1996} to compute the Poincaré polynomial of
$\M_{0,n}$.  The Poincaré-Serre polynomial of $M_{2,2}$ follows as a
special case of Corollary~III.2.2 in Tommasi's \cite{tommasi:phd}; the
results also follows by combining \cite[p.~22]{getzler:1998} with
\cite[Appendix~A]{hulek+tommasi}.  The rational cohomology of
$\M_{1,4}$ is completely described in Theorem~1 of
\cite{gorinov:rational-cohomology-M1n}; the Betti numbers were already
present in \cite[p.~22]{getzler:1999}.  In all these cases, the
numerical results agree with the values in Table~\ref{tab:betti}.

An internal verification step in the code computes the classical and
virtual Euler characteristics of the fatgraph complex; the computed
values match those published in
\cite{harer-zagier;euler-characteristic,bini-harer,bini-gaiffi-polito},
where they are derived by theoretical means.

As a side effect of the computation, the entire family of fatgraphs
$\R_{g,n}$ (with $2g + n \leq 6$) has been computed, and for each
fatgraph the isomorphism group is known.  The full list of fatgraphs
and their isomorphisms is too long to print here, but the data is
publicly available at \url{http://fatghol.googlecode.com/download/list}.
Tables~\ref{tab:abs-fatgraphs} and~\ref{tab:marked-fatgraphs} provide
a numerical summary of the results.
\begin{table}
  \centering
  \begin{tabular}{lccccccccccc|c}
    \toprule
    {\small No.\ of edges:}
                 &\it 12&\it 11&\it 10& \it 9 & \it 8 & \it 7 & \it 6 & \it 5 & \it 4 & \it 3 & \it 2   & {Total}\\
    \midrule
    $g=0$, $n=3$ &      &      &      &       &       &       &       &       &       &      2&      1  &       3\\
    $g=0$, $n=4$ &      &      &      &       &       &       &      6&      6&      7&      6&         &      25\\
    $g=0$, $n=5$ &      &      &      &     26&     26&     72&    103&     65&     21&       &         &     313\\
    $g=0$, $n=6$ &   191&   191&   866&   1813&   1959&   1227&    418&     76&       &       &         &    6741\\
    $g=1$, $n=1$ &      &      &      &       &       &       &       &       &       &      1&      1  &       2\\
    $g=1$, $n=2$ &      &      &      &       &       &       &      5&      5&      8&      8&         &      26\\
    $g=1$, $n=3$ &      &      &      &     46&     46&    162&    256&    198&     72&       &         &     780\\
    $g=1$, $n=4$ &   669&   669&  3442&   7850&   9568&   6752&   2696&    562&       &       &         &   32208\\
    $g=2$, $n=1$ &      &      &      &      9&      9&     29&     52&     45&     21&       &         &     165\\
    $g=2$, $n=2$ &   368&   368&  2005&   4931&   6543&   5094&   2279&    546&       &       &         &   22134\\
    \bottomrule
  \end{tabular}
  \caption{\label{tab:abs-fatgraphs}%
    Number of distinct abstract fatgraphs with the given genus $g$ and
    number of boundary cycles $n$.  For readability, null
    values have been omitted and the corresponding entry left blank.}
\end{table}

\begin{table}
  \centering
  \begin{tabular}{@{\it}r|cccccccccc}
    \toprule
    {\normalfont No.\ of edges}
      & $\M_{0,3}$ & $\M_{0,4}$ & $\M_{0,5}$ & $\M_{0,6}$ & $\M_{1,1}$ & $\M_{1,2}$ & $\M_{1,3}$ & $\M_{1,4}$ & $\M_{2,1}$ & $\M_{2,2}$ \\
    \midrule
    12&           &           &           &   122880  &           &           &           &    14944  &           &      713  \\
    11&           &           &           &   616320  &           &           &           &    81504  &           &     3983  \\
    10&           &           &           &  1274688  &           &           &           &   185760  &           &     9681  \\
    9 &           &           &   2240    &  1359840  &           &           &    236    &   227564  &      9    &    12927  \\
    8 &           &           &   8160    &   862290  &           &           &    918    &   160128  &     28    &    10077  \\
    7 &           &           &  11280    &   294480  &           &           &   1440    &    63756  &     43    &     4519  \\
    6 &           &     64    &   7260    &    49800  &           &      9    &   1112    &    13000  &     39    &     1057  \\
    5 &           &    144    &   2112    &     3024  &           &     15    &    408    &     1008  &     20    &       97  \\
    4 &           &     99    &    210    &           &           &     10    &     54    &           &      3    &           \\
    3 &          4&     20    &           &           &      1    &      3    &           &           &           &           \\
    2 &          3&           &           &           &      1    &           &           &           &           &           \\
    \midrule
    {\em Total}& 7&    327    &  31262    &  4583322  &      2    &     37    &   4168    &   747664  &    142    &    43054  \\
    \bottomrule
  \end{tabular}
  \caption{\label{tab:marked-fatgraphs}%
    Number of distinct orientable marked fatgraphs in the Penner-Kontsevich
    complex of each of the indicated $\M_{g,n}$ spaces.  For readability, null
    values have been omitted and the corresponding entry left blank.}
\end{table}

\subsection{Performance}
\label{sec:performance}

\begin{figure}
  \centering
  \includegraphics[width=\linewidth]{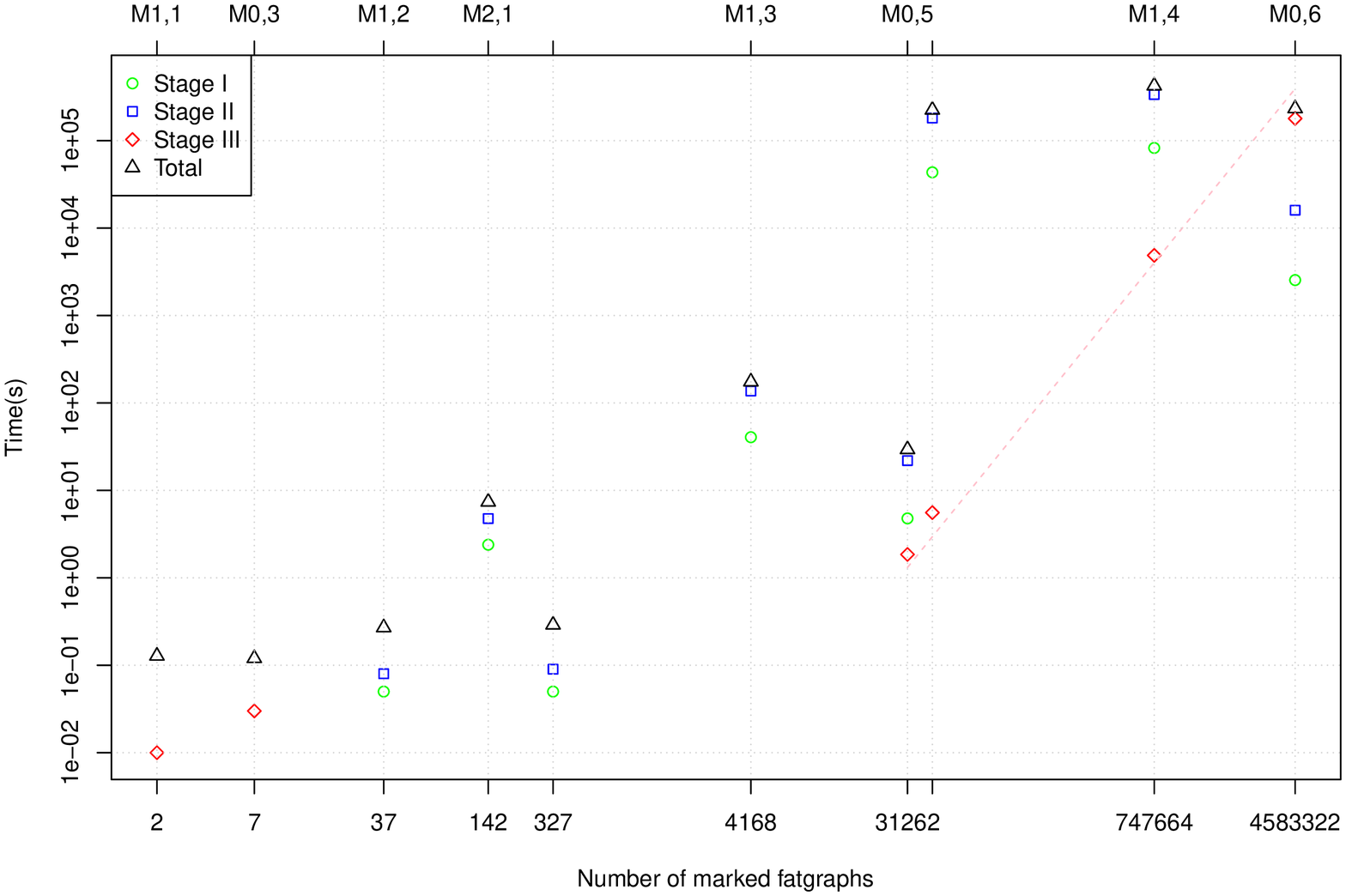}
  \caption{%
    Scatter plot of the data in Table~\ref{tab:times}.  Both axes use
    log-scale.  Note how Stage~III (computation of the boundary
    operators rank) becomes the dominant task as the number of marked
    fatgraphs increases.}
  \label{fig:times}
\end{figure}

Table~\ref{tab:times} gives a summary of the running times obtained on
the \texttt{idhydra.uzh.ch} cluster at the University of Zurich;
Figure~\ref{fig:times} provides a graphical representation of the same
data. The computational demands of the code are such that the homology
of $\M_{g,n}$ can actually be computed on desktop-class hardware for
$2g + n < 6$.
\begin{table}
  \centering
  \begin{tabular}{p{3em}rrrr}
    \toprule
    \multicolumn{2}{l}{Time (s):
      \hfill
    Stage~I  }&
    Stage~II  &
    Stage~III &
    Total     \\
    \midrule
    $\M_{0,3}$ &    $<1$ms &    $<1$ms &      0.03 &      0.12 \\
    $\M_{0,4}$ &      0.05 &      0.09 &    $<1$ms &      0.29 \\
    $\M_{0,5}$ &      4.78 &     21.91 &      1.85 &     29.43 \\
    $\M_{0,6}$ &   2542.56 &  16011.70 & 179157.39 & 233007.06 \\
    $\M_{1,1}$ &    $<1$ms &    $<1$ms &     0.010 &     0.128 \\
    $\M_{1,2}$ &      0.05 &      0.08 &    $<1$ms &      0.27 \\
    $\M_{1,3}$ &     40.56 &    136.88 &    $<1$ms &    174.75 \\
    $\M_{1,4}$ &  82486.51 & 336633.75 &   4872.69 & 424615.85 \\
    $\M_{2,1}$ &      2.39 &      4.76 &    $<1$ms &      7.39 \\
    $\M_{2,2}$ &  43402.18 & 181091.11 &      5.57 & 224694.61 \\
    \bottomrule
  \end{tabular}
  \caption{\label{tab:times}%
    Total CPU time (seconds) used by the Betti numbers computation for the
    indicated $\M_{g,n}$ spaces.  The \Cpp{} library \textsc{LinBox} 
    \cite{linbox:website,linbox} was used for the rank computations in Stage~III. 
    Running time was sampled
    on the \texttt{idhydra.uzh.ch} computer of the University of Zurich, 
    equipped with 480GB of RAM and Intel Xeon CPUs model X7542 running at 2.67GHz; 
    Python version 2.6.0 installed on the SUSE Linux Enterprise Server 11 64-bits 
    operating system was used to execute the program.  The system timer has a 
    resolution of $1$ms, but times are less accurate than that, because of the 
    scheduling jitter in multitasking systems. 
    The ``Total'' column does not just report the sum of the three stages, 
    but also accounts for the time the program spent in I/O and memory management.
  }
\end{table}

The scatter plot in Figure~\ref{fig:times} shows that the time
spent in computation of the $\+D^{(m)}$ matrix ranks done in Stage~III
can become the dominant contribution to the total running time as the
number of fatgraphs increases.  This highlights a limitation of the
program: the large number of fatgraphs in the Kontsevich complex might
turn out to be a challenge for today's sparse linear algebra software.

However, the set of fatgraphs for a given $(g,n)$ pair has to be
generated prior to computing the matrices $\+D^{(m)}$: a very large
set of graphs can exhaust the computer's memory long before
computation time becomes a blocking issue.

\subsection{Application to other fatgraph complexes}
\label{sec:other}

In \cite{godin:homology}, V.~Godin defined a ``bordered fatgraph
complex'', which computes the integral homology of the moduli
spaces of Riemann surfaces with boundaries.  Godin's fatgraphs extend
the abstract fatgraph by requiring that a \emph{leaf} (i.e., a
univalent vertex), and only one, is present in each boundary cycle.
The bordered fatgraph complex is then constructed exactly as the
fatgraph complex presented here, with the proviso that an edge ending
in a univalent vertex is never contracted: hence, the differential $D$
is given by the sum of contraction of non-loop non-leaf edges.

The algorithms of this paper can easily be adapted to compute the
homology of Godin's bordered fatgraph complex: after generating the
family of marked fatgraphs of a given $(g,n)$ pair, we decorate each
marked fatgraph with leaves; compute the matrix form of the
differential $D$ and then reduce it to Smith normal form to reckon the
\emph{integral} homology modules of the moduli space of bordered
surfaces.

There is no need for checking duplicates in the set of bordered
fatgraphs so generated,\footnote{%
  If two ``bordered fatgraphs'' were isomorphic, they would remain
  such if we remove the leaves and the edge supporting them, which
  would give us isomorphic marked fatgraphs}\space%
therefore the decoration step can be implemented efficiently.  A
shortcut can also be taken in computing the matrix $D$: since leaf
edges are never contracted, the differential on bordered fatgraphs can
be deduced easily from the differential on marked fatgraphs.  However,
the number of bordered fatgraphs is much larger than the number of
marked fatgraphs;\footnote{%
  A leaf may be regarded as a choice of an edge or a vertex along a
  boundary cycle; if there are $p_i$ vertices (counted with
  multiplicities) and $q_i$ edges along the $i$-th boundary cycle,
  then the number of ways we could possibly add leaves to a marked
  fatgraph $G$ is $r_1r_2 \cdots r_n$, with $r_i = p_i + q_i$ so that:
  \begin{equation*}
    r_1 + r_2 + \cdots + r_n 
    = \sum_i p_i + \sum_i q_i 
    = \sum_{v\in V(G)} z_v + 2m
    = 4m,
  \end{equation*}
  where $m$ is the total number of edges and $z_v$ is the valence of
  vertex $v$.
}%
this means that the final linear algebra computations require even more
computational resources than they do for $\M_{g,n}$ computations.

\subsection{Future development directions}
\label{sec:future}

There are a number of directions in which the current algorithms and
code could be improved.

As already noted, the generation algorithms produce quite a number of
duplicates, that have to be removed using a quadratic-complexity
procedure.  A variant of the ``isomorph-free generation'' algorithm of
McKay \cite{mckay:isomorph-free} could replace the naive
"MgnTrivalentGraphs" code; the question of which algorithm would be
faster has probably to be sorted out empirically, the critical
performance factor being the number of times the ``isomorphism'' test
is invoked.  

Another approach would be to turn the generation procedure ``upside
down'': instead of starting with trivalent graphs and contracting
edges, one could start with $(g,n)$-fatgraphs with one vertex and
expand those until the whole set of fatgraphs is generated.  This
would have the advantage that the chromatic fatgraph polynomial of
Bollobás and Riordan
\cite{bollobas-riordan:chromatic-polynomial-fatgraphs} is available as
an invariant to speed up the isomorphism procedure.  On the other
hand, the number of fatgraphs generated this way seems consistently
larger than the number of fatgraphs generated with the procedure
adopted here (see Section~\ref{sec:N2}).

So far, the major obstacle to applying the algorithms of this paper to
a wider range of moduli spaces has been the large number of fatgraphs
involved: it affects both the total run time and memory consumption of
the code.  Most algorithms described here lend themselves naturally to
parallelization, so it would be possible to rewrite the program to
exploit several processors and distributed memory, which could solve
both issues. However, the number of generated fatgraphs grows
super-exponentially in the asymptotic limit
\cite{bender-canfield,bender-gao-richmond}, so any implementation of
the algorithms outlined here will soon hit the limit of any
present-day computing device.  The question remains open, whether more
significant result could be obtained before hitting the limits of
today's computers.

\section*{Acknowledgments}

The computational power necessary for the production runs of the
software was kindly donated by the Informatik Dienste of the
University of Zurich and by the Organisch-Chemisches Institut (and
particularly the ``Theoretical Chemistry and Computational Grid
Applications'' research group led by Kim Baldridge), whose support is
here gratefully acknowledged.

I am specially indebted, and would like to express my gratitude, to
Enrico Arbarello, Domenico Fiorenza and Gilberto Bini, for their
constant encouragement and manifold support during the long
preparation of this paper.

I also would like to thank Orsola Tommasi and Dan Petersen, who
pointed out inaccuracies in an earlier versions of this text and
suggested improvements.

\appendix

\section{Comparison of fatgraph generation methods}
\label{sec:N}

This section compares three different approaches to generating
trivalent fatgraphs: namely, we compare the "MgnTrivalentGraphs"
algorithm described in Section~\ref{sec:MgnTrivalentGraphs} with two
alternatives.\footnote{%
  The author is aware of no other algorithm for generating the set of
  all fatgraphs.  The comparison here is taken with the solutions used
  in earlier attempts of implementation of the
  \href{http://code.google.com/p/fatghol}{FatGHoL} software.  }
Table~\ref{tab:N} presents a summary of results.

None of the suggested algorithms is capable of directly producing an
isomorph-free set of distinct fatgraphs; they all produce a larger set
of fatgraphs that must be reduced by taking only one representative
per isomorphism class of fatgraphs.  Therefore, Table~\ref{tab:N} also
reports the actual number of distinct fatgraphs for a given $g,n$ pair;
not all counts are known: a cell is left empty when the corresponding
count has not yet been computed. From the results gathered so far, it
is apparent that \emph{all} algorithms overestimate the actual
number of fatgraphs.
\begin{table}[tbh]
  \def\C{\multicolumn{1}{c}}
  \centering
  \begin{tabular}{cc|rrrrr}
    \toprule
    $g$&$n$&\C{$N$}&      {$N_1^+$}&     {$N_2^+$}&               \C{$N_3$} \\
    \midrule
    0 &3 &    2 &  \,   \,   \,$-$ &  \,   \, 15 &$                5\,760 $\\
    0 &4 &    6 &  \,   \,   \, 84 &  \,   \,630 &$ 1.072964\times 10^{13} $\\
    0 &5 &   26 &  \,   \,   \,936 &  \, 15\,015 &$ 4.593811\times 10^{24} $\\
    0 &6 &  191 &  \,   \,  8\,892 &  \,306\,306 &$ 6.326929\times 10^{37} $\\
    0 &7 &      &  \,   \,114\,600 & 5\,819\,814 &$ 1.132261\times 10^{52} $\\
    1 &1 &    1 &  \,   \,   \,$-$ &  \,   \, 15 &$                5\,760 $\\
    1 &2 &    5 &  \,   \,   \,114 &  \,   \,630 &$ 1.072964\times 10^{13} $\\
    1 &3 &   46 &  \,   \,  1\,644 &  \, 15\,015 &$ 4.593811\times 10^{24} $\\
    1 &4 &  669 &  \,   \, 24\,156 &  \,306\,306 &$ 6.326929\times 10^{37} $\\
    1 &5 &      &  \,   \,511\,416 & 5\,819\,814 &$ 1.132261\times 10^{52} $\\
    2 &1 &    9 &  \,   \,  6\,336 &  \, 15\,015 &$ 4.593811\times 10^{24} $\\
    2 &2 &  368 &  \,   \, 17\,982 &  \,306\,306 &$ 6.326929\times 10^{37} $\\
    2 &3 &      &  \,   \,606\,144 & 5\,819\,814 &$ 1.132261\times 10^{52} $\\
    3 &1 &      & 1\,065\,718\,368\makebox[0pt]{\space${}^\dagger$}
                                   & 5\,819\,814 &$ 1.132261\times 10^{52} $\\
    \bottomrule
  \end{tabular}    
  \caption{\label{tab:N}%
    Number of (non-unique) trivalent fatgraphs generated according to
    different algorithms. The $N$ column reports the actual number of distinct fatgraphs
    for the given $g$, $n$; empty cells mean the corresponding number has not been
    computed.  The $N_1^+$ column lists upper bounds for the recursive generation algorithm
    \q{MgnTrivalentGraphs} (see Listing~\ref{algo:MgnTrivalentGraphs}); values marked with 
    the ``$\dagger$'' symbol are estimated using earlier values of $N_1^+$ because 
    the corresponding values of $N$ are not available.
    The $N_2^+$ values bound from above the number of fatgraphs generated by grafting binary
    trees into clovers. 
    Finally, $N_3$ is the count of fatgraphs generated by enumerating
    pairs of permutations (as per combinatorial definition of fatgraph).
  }
\end{table}

In what follows, let $N(g,n) := \card{\R_{g,n}}$ be the number of
distinct $(g,n)$-fatgraphs; also define:
\newcommand\mmax{m_{\text{max}}}
\newcommand\mmin{m_{\text{min}}}
\begin{align}
  \xi(g,n)   &:= 2g +n,                     \label{eq:xi}
  \\
  \mmax(g,n) &:= 6g + 3n - 6 = 3\xi - 6,    \notag
  \\
  \mmin(g,n) &:= 2g + n - 1 = \xi -1.       \notag
\end{align}
It is trivial to check that $\mmax$ and $\mmin$ are the maximum and
minimum number of edges that a $(g,n)$-fatgraph can have.

\subsection{Generation by recursive edge addition}
\label{sec:N1}

The algorithm \q{MgnTrivalentGraphs} described in Section~\ref{sec:MgnTrivalentGraphs}
produces a $(g,n)$-fatgraph by adding an edge to fatgraphs with lower
$(g,n)$; the procedure can then be applied recursively.

Let $N_1(g,n)$ be the number of (non distinct) fatgraphs returned by
"MgnTrivalentGraphs($g$,$n$)".  According to
Section~\ref{sec:MgnTrivalentGraphs}, this can be written as:
\begin{equation*}
  N_1(g,n) = N_{1,A}(g,n) + N_{1,B}(g,n) + N_{1,C}(g,n),
\end{equation*}
where $N_{1,A}$, $N_{1,B}$, $N_{1,C}$ are the numbers of fatgraphs
constructed in cases {\sl A)}, {\sl B)}, {\sl C')} of
Algorithm~\ref{algo:MgnTrivalentGraphs}. 

In case {\sl A)}, we have 1 generated $(g,n)$-fatgraph per each pair
formed by a $(g,n-1)$-fatgraph and one of its \emph{oriented} edges,
modulo the action of the automorphism group $\Aut(G)$.  However, the
number of orbits of this $\Aut(G)$-action is difficult to estimate.
Since the generic fatgraph only has one automorphism, an upper bound
can instead be given by considering all pairs formed by a fatgraph and
an oriented edge:
\begin{equation*}
  N_{1,A}(g,n) \leq N_{1,A}^+(g,n) := 2 \cdot \mmax(g,n-1) \cdot N(g,n-1).
\end{equation*}

In case {\sl B)}, the algorithm generates one $(g,n)$-fatgraph per
each triplet formed by a $(g,n-1)$-fatgraph and two oriented edges,
not necessarily distinct (modulo the action of $\Aut G$); a similar
remark about the upper bound applies:
\begin{equation*}
  N_{1,B}(g,n) \leq N_{1,B}^+(g,n) := (2 \cdot \mmax(g,n-1))^2 \cdot N(g,n-1) 
\end{equation*}

In case {\sl C')}, the computation is exactly the same, except we
apply the $q$ construction to fatgraphs belonging in $\R_{g-1,n+1}$:
\begin{equation*}
  N_{1,C}(g,n) \leq N_{1,C}^+(g,n) 
  := 4 \cdot \mmax(g-1,n+1)^2 \cdot N(g-1,n+1).
\end{equation*}
Table~\ref{tab:N} shows the upper bound given by
\begin{equation*}
  N_1^+(g,n) := N_{1,A}^+(g,n) + N_{1,B}^+(g,n) + N_{1,C}^+(g,n).
\end{equation*}

According to Table~\ref{tab:N}, the "MgnTrivalentGraphs" algorithm
outperforms the alternative procedures when $2g+n < 7$, and apparently
generates a much larger set of fatgraphs when $2g + n > 7$.  However,
the values were obtained using $N_1^+(g,n)$ instead of $N(g,n)$ in
recursive computations when the actual value of $N(g,n)$ is not known;
therefore $N_1^+(g,n)$ might grossly overestimate the number of graphs
considered by "MgnTrivalentGraphs" for $2g + n > 6$.  Further
investigation is needed to ascertain whether this is due to the
algorithm of Section~\ref{sec:N2} being asymptotically faster, or to
the estimate for $N_1(g,n)$ being grossly imprecise when no data about
the real number of trivalent fatgraphs in the recursion step is known.
However, the author conjectures that this estimate holds:
\begin{equation}
  \label{eq:N1-conj}
  N_1(g,n) \leq O(\xi^3) \cdot N(g,n)
\end{equation}

\subsection{Generation by insertion of binary trees}
\label{sec:N2}

A different approach is the following:
\begin{itemize}
\item Generate all distinct $(g,n)$-fatgraphs with 1 vertex; each such
  fatgraph has $\mmin(g,n)$ edges, hence the vertex has valence $2
  \cdot \mmin(g,n)$.
\item Given any such fatgraph $G_0$, build a trivalent
  $(g,n)$-fatgraph $G$ by replacing the vertex with a full binary tree
  on $2\cdot\mmin(g,n)$ leaves.
\end{itemize}
Call a fatgraph with only one vertex a \emph{clover}.  Let $N_2'(g,n)$
be the number of distinct $(g,n)$-clovers; we can estimate it as
follows.

\begin{lemma}
  \label{lemma:2}
  The number of isomorphic clovers is equal to the number
  of orbits of the adjoint action of $(1 2 \ldots 2m)$ over the
  set of self-conjugate permutations $\{ \sigma_1 \in \Perm{2m} :
  \sigma_1^2 = \idmap \}$.
\end{lemma}
\begin{proof}
  Let $G_0 = (L; \sigma_0, \sigma_1, \sigma_2)$ be a $(g,n)$-fatgraph
  given in combinatorial form, where $L = \{1, \ldots, 2m\}$ and
  $\sigma_i \in \Perm{m}$.  If $G_0$ is a clover, then
  $\sigma_0$ is a permutation formed by just one cycle; without loss
  of generality we may assume $\sigma_0$ is the rotation $(1 2 \ldots
  2m)$.  Let $G_0' = (L; \sigma_0', \sigma_1', \sigma_2)$ be another
  $(g,n)$-clover: by the same reasoning we have
  $\sigma_0' = \sigma_0 = (1 2 \ldots m)$; if $f\colon G_0 \to G_0'$ is an
  isomorphism, then $f$ commutes with $\sigma_0$ hence $f =
  \sigma_0^j$ for some $j | 2m$.  Therefore, from
  $\sigma_1 \circ f = f \circ \sigma_1'$ we get $\sigma_1' =
  \sigma_0^{-j} \circ \sigma_1 \circ \sigma_0^j$.  This proves the
  claim.
\end{proof}
\begin{lemma}
  \label{lemma:4}
  Let $L$ be a finite set of $l = p \cdot q$ elements. The number of permutations
  of $L$ which can be expressed as product of $q$ disjoint $p$-cycles
  is:
  \begin{equation}
    \label{eq:1}
    C(p,q) = \prod_{i=1}^{q} \prod_{j=1}^{p-1} (pi - j).
  \end{equation}
\end{lemma}
\begin{proof}
  Without loss of generality we can assume $L = \{ 1, \ldots, pq \}$;
  let $\tau \in \Perm{pq}$ be a permutation composed of $q$ disjoint
  $p$-cycles.
  We can give a ``canonical'' form to $\tau$ if we order its cycles
  by stipulating that:
  \begin{itemize}
  \item a cycle $(a_1 a_2 \ldots a_p)$ is always written such that
    $a_1 = \min a_i$;
  \item $(a_1 a_2 \ldots a_p)$ precedes $(b_1 b_2 \ldots b_p)$ iff
    $\min a_i < \min b_i$.
  \end{itemize}

  Now assume $\tau$ is written in this canonical form; then $a_1 = 1$
  and we have $pq-1$ choices for the element $a_2 = \tau(a_1)$
  following $a_1$ in the cycle, $pq-2$ choice for the next element
  $a_3 = \tau(a_2)$, and so on until the final element $a_p$ of the
  first cycle.  Then starting element $a_{p+1}$ of the second cycle
  has to be the minimum element of $L \setminus \{ a_1, a_2, \ldots,
  a_p \}$, but we have $(p-1)q - 1$ choices for $a_{p+2} =
  \tau(a_{p+1})$: an iterative argument proves the assertion.
\end{proof}
\begin{lemma}
  \label{lemma:3}
  The number of distinct self-conjugate permutations on a set of $l$
  elements is given by $(l-1)!! := (l-1) \cdot (l-3) \cdot \ldots
  \cdot 1$.
\end{lemma}
\begin{proof}
  A self-conjugate permutation $\tau$ on a set $L$ of $l = 2m$
  elements is the product of $m$ disjoint 2-cycles, and the the result
  follows from Lemma~\ref{lemma:4}.
\end{proof}
Combining Lemma~\ref{lemma:2} and~\ref{lemma:3}, we immediately get
the following estimate:
\begin{equation*}
  \label{eq:N2'}
  \frac{(2m-1)!!}{2m} \leq N_2'(g,n) \leq (2m-1)!!,
  \qquad m = \mmin(g,n),
\end{equation*}
where the upper bound comes from assuming that no two clovers can
be transformed one into the other by a rotation, and the lower bound
comes from considering all clovers as part of the same equivalence class.

In order to create a trivalent fatgraph from a clover, we replace the
vertex with a full binary tree with $l = 2m$ leaves; equivalently, we
identify the leaves of the tree according to the same ``gluing
pattern'' that identifies half-edges in the clover.  

More precisely, let $G_0 = (L; \sigma_0, \sigma_1, \sigma_2)$ be a
clover, with $L = \{ 1, \ldots, 2m \}$ and $\sigma_0 = (1 2 \ldots
2m)$ as above.  Let $L'$ be set of leaves of a chosen binary tree $T$
and $f\colon L' \to L$ a bijection.  Now $\tau := f\inv \circ \sigma_1
\circ f$ is a fixed-point free involution on $L'$: by identifying
leaves of $T$ according to $\tau$, we get a trivalent fatgraph
$G$, which we say is obtained by plugging $T$ into $G_0$ (by means of
$f$).

Given a permutation $\phi'$ on $L'$, the map $f' = f \circ \phi'$ is
a bijection and we have:
\begin{equation*}
  \tau' = {f'}\inv \circ \sigma_1 \circ f' = {\phi'}\inv \circ
  (f\inv \circ \sigma_1 \circ f) \circ \phi' = {\phi'}\inv \circ \tau
  \circ {\phi'},
\end{equation*}
which is an involution on $L'$ conjugate to $\tau$.  Conversely, if
$\sigma_1' = \phi\inv \circ \sigma_1 \circ \phi$ is conjugate to
$\sigma_1$, then $f' = \phi \circ f\colon L' \to L$ is again a bijection,
hence:
\begin{equation*}
{f}\inv \circ \sigma_1' \circ f = (f\inv \circ \phi\inv)
\circ \sigma_1 \circ (\phi \circ f)
= {f'}\inv \circ \sigma_1 \circ {f'},
\end{equation*}
which is the involution defining the attachment map of $T$ to $G_0$ by
means of $f'$.  Since any two involutions are conjugate, we can fix
the map $f$ once and for all binary trees with the same number of
leaves, and only let the involution $\sigma_1$ (i.e., the clover
$G_0$) vary.

Therefore $N_2(g,n) = N_2'(g,n) \cdot Y({\mmin(g,n))}$, where $Y(l)$ is
the count of full binary trees with $l$ leaves. The number $Y(l)$ is
given by the $(l-1)$-th Catalan number:
\begin{equation*}
  Y(l) = \frac{(2l-2)!}{(l-1)! \cdot l!}.
\end{equation*}
Hence from \eqref{eq:N2'} we get:
\begin{equation*}
  N_2^-(g,n) \leq N_2(g,n) \leq N_2^+(g,n),
\end{equation*}
where:
\begin{gather}
  N_2^-(g,n) := \frac{1}{2m} \cdot \frac{(4m-2)!}{(2m-2)!! (2m)!}, 
  \notag
  \\
  \label{eq:N2+}
  N_2^+(g,n) := \frac{(4m-2)!}{(2m-2)!! (2m)!},
  \\
  m := \mmin(g,n). 
  \notag
\end{gather}

\subsection{Generation from permutations}
\label{sec:N3}

As in the previous section, represent a fatgraph $G$ in combinatorial
form as $(L; \sigma_0, \sigma_1, \sigma_2)$ where $L = \{ 1, \ldots,
2m \}$.  Here we count the number of trivalent fatgraphs that
are generated by naively constructing a fatgraph from its
combinatorial definition.

If $G$ is trivalent, then $\sigma_0$ is a product of disjoint
3-cycles; by Lemma~\ref{lemma:4}, the number of such $\sigma_0$ is:
\begin{equation}
  \label{eq:N3a}
  C(3,k) = (l-1)(l-2) \cdot (l-4)(l-5) \cdot \ldots \cdot 2 \cdot 1,
  \qquad l = 2m = 3k
\end{equation}
For each chosen $\sigma_0$, each choice of a self-conjugate
permutation $\sigma_1$ gives rise to a trivalent $(g,n)$-fatgraph; by
Lemma~\ref{lemma:3} there are exactly $(2m-1)!!$ such choices.
Therefore, we have:
\begin{equation}
  \label{eq:N3}
  \begin{split}
    N_3(g,n) = (2m-1)!! \cdot C(3, 2m/3) = (2m-1)!! \cdot (2m-1)(2m-2) \cdot 
    \\ \cdot (2m-4)(2m-5) \cdot \ldots \cdot 2 \cdot 1,
  \end{split}
\end{equation}
where $m = \mmax(g,n)$.

\section{Pseudo-code notation}
\label{sec:pseudo-code}

Blocks of code are marked by indentation (rather than delimited by
specific keywords).  

The `"def"'
keyword is used to mark the beginning of a function definition.

The notation `"for $x$ in $S$"' is used to loop over all the items $x$
in a set or sequence $S$; sometimes the notation `"for x in $a,
\ldots, b$"' is used instead.  The form `"for $i, x$ in enumerate($S$)"' 
is used for keeping track of the iteration number when looping over
the elements of $S$: as $x$ runs over the items in $S$, $i$ orderly
takes the values $0$, $1$, \ldots, up to $\card{S} - 1$.

\subsection{Basic types}

Numbers and basic data structures (arrays, lists, sets; see below) are
considered basic types, together with the logical constants "True" and
"False", and the special value "None".

\subsection{Objects}

The word ``object'' is used to denote a kind of aggregate type: an
\emph{object} is a tuple `"($a_1$, $a_2$, $\ldots$, $a_N$)"', where
each of the slots "$a_i$" can be independently assigned a value; the
values assigned to different $a_i$'s need not be of the same type.  We
write $X.a_i \gets b$ to mean that the slot $a_i$ of object $X$ is
assigned the value $b$. Unless otherwise noted, object slots are
mutable, i.e., they can be assigned different values over the course
of time.

An objects' \emph{class} is the tuple `"($a_1$, $\ldots$, $a_N$)"' of
slots names that defines the object; the actual tuple of values is
called an \emph{instance}.\footnote{%
  Readers familiar with object-oriented programming will note that
  this is an over-simplified version of the usual object-oriented
  definition of objects and classes; this originates in the fact that
  the concrete implementation of the algorithms was done in
  object-oriented Python, but object-orientation is by no means
  essential to the implementation.}


\subsection{Arrays, lists, sets}

A few types of basic data structures are used in the code: arrays,
lists and sets.  They are distinguished only for clarity, and we make
no assumption that these are primitive: for instance, each of these
data structures could be implemented on top of the ``list'' type
defined here.

An ``array'' is a fixed-size collection of elements of the same type;
the number and type of elements stored in an array will be stated when
the array is first created.  Items in an array can be accessed by
position: if $a$ is an array, then its $k$-th element will be accessed
as $a[k]$.  Array elements can be mutated; we write $a[k] \gets b$ to
mean that object $b$ is stored into the $k$-th place of array $a$.

A ``list'' is a variable-size collection of objects. Two features
distinguish lists from arrays:
\begin{inparaenum}
\item lists can grow and shrink in size, and
\item lists can store items of different types.
\end{inparaenum}
If $l$ is a list with $n$ elements, the notation "$l$.append($x$)"
will be used to mean that $x$ should be added as $(n+1)$-th item in list
$l$.  Again, the square bracket notation $l[k]$ is used to denote the
value stored in the $k$-th place in $l$, and $l[k] \gets x$ means that
the $k$-th slot of $l$ is mutated to the value $x$.  The operator
``$+$'' stands for concatenation when applied to lists.

A ``set'' is a mutable unordered collection of objects
of the same type.  The only relevant difference with sets in the
mathematical sense of the word is that set variables are mutable:
if $s$ is a set, then "$s$.add($x$)" will be used to specify that $s$
should be mutated into the set $s \cup \{x\}$.  No duplicates
are admitted: if $x \in s$ and $x = y$, then "$s$.add($y$)" does not
alter $s$ in any way.

The word ``sequence'' will be used to denote any one of the above
three.  When $S$ is a sequence, we define "size($S$)" as the number of
elements in $S$; if $S$ is a "list" or "array" object, valid indices
into $S$ range from $0$ to "size($S$)-1".

\subsubsection{List comprehensions}

A special syntax is used to form a list when its items can be gotten
by applying a function or operation to the elements of another
sequence.

The notation `"L = [$f(x)$ for $x$ in $S$]"' makes $L$ into the list
formed by evaluating function $f$ on each element in $S$, analogously
to the usual notation $\{f(x) : x \in S\}$ for sets.

As an extension, the expression 
`"L = [$f(x)$ for $x$ in $S$ if $P(x)$]"' 
makes $L$ into the list of values of $f$ over the set $S'$
of elements of $S$ for which the predicate $P(x)$ is true: 
$S' = \{ x : x \in S \wedge P(x) \}$.

\subsection{Operators}

The ``"
division: for integers $k$ and $n>0$, the expression "($k$ 
evaluates to the smallest non-negative residue of $k \mod n$.

The ``$+$'' operator normally denotes addition when applied to
numbers, and concatenation when applied to lists.

Any other operator keeps it usual mathematical meaning.

\bibliographystyle{plain}
\bibliography{bib/math.AG,bib/math.NA,bib/cs}

\end{document}


